\documentclass{amsart}
\usepackage{amssymb}
\usepackage{amsthm}
\usepackage{epsfig}
\usepackage{hyperref}

\usepackage{epstopdf, color}

\usepackage{xypic}
\newcommand{\cd}[1]{\begin{equation*}{\xymatrix{#1}}\end{equation*}}
\newcommand{\cdlabel}[2]{\begin{equation}\label{#1}{\xymatrix{#2}}\end{equation}}

\def\bar{\overline}

\def\co{\colon\thinspace}
\def\R{\mathbb{R}}
\def\Z{\mathbb{Z}}

\def\H{\mathbb{H}}

\def\N{\mathbb{N}}
\def\Q{\mathbb{Q}}
\def\P{\mathcal P}

\newcommand{\mc}[1]{\mathcal{{#1}}}

\def\pp{p}
\def\Hbound {\partial_{\mc{H}}X}

\newtheorem{theorem}{Theorem}[section]
\newtheorem{lemma}[theorem]{Lemma}
\newtheorem{corollary}[theorem]{Corollary}
\newtheorem{proposition}[theorem]{Proposition}
\newtheorem{question}[theorem]{Question}

\newtheorem{claim}[theorem]{Claim}

\theoremstyle{definition}
\newtheorem{remark}[theorem]{Remark}
\newtheorem{definition}[theorem]{Definition}

\newtheorem{assumption}[theorem]{Assumption}

\newtheorem{case}{Case}

\newtheorem{observation}[theorem]{Observation}
\newtheorem{torsionremark}[theorem]{Remark (about torsion)}
\newtheorem{notation}[theorem]{Notation}

\newcommand{\llangle}{\langle\negthinspace\langle}
\newcommand{\rrangle}{\rangle\negthinspace\rangle}
\newcommand{\CComp}{\widehat{C}}

\newcommand{\Xgraph}{\ensuremath{X}}

\newcommand{\Gmk}{\ensuremath{G/K}}

\newcommand{\twid}{\widetilde}

\newcommand{\extend}[1]{\check{{#1}}}
\newcommand{\straightedge}[1]{{#1}_{\mc{T}}}

\def\e {e}   
\def\geod {\gamma}  
\newcommand{\skel}{\mathrm{Skel}}
\newcommand{\skelmap}{\straightedge{\ddot{\extend{\phi}}}}
\newcommand{\Lk}{\mathrm{Lk}}

\def\version{17 July, 2007}

\begin{document}
\title{Dehn filling in relatively hyperbolic groups}

\author[D. Groves]{Daniel Groves}
\address{California Institute of Technology}
\email{groves@caltech.edu}

\author[J. F. Manning]{Jason Fox Manning}
\address{University at Buffalo, SUNY}
\email{j399m@buffalo.edu}

\date{\version}

\thanks{The first author was supported in part by NSF Grant DMS-0504251.  The second author was supported in part by an NSF Mathematical Sciences Postdoctoral Research Fellowship.  Both authors thank the NSF for their support.  Most of this work was done
while both authors were Taussky-Todd Fellows at Caltech.}

\begin{abstract}
We introduce a number of new tools for the study of relatively
hyperbolic groups.  First, given a relatively hyperbolic
group $G$, we construct a nice combinatorial
Gromov hyperbolic model space acted on properly by $G$, which reflects
the relative hyperbolicity of $G$ in many natural ways.  Second, we
construct two useful bicombings on this space.  The first of these,
{\em preferred paths}, is combinatorial in nature and allows us to
define the second, a relatively hyperbolic version of a construction
of Mineyev.

As an application, we prove a group-theoretic analog of the
Gromov-Thurston $2\pi$ Theorem in the context of relatively hyperbolic
groups. 
\end{abstract} 

\maketitle

\tableofcontents

\section{Introduction}
A finitely generated group is {\em word hyperbolic} \cite{gromov:hyp}
if it acts properly and cocompactly on a metric space (e.g., its
Cayley graph)
satisfying a certain coarse geometric property called Gromov
hyperbolicity (Definition \ref{d:hyperbolic} below).  Two spaces acted
on properly and cocompactly by the same 
group will have the same coarse geometry, so word hyperbolicity
depends only on the group in question. 
For example, the
fundamental group of a compact hyperbolic $n$-manifold acts properly
and cocompactly on the hyperbolic space $\H^n$, and so is word
hyperbolic.  
If $G$ is the fundamental group of a non-compact but finite volume
hyperbolic $n$-manifold, then $G$ also has an apparently natural
action on $\H^n$.  We wish to say that groups with such actions on
Gromov hyperbolic spaces are {\em relatively hyperbolic}.
Considering the case of fundamental groups of
finite area hyperbolic surfaces, we see that it is necessary first to
specify a collection of ``peripheral'' subgroups to determine an
action up to some kind of coarse equivalence.  Thus it 
makes no real sense to ask whether a group is {\em relatively
hyperbolic}, but only to ask whether a group is hyperbolic {\em
relative} to a collection of subgroups.
There are
several competing ways to say what it means for a group to be (strongly)
hyperbolic relative to a collection of subgroups
\cite{gromov:hyp,farb:relhyp,bowditch:relhyp}.
These
definitions are now all known to be equivalent. \footnote{In 
\cite{osin:relhypbook}, Osin
gives a more general definition where it is not assumed that parabolic
subgroups are finitely generated, or that there are finitely many conjugacy
classes of parabolic subgroups.  Since the current paper appeared
as a preprint, Chris Hruska \cite{Hruska} has extended the definition proposed here (Theorem \ref{t:tfae}.\eqref{tfae5})
to the infinitely generated setting, and proved its equivalence
to Osin's definition.} 
(See Section \ref{ss:relhyp} for precise definitions and 
more examples of relatively hyperbolic groups.)

This paper has three main purposes.
First, we introduce a new space (the ``cusped space'') for studying
relatively hyperbolic groups (Section \ref{section:cusped}).  Second, we
construct a pair of useful bicombings on this space (Sections
\ref{section:pp} and \ref{section:bicombing}).
Third, we extend Thurston's Hyperbolic Dehn Surgery Theorem to the
context of (torsion-free) relatively hyperbolic groups (Part 2).  
We discuss these now in turn.

\subsection{A new geometry for relatively hyperbolic groups}
Roughly speaking, a group $G$ is hyperbolic relative to a collection
$\mc{P}$ of subgroups if the geometry of $G$ is $\delta$-hyperbolic,
except for that part corresponding to the subgroups in $\mc{P}$.  
In Subsection \ref{ss:combhoro} we
define {\em combinatorial horoballs} which are a way of embedding
any graph into a $\delta$-hyperbolic graph.  We build the cusped space by 
gluing onto the Cayley graph of $G$ a collection of combinatorial
horoballs, one for each coset of each subgroup in $\mc{P}$.  The resulting
space is Gromov hyperbolic if and only if $G$ is hyperbolic relative to
$\mc{P}$.  In case the space is Gromov hyperbolic, the action of $G$
on the space satisfies the conditions given in Gromov's original
definition of relative hyperbolicity (Definition \ref{d:gromovrh}).
The cusped
space thus combines the most useful combinatorial and coarse
geometric aspects of relatively hyperbolic groups.  A closely related
construction appears in work of Cannon and Cooper \cite{CannonCooper} (see
also \cite{rebbechi} for another related construction). 
As in \cite{CannonCooper}, if $G$ is
the fundamental group of a finite volume hyperbolic $n$-manifold, and 
$\mc{P}$ consists of the fundamental groups of the cusps, the space
we build is quasi-isometric to $\mathbb H^n$, though we do not provide
a proof of this here.

Part of the novelty
of the cusped space compared to the one in \cite{CannonCooper} 
is that it is a graph (metrized so the length of each edge is 
$1$), so that the metric and combinatorial aspects harmonize with each
other more easily.
This allows us to more easily adapt a number of constructions and
results in word hyperbolic groups to the relative setting.
In particular, we consider 
combinatorial isoperimetric inequalities (Subsection \ref{s:combii}),
homological isoperimetric inequalities (Subsection \ref{s:hom}), and
a construction of Mineyev from \cite{min:str}.  Considering the different
types of isoperimetric inequalities and spaces, we get a number of new
characterizations of relatively hyperbolic groups.

Let $G$ be a finitely generated group which is finitely presented
relative to a collection $\mc{P} = \{ P_1 , \ldots , P_n \}$ of finitely
generated subgroups (see Subsection \ref{s:relpres}).  Let
$\hat\Gamma$ be the coned-off Cayley graph for $G$ with respect
to $\mc{P}$ (see Definition \ref{d:conedgraph}), $\hat{C}$ the
coned-off Cayley complex (see Definition \ref{d:conedcomplex}),
and let $X$ be the cusped space associated to $G$ and $\mc{P}$
(defined in Section \ref{section:cusped}).

Then we have

\medskip

\noindent{\bf Theorem \ref{t:tfae}}  
{\em The following are equivalent:
\begin{enumerate}
\item $G$ is hyperbolic relative to $\mc{P}$ in the sense of Gromov;
\item $G$ is hyperbolic relative to $\mc{P}$ (i.e. $\hat\Gamma$ is
Gromov hyperbolic and fine);
\item $\hat{C}$ satisfies a linear combinatorial isoperimetric inequality;
\item $\hat{C}$ satisfies a linear homological isoperimetric inequality;
\item $X^{(1)}$ is Gromov hyperbolic;
\item $X$ satisfies a linear combinatorial isoperimetric inequality; 
\item $X$ satisfies a linear homological isoperimetric inequality.
\end{enumerate}
}

\medskip

See Subsection \ref{s:combii} and Definition \ref{d:hlii} for definitions
of linear isoperimetric inequalities (both combinatorial and homological).

The equivalence of \eqref{tfae1} and \eqref{tfae2} is by now
well known (see, for example, \cite[Appendix]{Dah:thesis}).
As far as we are 
aware, the equivalence of \eqref{tfae2} and \eqref{tfae3} has
not appeared elsewhere, though it is implicit in \cite{osin:relhypbook}.
What is really novel about Theorem \ref{t:tfae} is the space $X$ and the use
of homological isoperimetric inequalities.  

\subsection{Bicombings on relatively hyperbolic groups}
The second main purpose of this paper is the construction
in Sections \ref{s:horoballs}--\ref{section:bicombing} of a pair of
useful bicombings on the cusped space.  In Section
\ref{s:horoballs} we prove a general result about convex sets
and between-ness in a Gromov hyperbolic space $\Upsilon$.  
Given a family
$\mc{G}$ of `sufficiently separated' convex sets we construct,
for any pair of points $a,b \in \Upsilon$, a
canonical collection of sets in $\mc{G}$ which are `between' $a$ and $b$.  
These collections satisfy a number of
axioms, and allow a great deal of combinatorial control over triangles
in $\Upsilon$ built from quasi-geodesics using our construction.

This analysis is undertaken in Section \ref{section:pp},
where we define {\em preferred paths} for our space $X$.  
These give a $G$-equivariant bicombing of $X$ 
by uniform quasi-geodesics, whose intersection
with horoballs is very controlled (this is where the results from Section
\ref{s:horoballs} are used).  In particular (possibly partially ideal)
triangles whose sides are preferred paths have very well controlled
combinatorial structure (see Subsection \ref{s:preftri} for details on this).  

We expect that the construction in Section \ref{s:horoballs}
and that of preferred paths in Section \ref{section:pp} will
have many applications.  The first is the bicombing $q$ which is 
defined in Section \ref{section:bicombing}.  This gives a relatively
hyperbolic version of a construction of Mineyev from \cite{min:str}.
Applications of Mineyev's construction are myriad (see, for example,
\cite{min:str,min:bcchg,min:flow,minmosha,minyu,yu}).  It can reasonably be expected that many of
these results can be extended to the relatively hyperbolic setting
using the bicombing from Section \ref{section:bicombing} of this paper,
or variations on it.  In particular, in \cite{GM:BoundCoh}, we define a
homological bicombing on the coned-off Cayley graph of a relatively
hyperbolic group (using the bicombing from this paper in an essential
way) in order to investigate
relative bounded cohomology and relatively hyperbolic groups, in analogy with \cite{min:str} and \cite{min:bcchg}.  Also, in \cite{G:Bolic},
the first author proves that if the parabolic subgroups of $G$ act `nicely'
on a strongly bolic metric space (as defined by Lafforgue \cite{Lafforgue}) then so does $G$.  Using the work of Lafforgue \cite{Lafforgue} and Dru\c{t}u and Sapir \cite{DS:RD}, this has implications
for the Baum-Connes conjecture for certain relatively hyperbolic groups.

It is also worth noting that in Part 2 of this paper, the major tool is 
preferred paths.  The only time when we need the homological
bicombing (which is the analogue of Mineyev's construction) is in
the proof of Theorem \ref{t:Zhlii}.  Otherwise, we use only the results
from Section \ref{section:pp}, which have no relation to Mineyev's
construction.

\subsection{Relatively hyperbolic Dehn filling}
Part 2 of this paper is devoted to another application of the constructions
of Part 1.  We investigate a group theoretic analog of Dehn filling,
which is the third and final purpose of this paper.

We first briefly remind the reader what is meant by ``Dehn filling''
in the context of $3$-manifolds.  Suppose that $M$ is a compact 
$3$-manifold,
with some component $T$ of $\partial M$ homeomorphic to a torus.  Let
$\alpha$ be an essential simple closed curve in $T$.  Let $W$ be a
solid torus, and $\mu$ a meridian for $W$ (a curve which bounds an
embedded disk in $W$), and let $\phi\co (\partial W,\mu)\to (T,\alpha)$ be a
homeomorphism of pairs.  The $3$-manifold
\[ M(\alpha) = M \cup_\phi W \] obtained by gluing using $\phi$ is
called the {\em Dehn filling of $M$ along $\alpha$}, and depends up to
homeomorphism only on the homotopy class of $\alpha$ in $T$.
Thurston's Hyperbolic Dehn Surgery Theorem \cite{thurston:79} says that if the
interior of $M$ admits a hyperbolic metric, then so does the interior
of $M(\alpha)$, for all but finitely many choices of $\alpha$.  The
number of curves to be excluded and the relationship between the
geometry of $M$ and that of $M(\alpha)$ can be made quite precise (see
for example \cite{hodgsonkerckhoff:ub}).

On the level of fundamental groups, $\pi_1(M(\alpha)) =
\pi_1(M)/\llangle a\rrangle$, where $a$ is an element of $\pi_1(M)$
whose conjugacy class is represented by $\alpha$.  One of the group
theoretic statements implied the Hyperbolic Dehn Surgery Theorem is
\begin{theorem} \label{HDST} \cite{thurston:79}
Let $G=\pi_1(M)$ where $M$ is a hyperbolic $3$-manifold with a single
torus cusp $C$, and let $P=\pi_1(C)\cong \Z\oplus\Z$ be the cusp
subgroup.  Then for all but finitely many $a\in P$,
$G/\llangle a \rrangle$ is infinite, non-elementary, and word
hyperbolic.
\end{theorem}
This group theoretic statement
finds its strongest quantitative formulation
in the ``$6$ Theorem'' independently due to Lackenby and Agol
\cite{lackenby:whds,agol:whds}.  

Theorem \ref{HDST} puts the conclusion of the Hyperbolic Dehn
Surgery Theorem in an algebraic context.  Koji Fujiwara asked whether
there was an algebraic analogue of this theorem which also puts
the hypotheses into an algebraic context.  We learned of this question from Danny Calegari.

As an answer to this question, we provide the following result (where
$|K_i|_{P_i}$ denotes the minimal length of a nontrivial element of
$K_i$ using the word metric on $P_i$ with respect to the generating
set $S \cap P_i$;  see Definition \ref{d:slopelength}):

\medskip

\noindent{\bf Theorem \ref{t:rhds}}  
{\em
Let $G$ be a torsion-free group, which is hyperbolic
relative to a collection $\mc{P} = \{ P_1, \ldots , P_n \}$
of finitely generated subgroups.  Suppose that $S$
is a generating set for $G$ so that for each $1 \le i \le n$ 
we have $P_i = \langle P_i \cap S \rangle$.

There exists a constant $B$ depending only on 
$(G, \mc{P})$ so that for any collection $\{ K_i \}_{i=1}^n$ 
of subgroups satisfying
\begin{itemize}
\item $K_i \unlhd P_i$; and
\item $| K_i |_{P_i} \ge B$,
\end{itemize}
then the following hold, where $K$ is the normal closure in $G$ of
$K_1\cup \cdots \cup K_n$. 
\begin{enumerate}
\item The map $P_i/K_i\xrightarrow{\iota_i}G/K$ given by $pK_i\mapsto pK$
  is injective for each $i$.
\item $G/K$ is hyperbolic relative to the collection
  $\mc{Q}=\{\iota_i(P_i/K_i)\mid 1\leq i\leq n\}$.
\end{enumerate}
}

\medskip

It is well known (see, for example, \cite[Theorem 3.8]{farb:relhyp})
that a group which is hyperbolic relative to a collection
of word hyperbolic subgroups is word hyperbolic.  Thus, an immediate corollary of Theorem \ref{t:rhds} is the following:

\begin{corollary} \label{c:hypcase}
Under the hypotheses of Theorem \ref{t:rhds}, if each of the 
$P_i/K_i$ are themselves word hyperbolic, then $G/K$ is 
word hyperbolic.
\end{corollary}

Together with Theorem
\ref{t:nonelementary} (non-elementariness), Corollary \ref{c:hypcase}
unifies a number of known results: 
\begin{enumerate}
\item
We have already remarked that Theorem \ref{t:rhds} generalizes
Theorem \ref{HDST}.  However, since we make no assumptions
about the subgroups $P_i$, Theorem \ref{t:rhds} also generalizes
the results of Lackenby from \cite{lackenby:handlebodies} about
filling $3$-manifolds some of whose 
boundary components are higher genus surfaces (although we do not
obtain such nice quantitative bounds as Lackenby).
\item 
Corollary \ref{c:hypcase} also generalizes some known results about
hyperbolic groups.  For example, modulo the extra torsion-free hypothesis, 
Theorem \ref{t:rhds} is a generalization of statements
(1)--(3) of \cite[Theorem 5.5.D, p.149]{gromov:hyp}.  Of interest
in this context is that we make no use whatsoever of small
cancellation techniques.  
\item 
Much of the group-theoretic content of many ``CAT$(-1)$'' or ``CAT$(0)$ with
isolated flats'' filling constructions on hyperbolic manifolds with
torus cusps is also contained in Theorem \ref{t:rhds}.  Examples of
this are in \cite{schroeder:cuspclosing, moshersageev}.  (See also
Remark \ref{r:catminusone}.)
\end{enumerate}

We now make a few more remarks about Theorem \ref{t:rhds}.

\begin{remark}
The ``short'' fillings really must be excluded in Theorem
\ref{t:rhds}, as can be seen, for instance, from the many examples of exceptional
fillings of hyperbolic $3$-manifolds. 
By considering fillings of the Hopf link, we can see that
it is also important that
{\em each} of the lengths $|K_i|_{P_i}$ is large.  

An even simpler example is given by $G$ equal to the free group
$\langle x,y\rangle$, and $\mc{P}= \{P_1=\langle x\rangle,P_2=\langle
y\rangle,P_3=\langle xy \rangle\}$.  Choosing $K_1=\langle
x^p\rangle$, $K_2 = \langle y^q\rangle$, and $K_3 = \langle
(xy)^r\rangle$, the quotient $G/K$ will be infinite and word
hyperbolic if and only if $\frac{1}{p}+\frac{1}{q}+\frac{1}{r}<1$.  This
occurs if all three of $p$, $q$, and $r$ are at least $4$, but of
course if $p=q=2$, then $r$ can be arbitrarily large while $G/K$
remains finite.
\end{remark}

\begin{remark}
Since we are, from the very beginning, working in the coarse world of $\delta$-hyperbolic
spaces, we have no hope of obtaining the fine control over
constants, as obtained in \cite{lackenby:whds, agol:whds,
hodgsonkerckhoff:ub}.  Therefore, we have made very little
attempt throughout this paper to make our constants optimal.

However, it is worth remarking that there are some delicate
interdependencies between some of the constants we use.
\end{remark}

\begin{remark}\label{r:osin}
Denis Osin has independently proved Theorem \ref{t:rhds}; see
\cite{osin:peripheral}.  In fact, Osin works in a somewhat more
general setting, in two respects.

First, Osin has a more general
notion of relative hyperbolicity, which allows infinitely many
(possibly infinitely generated) parabolics.  It can be shown (and this
is really implicit in Osin's proof) that the appropriate
statement for infinitely
many infinitely generated parabolics 
\emph{follows} from the statement
for finitely many finitely generated parabolics.\footnote{For more on
  this, see \cite{GM2}.}

Second, and more seriously, Osin makes no assumption of
torsion-freeness in \cite{osin:peripheral}.  We believe that our
methods apply (with appropriate modification) to the non-torsion-free
case, at least to prove the analog of Theorem \ref{t:rhds}, but at
such a cost in clarity and brevity that we have elected 
to deal only with the torsion-free case.  We have tried to make
explicit our use of torsion-freeness and how one might go about
avoiding it (see Remarks \ref{r:tf1}, \ref{r:tf2}, \ref{r:tf3}, \ref{r:tf4}, 
\ref{r:tf5}, \ref{tr:bonus},\ref{r:tf6} and \ref{r:tf7};  on first reading,
we recommend ignoring all of these comments; to facilitate this,
they are all labelled as `Remark (about torsion)').
Certain of our results, which are not used in the proof of Theorem
\ref{t:rhds}, must be considerably rephrased in the presence of
torsion (see especially Theorem \ref{t:dontintersect}).

Note also that Osin's main theorem \cite[Theorem 1.1]{osin:peripheral}
states that given a finite set $\mc{F} \subset G$, there is a $B$
so that under the conditions of Theorem \ref{t:rhds} the map
$G \to G/K$ is injective on $\mc{F}$.  This is Corollary 
\ref{c:injectonfinite} below.
\end{remark}

\begin{remark}\label{r:catminusone}
It is worth noting that, even when starting with a rank one locally symmetric
space, our group-theoretic version of filling produces hyperbolic
groups in many situations where the existence a locally CAT$(-1)$ filling is
not at all clear; see the non-existence results of \cite{hummelschroeder}.
The advantage of a CAT$(-1)$ filling is that information about the
fundamental group can be obtained from local information about the
locally CAT$(-1)$ model.  For more on this, see
\cite{fujiwaramanning}, which is in preparation.
\end{remark}

\section{Preliminaries} \label{s:prelim}

\subsection{Coarse geometry}
All metric spaces will be assumed to be complete geodesic metric
spaces, and the distance between two points $x$ and $y$ will usually
be denoted $d(x,y)$. 

\begin{definition}
If $X$ is a metric space, $A \subset X$ and $R \ge 0$, then let 
$N_R(A)$ be the $R$-neighborhood of $A$ in $X$.
\end{definition}

\begin{definition}
If $X$ and $Y$ are metric spaces, $K\geq 1$ and $C\geq 0$,
a \emph{$(K,C)$-quasi-isometric embedding} of $X$ into $Y$ is a function
$q\co X\to Y$ so that
 For all $x_1$, $x_2\in X$
\[\frac{1}{K}d(x_1,x_2)-C\leq d(q(x_1),q(x_2))\leq Kd(x_1,x_2)+C\]

If in addition the map $q$ is \emph{$C$-coarsely onto} -- \emph{i.e.},
$N_C(q(X))=Y$ -- $q$ is called a
\emph{$(K,C)$-quasi-isometry}.
The two metric spaces $X$ and $Y$ are then
said to be \emph{quasi-isometric} to one
another.  This is a symmetric condition.
\end{definition}
\begin{definition}
A \emph{$(K,C)$-quasi-geodesic} in $X$
is a $(K,C)$-quasi-isometric embedding $\gamma\co\R\to X$.
We will occasionally abuse
notation by referring to the image of $\gamma$ as a quasi-geodesic.

A {\em $(K,C)$-quasi-geodesic ray} is a $(K,C)$-quasi-isometric
embedding $p \co \R_{\ge 0} \to X$.
\end{definition}

\subsection{Gromov hyperbolic spaces}
\begin{figure}[htbp]
\begin{center}
\begin{picture}(0,0)%
\includegraphics{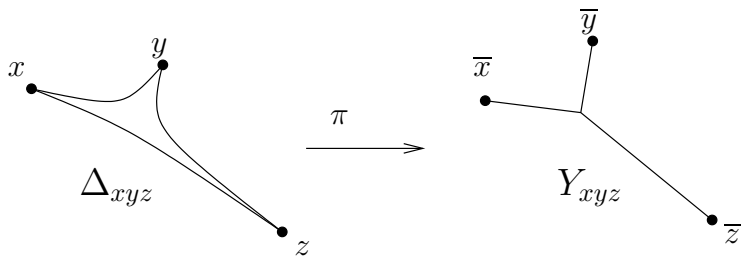}%
\end{picture}%
\setlength{\unitlength}{3947sp}%
\begingroup\makeatletter\ifx\SetFigFont\undefined%
\gdef\SetFigFont#1#2#3#4#5{%
  \reset@font\fontsize{#1}{#2pt}%
  \fontfamily{#3}\fontseries{#4}\fontshape{#5}%
  \selectfont}%
\fi\endgroup%
\begin{picture}(4530,1665)(436,-1930)
\put(451,-736){\makebox(0,0)[lb]{\smash{{\SetFigFont{12}{14.4}{\familydefault}{\mddefault}{\updefault}{\color[rgb]{0,0,0}$x$}%
}}}}
\put(1351,-586){\makebox(0,0)[lb]{\smash{{\SetFigFont{12}{14.4}{\familydefault}{\mddefault}{\updefault}{\color[rgb]{0,0,0}$y$}%
}}}}
\put(2251,-1861){\makebox(0,0)[lb]{\smash{{\SetFigFont{12}{14.4}{\familydefault}{\mddefault}{\updefault}{\color[rgb]{0,0,0}$z$}%
}}}}
\put(3376,-736){\makebox(0,0)[lb]{\smash{{\SetFigFont{12}{14.4}{\familydefault}{\mddefault}{\updefault}{\color[rgb]{0,0,0}$\overline{x}$}%
}}}}
\put(4051,-436){\makebox(0,0)[lb]{\smash{{\SetFigFont{12}{14.4}{\familydefault}{\mddefault}{\updefault}{\color[rgb]{0,0,0}$\overline{y}$}%
}}}}
\put(4951,-1786){\makebox(0,0)[lb]{\smash{{\SetFigFont{12}{14.4}{\familydefault}{\mddefault}{\updefault}{\color[rgb]{0,0,0}$\overline{z}$}%
}}}}
\put(2476,-1036){\makebox(0,0)[lb]{\smash{{\SetFigFont{12}{14.4}{\familydefault}{\mddefault}{\updefault}{\color[rgb]{0,0,0}$\pi$}%
}}}}
\put(901,-1486){\makebox(0,0)[lb]{\smash{{\SetFigFont{14}{16.8}{\familydefault}{\mddefault}{\updefault}{\color[rgb]{0,0,0}$\Delta_{xyz}$}%
}}}}
\put(3901,-1486){\makebox(0,0)[lb]{\smash{{\SetFigFont{14}{16.8}{\familydefault}{\mddefault}{\updefault}{\color[rgb]{0,0,0}$Y_{xyz}$}%
}}}}
\end{picture}%
\caption{A triangle and its comparison tripod}
\label{fig:tripdef}
\end{center}
\end{figure}

Given a geodesic
triangle $\Delta(x,y,z)$ in
any geodesic metric space, there is a unique \emph{comparison tripod},
$Y_{xyz}$, a metric tree so
that the distances between the three extremal points of the tree,
$\overline{x}$, $\overline{y}$ and $\overline{z}$ , are the same as
the distances between $x$, $y$ and $z$ (See Figure \ref{fig:tripdef}.).
There is a unique map 
\[\pi\co \Delta(x,y,z)\to Y_{xyz}\] 
which takes $x$ to
$\overline{x}$, $y$ to $\overline{y}$ and $z$ to $\overline{z}$, and
which restricts to an isometric embedding on
each side of $\Delta(x,y,z)$.

\begin{definition}
Let $\delta\geq 0$.
The triangle $\Delta(x,y,z)$ is \emph{$\delta$-thin} if the diameter of
$\pi^{-1}(p)$ is at most $\delta$ for every point $p\in Y_{xyz}$.
\end{definition}

\begin{remark}
Any triangle which is $\delta$-thin is also \emph{$\delta$-slim}:
i.e., any side of the triangle 
is contained in the $\delta$-neighborhood of the union of the other
sides.  More generally, a geodesic $n$-gon is $(n-2)\delta$-slim.
It is also true that for every $\delta$
there is a $\delta'$ so that $\delta$-slim triangles are
$\delta'$-thin (see \cite[Proposition III.H.1.17,
  p. 408]{bridhaef:book}).  
\end{remark}

\begin{definition} \label{d:hyperbolic}
A geodesic metric space $X$ is \emph{$\delta$-hyperbolic} if every
geodesic triangle in $X$ is $\delta$-thin.
If $\delta$ is
unimportant we may simply say that $X$ is \emph{Gromov hyperbolic}.
\end{definition}

See \cite[Chapter III.H]{bridhaef:book} for the background and many
basic results about Gromov hyperbolic spaces.

\begin{definition}
Let $x$, $y$, $z\in X$.  The \emph{Gromov product} of $x$ and $y$ with 
respect to $z$
is $(x,y)_z=\frac{1}{2}(d(x,z)+d(y,z)-d(x,y))$.  Equivalently, $(x,y)_z$ is the
distance from $\overline{z}$ to the central vertex of the comparison tripod
$Y_{xyz}$ for any geodesic triangle $\Delta(x,y,z)$.
\end{definition}

\begin{definition}
Fix some $z\in X$, where $X$ is some Gromov hyperbolic metric space. 
We say that a sequence $\{x_i\}$
\emph{tends to infinity} if $\liminf_{i,j\to\infty}(x_i,x_j)_z =
\infty$.  On the set of such sequences we may define an equivalence
relation: $\{x_i\}\sim\{y_i\}$ if $\liminf_{i,j\to\infty}(x_i,y_j)_z =
\infty$.  
The \emph{Gromov boundary} of $X$, also written $\partial X$,
is  the
set of equivalence classes of sequences tending to infinity.  The
Gromov boundary does not depend on the choice of $z$ (see \cite[Proposition III.H.3.7]{bridhaef:book}).  
\end{definition}
\begin{remark}
We may topologize $X\cup\partial X$ so that if $\{x_i\}$ tends to
infinity, then $\lim_{i\to\infty} x_i =[\{x_i\}]$.
Furthermore, if $\gamma\co{[0,\infty)}\to X$ is a quasi-geodesic ray,
then for any sequence $\{t_i\}$ with $\lim_{i\to\infty}t_i=\infty$,
the sequence $\{\gamma(t_i)\}$ tends to infinity.  
The point $\{\gamma(t_i)\}\in\partial X$ does not depend on
the choice of $\{t_i\}$.  An isometric action on $X$ extends to a
topological action on the boundary.
\end{remark}

\begin{remark} \label{r:deltaprops}
We will implicitly assume, whenever we say that a space
is $\delta$-hyperbolic, that triangles are $\delta$-thin and
$\delta$-slim. 

This can be achieved by replacing $\delta$ by some larger 
constant.

We will also assume that $\delta$ is an integer.
\end{remark}

One can also consider geodesic ``triangles'' in a Gromov hyperbolic
space $X$
whose vertices are {\em  ideal}, i.e., points in $\partial X$.
The following is a simple exercise in $\delta$-slim triangles.
\begin{lemma} \label{l:idealthin}
A geodesic triangle in a $\delta$-hyperbolic space
with some or all vertices ideal is 
$3\delta$-slim.
\end{lemma}

\subsection{Cayley graph}
Although there is little difference between various definitions of the
Cayley graph, it is convenient to fix one here.  By ``$G$ is generated
by $S$'', we mean that there is a surjection
\[\pi\co \mc{F}(S)\to G,\]
where $\mc{F}(S)$ is the free group on the set $S$. 
If $S\subset G$, then we implicitly assume $\pi$ is the homomorphism
induced by inclusion.
The {\em Cayley graph of $G$ with respect to $S$}, written
$\Gamma(G,S)$, is the graph with vertex set $G$, and edge set $G\times
S$.  The edge $(g,s)$ connects the vertices $g$ and $g\pi(s)$.

\subsection{Relative presentations} \label{s:relpres}

We recall the following definitions of Osin. (We change the notation
slightly.)

\begin{definition} \label{d:relgen}
\cite[Definition 2.1]{osin:relhypbook}
Let $G$ be a group, $\{ H_\lambda \}_{\lambda \in \Lambda}$ a 
collection of subgroups of $G$, and $A$ a subset of $G$.  We say that
$A$ is {\em a relative generating set for $G$ with respect to 
$\{ H_\lambda \}_{\lambda \in \Lambda}$} if $G$ is generated by
\[ A \cup \left( \cup_{\lambda \in \Lambda} H_\lambda \right) .	\]
\end{definition}

We will be concerned with situations where the index set $\Lambda$
is finite.

\begin{definition} \label{d:relpres}
[Osin]
Suppose that $G$ is generated by $\mc{A}$ with respect to
$\{ H_\lambda \}_{\lambda \in \Lambda}$.  Then $G$ is a quotient
of 
\[ F = 	F(\mc{A}) \ast \left( \ast_{\lambda \in \Lambda} H_\lambda \right),	\]
where $F(\mc{A})$ is the free group on the alphabet $\mc{A}$.
Suppose that $N$ is the kernel of the canonical quotient map
from $F$ to $G$.  If $N$ is the normal closure of the set $\mc{R}$ then
we say that
\[	\langle \mc{A}, \{ H_\lambda \}_{\lambda \in \Lambda} \mid
\mc{R} \rangle	,	\]
is a {\em relative presentation} for $G$ with respect to 
$\{ H_\lambda \}_{\lambda \in \Lambda}$.

We say that $G$ is {\em finitely presented relative to
$\{ H_\lambda \}_{\lambda \in \Lambda}$} if we can choose
$\mc{R}$ to be finite.
\end{definition}

The following lemma is essentially contained in Theorem 2.44 of \cite{osin:relhypbook}.
\begin{lemma}\label{l:Hfg}
If $G$ is a finitely generated group, finitely presented relative to a collection of nontrivial 
subgroups $\{H_\lambda\}_{\lambda\in\Lambda}$, then $\Lambda$ is finite and each $H_\lambda$ is finitely generated.
\end{lemma}
\begin{proof}
By Theorem 2.44 of \cite{osin:relhypbook}, there is a finite subset $\Lambda_0\subseteq \Lambda$ so that 
\[G\cong (\ast_{\lambda\in\Lambda\smallsetminus\Lambda_0}H_\lambda)\ast G_0,\]
where $G_0$ is generated by $\bigcup_{\lambda\in\Lambda_0}H_\lambda$ and the finite relative generating set for $G$.

By Grushko's Theorem, and the fact that $G$ is finitely generated, $\Lambda\smallsetminus \Lambda_0$ is also finite.  Moreover, each $H_\lambda$ for $\lambda\in \Lambda\smallsetminus\Lambda_0$ is finitely generated, as is $G_0$. 

Osin's theorem further asserts that $G_0$ has the structure of a tree of groups, where each edge group is finitely generated.  The collection of vertex groups is $\{H_\lambda\}_{\lambda\in\Lambda_0}\cup\{Q\}$ where $Q$ is some finitely generated group.  Since $G_0$ is finitely generated, an elementary application of Bass-Serre theory shows that each $H_\lambda$ is finitely generated.
\end{proof}

\begin{definition} \label{d:compatible}
Suppose that 
$G$ is finitely generated, and also finitely presented relative
to $\{ H_1, \ldots , H_m \}$.  By Lemma \ref{l:Hfg}, each of the $H_i$ is
finitely generated.  In this situation,
we  usually  fix a finite generating set $S$ for $G$ so that 
for each $1 \le i \le m$ we have $\langle S \cap H_i \rangle = H_i$.
We call such a set $S$ a {\em compatible} generating set for $G$.
\end{definition}

\begin{definition} \label{d:relCayleycomplex}
[Relative Cayley complex]
Suppose that $G$ is finitely presented relative to
the finitely generated subgroups
$\{ H_1 , \ldots , H_m \}$, and that $S$ is a finite
compatible generating set for $G$.  Let 
\[	\langle \mc{A}, \{ H_1 \ldots , H_m \} \mid
\mc{R} \rangle	,	\]
be a finite relative presentation for $G$, where
$\mc{A} \subset S$ is a set of relative generators for
$G$.

Let $\Gamma = \Gamma(G,S)$ be the Cayley graph
of $G$ with respect to $S$.  The elements of $\mc{R}$
correspond to loops in $\Gamma$.  We glue a $2$-cell
to each such loop, in a manner equivariant under the 
$G$-action on $\Gamma$.  The resulting $2$-complex
is called the {\em relative Cayley complex} of $G$ with
respect to $\langle S, \{ H_1 \ldots , H_m \} \mid
\mc{R} \rangle$, and is denoted $C(G,S,\mc{R})$.
\end{definition}

In general, the relative Cayley complex of $G$ will not
be simply-connected.  This will only be the case if 
$\langle S \mid \mc{R} \rangle$ is already a presentation for 
$G$.

In Subsection \ref{ss:relhyp} below, we define a simply-connected
$2$-complex $\hat{C}$ associated to a finite relative presentation
of relatively hyperbolic group, starting
with the relative Cayley complex.

Also, in Section \ref{section:cusped}, we construct another
another simply-connected $2$-complex $X$, again starting
with a finite relative presentation of a relatively hyperbolic group.

The spaces $X$ and $\hat{C}$ both contain copies
of the relative Cayley complex.

\subsection{Combinatorial maps, chains, etc.}

\begin{definition}
Let $Y$ be a cell complex, and let $\omega = \sum_i \alpha_i \omega_i$
be a real cellular $n$-chain.  The {\em $1$-norm} of $\omega$ is
\[	|\omega|_1 = \sum_i |\alpha_i|.	\]
\end{definition}

\begin{definition}
Suppose that $Y$ is a $2$-complex and that $p \co I \to Y^{(1)}$
is a combinatorial path, in the sense that there is a cell structure
on $I$ so that $p$ sends each edge either to an edge or a vertex
of $Y$.

The {\em length}, or {\em $1$-norm} of $p$ is then the number of
$1$-cells in $I$ which are mapped onto edges of $Y$.  We denote
the $1$-norm of $p$ by $|p|_1$.

The map $p$ induces an obvious cellular $1$-chain $\underline{p}$
on $Y$ (using the orientation on $I$).  
\end{definition}

\begin{definition}
Suppose that $Y$ is a $2$-complex and $\Sigma$ an oriented,
cellulated surface.  A 
{\em combinatorial map} is a map $f \co \Sigma \to Y$
which sends each vertex
of $\Sigma$ to a vertex of $Y$, each edge of $\Sigma$ to an edge
or a vertex of $Y$ and each $2$-cell to a $2$-cell, an edge, or
a vertex of $Y$.
Furthermore, if $\sigma$ is a cell in $\Sigma$ which is sent to a cell
of the same dimension then the interior of $\sigma$ is
mapped homeomorphically by $f$ onto its image.

The {\em area} of a combinatorial map $f \co \Sigma \to Y$ is
the number of $2$-cells in $\Sigma$ which are mapped onto
$2$-cells.  When we refer to the `$1$-norm' of a map between
$2$-complexes, we mean the area.

As in the $1$-dimensional setting, a combinatorial map
$f \co \Sigma \to Y$ induces an integral, cellular 
$2$-chain $\underline{f}$ on $Y$.
\end{definition}

\begin{remark}
If $f$ is a combinatorial map of an interval or a surface into a $2$-complex $Y$ then $|\underline{f}|_1 \le |f|_1$.
\end{remark}

\begin{remark}
Throughout this paper we are somewhat cavalier about the 
difference between paths as maps, paths as subsets, and 
$1$-chains.  This should not cause any confusion.
\end{remark}

\subsection{Combinatorial isoperimetric inequalities}
\label{s:combii}

Any Gromov hyperbolic space satisfies a linear coarse
isoperimetric inequality (see 
\cite[Proposition III.H.2.7]{bridhaef:book}).
However, when we work with homological isoperimetric inequalities 
below, we need to use simply connected spaces.

Thus we pause in this paragraph to consider combinatorial
isoperimetric inequalities.  There is little novel here, but it is
worth noting that our spaces are not always locally finite, are
not uniformly locally finite even when they are locally finite,
and there will rarely be a proper and cocompact action on the
Gromov hyperbolic spaces in this paper.  Thus we need to be
slightly careful about the hypotheses in the results below.

\begin{proposition} \label{p:hyptocomb}
Let $X$ be a simply-connected $2$-complex and suppose 
that $X^{(1)}$
is $\delta$-hyperbolic.  Suppose further that for some
$K \ge 0$, any combinatorial loop of length at most $16\delta$ 
can be filled with a combinatorial disk of area at most $K$. Then 
any combinatorial loop $c$ in $X$ can be filled with a combinatorial 
disk whose area is at most $K|c|_1$.
\end{proposition}
\begin{proof}
This is essentially Dehn's algorithm.  See the proof of
\cite[Proposition III.H.2.7]{bridhaef:book}.
\end{proof}

\begin{proposition} \label{p:combtohyp}
Suppose that $X$ is a simply connected $2$-complex
and that there is a constant $M \ge 0$ so that the length
of the attaching map of any $2$-cell is at most $M$.

If $X$ satisfies a linear combinatorial isoperimetric inequality
then $X^{(1)}$ is $\delta$-hyperbolic for some $\delta$.

Furthermore, $\delta$ can be calculated in terms of $M$ and
the isoperimetric constant of $X$.
\end{proposition}
\begin{proof}
Follows from the proof of  \cite[Theorem III.H.2.9]{bridhaef:book}, or from a combinatorial version of the proof of
Theorem \ref{t:homtohyp} below.
\end{proof}

\subsection{Homological things} \label{s:hom}

In \cite{gersten}, Gersten proves that a group is hyperbolic
if and only if it has a linear {\em homological} filling function.
In this section we recall those notions and extend them slightly
in order to account for the fact that the actions on Gromov hyperbolic
spaces in this paper are rarely proper {\em and} cocompact (though
they are usually one or the other).  In this section we prefer
to work with rational coefficients.  However, real or complex 
coefficients would also work.

We recall some standard definitions:
\begin{definition}
Let $Y$ be a cell complex, $A$ a normed abelian group (with norm $||\cdot||$), and let $\mc{C}_n$ be the set
of $n$-cells of $Y$.
A {\em locally finite $n$-chain in $Y$ with coefficients in $A$} $\omega$
is a formal sum
\[ \omega = \sum_{c\in Y} \omega_c c,\]
where each $\omega_c\in A$.  The chain $\omega$ is said to be {\em summable}
if 
\[ |\omega|_1:=\sum_{c\in Y}||a_c||<\infty.\]
The quantity $|\omega|_1$ is the {\em norm of $\omega$}.  The {\em
support} of $\omega$ is the union in $Y$ of the $n$-cells $c$ for which
$\omega_c$ is nonzero.

The chains of compact support form the standard cellular chain group
$C_n(Y,\Z)$.
\end{definition}

We first recall a result from \cite{min:bcchg} about 
expressing $1$-chains as sums of paths.  Note that
this result is a generalization of \cite[Theorem 3.3]{allcockgersten}.

\begin{definition}
Suppose that $\Gamma$ is a graph and $T$ is a collection
of vertices.  A {\em $T$-path} is a directed path in $\Gamma$ 
whose initial and terminal vertices are in $T$ (or are equal).

A path is {\em simple} if it has no repeated vertices.
\end{definition}

Let $f$ be a summable $1$-chain in a graph
$\Gamma$.  Let $\Gamma(f)$ denote the directed graph
which is $\Gamma$ with an orientation chosen so that
$f(e) \ge 0$ for each edge $e$.  Let $\Gamma_+(f)$ be
the minimal subgraph of $\Gamma(f)$ containing all the edges $e$
with $f(e) \neq 0$.

\begin{theorem} \label{t:coherent} 
[Mineyev, Theorem 6, \cite{min:bcchg}]
Let $\Gamma$ be a graph, $T$ a set of vertices in $\Gamma$
and $f$ a summable $1$-chain on $\Gamma$ with coefficients
in $\Q$ and $\text{supp}(\partial f)\subset T$.
\begin{enumerate}
\item[(a)] There is a countable family $P = \{ p_1,p_2, \ldots \}$
of simple $T$-paths in $\Gamma_+(f)$ and a sequence
$\{ \alpha_i \}$ in $\Q \cap [0,\infty)$ so that (i) $f = \sum_i 
\alpha_i p_i$; and (ii) $|f|_1 = \sum_i \alpha_i |p_i|_1$.
\item[(b)] If $f$ has finite support then $P$ can be chosen to be 
finite.
\end{enumerate}
\end{theorem}

\begin{definition}
Suppose that $f = \sum_i \alpha_i p_i$ expresses the $1$-chain
$f$ as a sum of $1$-chains $\alpha_i p_i$, where $\alpha_i \in 
\Q_{\ge 0}$.
This sum is called {\em coherent} if $|f|_1 = \sum_i \alpha_i |p_i|_1$.
\end{definition}

\begin{definition}\label{d:hlii}
Suppose that $Z$ is a simply-connected $2$-complex.
We say that $Z$ satisfies a {\em linear homological 
isoperimetric inequality} if there is a constant $K \ge 0$
so that for any combinatorial loop $c$ in $Z$ there is some $\sigma 
\in C_2(Z;\Q)$ with $\partial \sigma = \underline{c}$, satisfying 
\[ |\sigma|_1 \le K |c|_1.	\]
\end{definition}

See \cite[Theorem 7]{min:bcchg} for (many) other notions of what
it might mean for a space to have a `linear isoperimetric inequality'.  In
this paper, we exclusively use the notion from Definition \ref{d:hlii}
above.

The next result follows immediately from the definitions.

\begin{lemma} \label{l:combtohom}
If a simply-connected $2$-complex $Z$ satisfies a 
linear combinatorial isoperimetric inequality then it
satisfies a linear homological isoperimetric inequality.
\end{lemma}

In \cite{gersten}, Gersten proves that the Cayley complex
of a finitely presented group $G$ satisfies a linear
homological isoperimetric inequality if and only if $G$ is
hyperbolic.  This is essentially a converse to Lemma
\ref{l:combtohom} above.  In this section we slightly
generalize Gersten's result.  However, most of our
work is in noting that the proof of Gersten's result from
\cite{min:bcchg} works in our setting.

\begin{theorem} \label{t:homtohyp}
Suppose that $Z$ is a simply-connected $2$-complex
and that there is a constant $M$ so that the attaching map
for each $2$-cell in $Z$ has length at most $M$.

Suppose further that $Z$ satisfies a linear homological
isoperimetric inequality. Then $Z^{(1)}$, the $1$-skeleton
of $Z$, is $\delta$-hyperbolic for some $\delta$.
\end{theorem}
\begin{proof}
The proof consists of noting that a number of other proofs
do not rely essentially on finite valence.

The first step is to prove that if $Z$ is not Gromov hyperbolic
then for any $\epsilon > 0$, there are $\epsilon$-thick
geodesic bigons in $Z$.  This is essentially 
\cite[Theorem 1.4]{pap}, modified in the obvious way.  Namely,
for each $M$, define
let $f(r) = \inf \{ d(\gamma(R +r)), \gamma'(R+r)) \}$ where
the infimum is taken over all positive integers $R$ and all
$\gamma$, $\gamma'$, geodesics such that $\gamma(0) = 
\gamma'(0)$ and $d(\gamma(R),\gamma'(R)) \ge 2M^2$.
This is not quite the function that Papasoglu uses, but it suffices
for the proof.  The remainder of the proof of this first step
is identical to that in \cite{pap}.

Now follow the proof of \cite[Proposition 8]{min:bcchg} to
prove that if there are $\epsilon$-thick bigons in $Z$ 
for all $\epsilon > 0$
then $Z$ does not have a linear homological isoperimetric
inequality.  Mineyev's argument only relies on
a bound on the length of attaching maps for $2$-cells 
(and not, for instance, on vertex transitivity or local finiteness).  It is
certainly worth remarking that Mineyev relies on a result of Gersten
from \cite{gersten}, which similarly only needs a bound on the
attaching maps of $2$-cells.
\end{proof}

\subsection{Homological bicombings}\label{s:hombicomb}
There are various notions of ``bicombing'' for graphs.  In particular,
one can define combings made up of paths or of $1$-chains.  

\begin{definition}\label{d:geodesicbicombing}
Let $\Gamma$ be a graph, and let $\mathrm{Geod}(\Gamma)$ be the set of
(oriented) geodesics in $\Gamma$.  That is, each element of
$\mathrm{Geod}(\Gamma)$ is a path $\sigma\co I\to \Gamma$,
where $I\subset \R$ and $\sigma$ is parametrized by arc length.
A {\em geodesic bicombing on $\Gamma$} is a
function 
\[
\gamma\co \Gamma^{(0)}\times\Gamma^{(0)}\to \mathrm{Geod}(\Gamma)\]
so that $\gamma(x,y)$ is a geodesic which begins at $x$ and ends at
$y$.
\end{definition}

\begin{definition}\label{d:geodesicextension}
Let $\Gamma$ and $\mathrm{Geod}(\Gamma)$ be as in Definition
\ref{d:geodesicbicombing}. 
Let $\bar{\Gamma}$ be some compactification (or bordification) of
$\Gamma^{(0)}$, and let $\Lambda\subseteq \Gamma$.  A {\em geodesic
  bicombing on $\Lambda$} is a function 
\[\gamma\co \Lambda\times\Lambda\smallsetminus \Delta \to
  \mathrm{Geod}(\Gamma),\]
where $\Delta = \{(x,x)\mid x\in \Lambda\}$ and 
\begin{enumerate}
\item If $x\in \Gamma^{(0)}$, then $\gamma(x,y)$ starts at $x$;
  otherwise, $\lim_{t\to -\infty}\gamma(x,y)(t)=x\in \Lambda\smallsetminus\Gamma$.
\item If $y\in \Gamma^{(0)}$, then $\gamma(x,y)$ ends at $y$; otherwise
  $\lim_{t\to\infty}\gamma(x,y)(t)=y\in \Lambda\smallsetminus\Gamma$.
\end{enumerate}
\end{definition}

\begin{definition}\label{d:homologicalbicombing}  \cite{min:str}
Suppose that $\Gamma$ is a graph, and $\mathbb A$ is a ring.  Let
$C_1(\Gamma;\mathbb A)$ be the group of finite formal sums of
$1$-cells in $\Gamma$ with coefficients in $\mathbb A$. 
A {\em homological
bicombing on $\Gamma$} is a function 
\[q\co \Gamma^{(0)}\times\Gamma^{(0)}\to C_1(\Gamma;\mathbb A)\]
so that $\partial q(a,b) = b-a$.
\end{definition}
\begin{remark} \label{r:geodtohom}
It is clear that a geodesic bicombing as in Definition
\ref{d:geodesicbicombing} gives rise to a homological bicombing as in
Definition \ref{d:homologicalbicombing}.  It is slightly less obvious
(but also true) that we can use a bicombing as in Definition
\ref{d:geodesicextension} to produce something homological.
\end{remark}

\begin{definition}\label{d:homologicalextension}
Let $\Gamma$ be a graph and $\mathbb A$ a ring, and let $\bar{C}_1(\Gamma;\mathbb A)$
be the group of (locally finite) formal sums of $1$-cells in $\Gamma$.
Let $\bar{\Gamma}$ be some compactification (or bordification) of
$\Gamma^{(0)}$, and let $\Lambda\subseteq \Gamma$.  Let $q$ be a 
function 
\[q \co \Lambda\times\Lambda \to
  \bar{C}_1(\Gamma;\mathbb A),\]
which is zero precisely on $\Delta = \{(x,x)\mid x\in \Lambda\}$. 
Given $R\in \N$, $x,y\in \Lambda$, and $z\in \Gamma^{(0)}$,
let $q_{z,R}(x,y)$ be the $1$-chain which is equal to $q(x,y)$ on
the ball of radius $R$ about $z$, and zero outside it.  The function
$q$ is a {\em homological bicombing on $\Lambda$} if it satisfies the
following condition:
For every
$x,y\in \Lambda$, and $z\in \Gamma^{(0)}$, there is an $R_0$
so that for every integer $R>R_0$, there exist $0$-chains $\xi_{R,+}$ and
$\xi_{R,-}$, each with coefficients summing to $1$ so that:
\begin{enumerate}
\item $\partial q_{z,R}(x,y)=\xi_{R,+}-\xi_{R,-}$ for all $R>R_0$,
\item any sequence $\{y_i\}_{i=R_0}^\infty$ with $y_i\in \xi_{i,+}$
  satisfies $\lim_{i\to\infty} y_i = y$, and any sequence
  $\{y_i\}_{i=R_0}^\infty$ with $x_i\in \xi_{i,-}$ 
  satisfies $\lim_{i\to\infty} x_i = x$.
\end{enumerate}
\end{definition}

\begin{definition}\label{d:quasigeodesic}
Let $\Gamma$ be a graph with a compactification $\bar\Gamma$ of
$\Gamma^{(0)}$, and let 
$\Lambda\subseteq \bar{\Gamma}$.  Let $\epsilon>0$.
A homological bicombing $q: \Lambda \times \Lambda \to 
\bar{C}_1(\Gamma;\R)$ 
is \emph{$\epsilon$-quasi-geodesic} if both

\begin{enumerate}
\item $q(a,b)$ has support in the $\epsilon$-neighborhood of some
  geodesic between $a$ and $b$, and 
\item If $a$, $b\in \Gamma^{(0)}$, then $|q(a,b)|_1\leq \epsilon
  d(a,b)$.  
\end{enumerate}
\end{definition}

\begin{remark} \label{r:wrongdef}
In general, one may also want to place constraints on the $1$-norms
of finite ``subsegments" of $q(a,b)$, where $a$ and $b$ are ideal points.
However, we do not need this refinement in this paper.
\end{remark}

\subsection{Relatively hyperbolic groups}
\label{ss:relhyp}

Relatively hyperbolic groups were first defined by
Gromov in \cite{gromov:hyp}.  Alternative definitions 
were given by Farb \cite{farb:relhyp}
and Bowditch \cite{bowditch:relhyp}.  These
definitions are all now known to be equivalent.
See \cite[Appendix]{Dah:thesis}.

Further characterizations of relatively hyperbolic groups
are given by Osin \cite{osin:relhypbook}, in terms of
relative Dehn functions, and Yaman \cite{yaman}, in 
terms of convergence group actions.

Recently there has been a large amount of interest in
relatively hyperbolic groups.  (See \cite{bumagin:rh}, 
\cite{dahmani:classifyingspace}, \cite{drutusapir}, \cite{groves_rh}, \cite{osin:relhypbook}, among many others).

Here is the original definition of Gromov's \cite[Section 8.6]{gromov:hyp}:
\begin{definition}\label{d:gromovrh}
Suppose that $G$ acts isometrically and properly on a proper,
geodesic, Gromov hyperbolic metric space $X$, so that the quotient is
quasi-isometric to a wedge of $n$ rays. 
Let
$\gamma_1,\ldots,\gamma_n$ be unit-speed geodesic rays in $X/G$ tending to
distinct points in the Gromov boundary of $X/G$, and choose lifts
$\twid{\gamma}_1,\ldots,\twid{\gamma}_n$ to $X$.  
For each $i$, let $e_i$ be the point in $\partial X$ to which
$\twid{\gamma}_i$ limits, and let $P_i$ be the stabilizer in $G$ of
$e_i$.  For each $i$ define a \emph{horofunction} $h_i\co X\to \R$ by
\[ h_i(x) = \limsup_{t\to\infty} d(x,\twid{\gamma}_i(t)) - t. \]
The \emph{$R$-horoballs} of $X$ are the sub-level sets
$B_i(R) = h_i^{-1}(-\infty,R)$ and their $G$-translates.
Assume that there exists a constant $R$ so that for any $g\in G$ and
any $i,j$, either $gB_i(R)\cap B_j(R)$ is empty or $i=j$ and $g\in
P_i$.  Finally, suppose that $G$ acts cocompactly on the complement of
the union of the horoballs.
Then we say that $G$ is {\em hyperbolic relative to 
$\mc{P}=\{P_1,\ldots,P_n\}$ in the sense of Gromov}.
\end{definition}

\begin{definition} \label{d:non-elem}
Suppose that $G$ is a relatively hyperbolic group acting on the $\delta$-hyperbolic space $X$ as in Definition \ref{d:gromovrh}.  An element
$g \in G$ is called {\em hyperbolic} if it does not have a bounded orbit
in $X$, and it fixes exactly two points in $\partial X$.

We say that $G$ is {\em non-elementary relatively hyperbolic} if there
are hyperbolic elements $g,h$ in $G$ so that 
$\mathrm{Fix}_{\partial X} (g) \cap \mathrm{Fix}_{\partial X}(h) = 
\emptyset$.
\end{definition}

\begin{remark}
By the usual Ping-Pong argument, if $g,h$ are as in Definition 
\ref{d:non-elem} above, then there is some $j \ge 1$ so that $g^j$ and
$h^j$ generate a free group.
\end{remark}

We now give another definition of relatively hyperbolic
groups, which is a hybrid of Farb's \cite{farb:relhyp} and 
one of Bowditch's \cite{bowditch:relhyp}.

\begin{definition}[Coned-off Cayley graph]
\label{d:conedgraph}
Suppose that  $G$ is a finitely generated
group, with finite generating set $S$.  Let 
$\Gamma(G,S)$ be the Cayley graph of $G$
with respect to $S$.

Suppose that $\P = \{ P_1, \ldots , P_k \}$ is a finite
collection of finitely generated subgroups of $G$.
We form a new graph containing $\Gamma(G,S)$,
called the {\em coned-off Cayley graph} and denoted
$\widehat{\Gamma}(G,\P,S)$ as follows:

For each $i \in \{ 1 , \ldots , k \}$ and each coset
$g P_i$ we add a new vertex, $v_{g,i}$ to
$\Gamma(G,S)$.  We also add a vertex from
each element of $gP_i$ to $v_{g,i}$.
\end{definition}

\begin{definition}
[Fine graphs; See {\cite[page 11]{bowditch:relhyp}}]
A (not necessarily locally finite) graph $\mathcal K$ is
{\em fine} if for every edge $e$ in $\mathcal K$ and each
integer $L > 0$, the number of simple simplicial loops
of length at most $L$ which contain $e$ is finite.
\end{definition}

\begin{definition} \label{d:rh}
Suppose that $G$ is a finitely generated group and that
$\P = \{ P_1 , \ldots , P_k \}$ is a collection of finitely
generated subgroups.  We say that $G$ is {\em hyperbolic
relative to $\P$} if the coned-off Cayley graph
$\widehat{\Gamma}(G,\P,S)$ is fine and $\delta$-hyperbolic
for some $\delta > 0$.
\end{definition}

\begin{remark}\label{r:bcp}
By now the class of groups which we call `relatively hyperbolic'
is standard.  However, we should point out that, in the terminology of
\cite{farb:relhyp}, $G$ is hyperbolic relative to $\mc{P}$ if 
and only if the
coned-off Cayley graph is $\delta$-hyperbolic.  
Farb's hypothesis of Bounded Coset Penetration is equivalent to
fineness of $\widehat{\Gamma}$ (See for instance 
\cite[Appendix]{Dah:thesis}).  
It is shown in \cite{bumagin:rh} and \cite[Appendix]{Dah:thesis}
that Definitions \ref{d:gromovrh} and \ref{d:rh} are equivalent.
\end{remark}

Whenever $G$ is hyperbolic relative to $\mc{P}$, we will always
assume that our (finite) generating set for $G$ is compatible,
in the sense of Definition \ref{d:compatible}.

We briefly list some examples of relatively hyperbolic groups (and the
subgroups they are hyperbolic relative to):
\begin{enumerate}
\item Hyperbolic groups are hyperbolic relative to the empty collection
of subgroups;
\item Fundamental groups of geometrically finite hyperbolic manifolds
are hyperbolic relative to the cusp subgroups;
\item Free products are hyperbolic relative to the free factors;
\item A group which acts properly and cocompactly on a CAT$(0)$ space
with isolated flats is hyperbolic relative to the stabilizers of maximal
flats (see \cite{HruskaKleiner}; and \cite{drutusapir});
\item Limit groups are hyperbolic relative to maximal noncyclic abelian
subgroups (see \cite{alibegovic} and \cite{dah:comb}).
\end{enumerate}

In Section \ref{section:cusped} below, we will introduce
a `cusped' space, $X(G,\P,S)$, associated to a group $G$ 
and finite collection $\P$ of finitely generated subgroups.
We will prove that $G$ is hyperbolic relative to $\P$ if
and only if $X(G,\P,S)$ is $\delta$-hyperbolic for some 
$\delta$.

\begin{remark}
One of the important features of the space $X$ defined in
Section \ref{section:cusped} below is that if $G$ is hyperbolic
relative to $\mc{P}$, then the action of $G$ on 
$X(G,\P,S)$ satisfies the requirements of Definition
\ref{d:gromovrh} (see Theorem \ref{t:tfae} below).  

There are a number of previous constructions of cusped
spaces (for example: Cannon and Cooper \cite{CannonCooper},
Bowditch \cite{bowditch:relhyp} and Rebbechi \cite{rebbechi}.
The novelty in our space $X$ is that it is a graph metrized
so that edge lengths are $1$.  Thus we can apply many
combinatorial constructions directly, such as isoperimetric
inequalities (combinatorial and homological).  Importantly,
in Section \ref{section:bicombing} below, we can also apply
a construction of Mineyev from \cite{min:str}.  

Another important feature of the space $X$ is that there
is a bound on the lengths of the attaching maps of $2$-cells.
\end{remark}

By Osin \cite[Theorem 1.5]{osin:relhypbook}, relatively 
hyperbolic groups
are always finitely presented relative to their parabolics. (This 
also follows from the construction of the {\em relative Rips
complex} in \cite{dahmani:classifyingspace}.)

\begin{theorem} \label{t:rhrelfinpres}
\cite[Theorem 1.5]{osin:relhypbook}
Let $G$ be a finitely generated group, $\{ H_1 ,\ldots , H_m \}$
a collection of subgroups of $G$.  The following are equivalent:
\begin{enumerate}
\item $G$ is finitely presented with respect to $\{ H_1, \ldots , H_m \}$
and the corresponding relative Dehn function is linear;
\item $G$ is hyperbolic relative to $\{ H_1, \ldots , H_m \}$.
\end{enumerate}
\end{theorem}

Recall by Lemma \ref{l:Hfg} that if $G$ is finitely generated, and finitely presented relative
to $\{ H_1, \ldots , H_m \}$, then each of the $H_i$ is finitely generated.

In this paper, we have no need for the `relative' Dehn functions of 
\cite{osin:relhypbook}.  Rather, we construct simply-connected
$2$-complexes
with linear combinatorial isoperimetric inequalities.  (They
also have linear homological isoperimetric inequalities; see
Theorem \ref{t:tfae} below.)

\begin{definition} [Coned-off Cayley complex]
\label{d:conedcomplex}
Suppose that $G$ is a finitely generated group, with
a collection $\P$ of finitely generated subgroups, and that 
$\langle \mc{A}, \P \mid \mc{R} \rangle$ is a finite
relative presentation for $G$.

Let $S$ be a (finite) compatible generating set for $G$ containing
$\mc{A}$.  Form a $2$-complex $\hat{C}(G,\P,S,\mc{R})$, called
the {\em coned-off Cayley complex} as follows:

Let $C= C(G,\P,S)$ be the coned-off Cayley graph (which contains
a copy of the Cayley graph $\Gamma(G,S)$.  Attach $2$-cells
to $C$ in a $G$-equivariant way, corresponding to the relations 
$\mc{R}$.  Also, add a $2$-cell to each loop of length three in $C$
which contains an infinite valence vertex.
\end{definition}

\begin{lemma}
The coned-off Cayley complex is simply-connected.
\end{lemma}

\begin{proposition} \label{p:conelii}
Suppose that $G$ is hyperbolic relative to $\P$, and let
$S$ be a finite generating set for $G$ containing generating sets
for each subgroup in $\P$.  Also, let $\langle S,\P \mid \mc{R} \rangle$
be a finite relative presentation for $G$.

Form the coned-off Cayley complex $\hat{C}(G,\P,S,\mc{R})$.

Then $\hat{C}(G,\P,S,\mc{R})$ has a linear combinatorial isoperimetric
inequality.
\end{proposition}
\begin{proof}
Since $G$ is hyperbolic relative to $\P$, the coned-off Cayley
graph is $\delta$-hyperbolic and fine.  Therefore there are only 
finitely many orbits
of simple loops of length at most $16\delta$.  Since $\hat{C}(G,\P,S,\mc{R})$
is simply-connected, this
implies that there exists a $K$ so that every combinatorial loop of 
length at most $16\delta$ can be filled with a combinatorial disk
in $\hat{C}(G,\P,S,\mc{R})$ of area at most $K$. 

The result now follows immediately from Proposition \ref{p:hyptocomb}.\end{proof}

We now prove the converse to Proposition \ref{p:conelii}.

\begin{proposition} \label{p:coneliitorh}
Suppose that $G$ is finite presented relative to $\P
= \{ P_1 ,\ldots  , P_n \}$ and that
$\langle S,\P \mid \mc{R} \rangle$ is a finite relative
presentation for $G$.

If the coned-off Cayley complex $\hat{C}=\hat{C}(G,\P,S,\mc{R})$ satisfies
a linear combinatorial isoperimetric inequality then
the coned-off Cayley graph $\hat\Gamma = \hat{\Gamma}(G,\P,S)$ is 
Gromov hyperbolic and fine.
\end{proposition}
\begin{proof}
That $\hat\Gamma$ is $\delta$-hyperbolic for some $\delta$ follows 
from Proposition \ref{p:combtohyp}, since there is certainly a bound
on the length of the attaching maps of $2$-cells in $\hat{C}$.

It remains to prove that $\hat\Gamma$ is fine.  Take an edge $e$
and a simple simplicial loop $c$ containing $e$ of length $L$.  Now,
there are only finitely many $2$-cells adjacent to each edge in $\hat{C}$.
A simple loop may be filled by a topological disk, and the isoperimetric
function gives a bound on the area of such a disk in terms of the length.
In particular, suppose a simple loop of length at most $L$ can be 
filled by a  topological disk of area at most $L'$.  There are only 
finitely many topological disks of area at most $L'$ containing $e$ on 
the boundary, so there are only finitely many simple loops of length 
at most $L$ containing $e$.  Therefore $\hat\Gamma$ is fine.
\end{proof}

It is a consequence of Theorem \ref{t:tfae} that the hypothesis in
Proposition \ref{p:coneliitorh} can be weakened to that of having a
{\em homological} linear isoperimetric inequality, but we do not know
a direct proof.  This raises a natural question (if the answer is
positive, then it can be used to give such a direct proof):
\begin{question}
Let $X$ be a simply connected $2$-complex with a homological (linear?)
isoperimetric inequality and a bound on the length of attaching maps
of $2$-cells and finitely many $2$-cells adjacent to any edge.  Must
$X$ be fine?
\end{question} 

\newpage
\part{The cusped space and preferred paths}

\section{The cusped space}\label{section:cusped}
The purpose of this section is to construct a space
$X$ from a finitely generated group $G$, and a finite
collection $\mc{P}$ of finitely generated subgroups.  The
utility of $X$ is that it is Gromov hyperbolic if and only if
$G$ is hyperbolic relative to $\mc{P}$ (see Theorem \ref{t:tfae}).

\subsection{Combinatorial horoballs} \label{ss:combhoro}
\begin{definition}
Let $\Gamma$ be any $1$-complex.
The \emph{combinatorial horoball based on $\Gamma$}, denoted
$\mc{H}(\Gamma)$, is the $2$-complex formed as follows:
\begin{itemize}
\item $\mc{H}^{(0)}= \Gamma^{(0)}\times \left( \{0\}\cup \N \right)$
\item $\mc{H}^{(1)}$ contains the following three types of edges.  The
  first two types are called \emph{horizontal}, and the last type is
  called \emph{vertical}.
\begin{enumerate}
\item If $e$ is an edge of $\Gamma$ joining $v$ to $w$ then there is a
  corresponding edge $\bar{e}$ connecting $(v,0)$ to $(w,0)$.
\item If $k>0$ and $0<d_{\Gamma}(v,w)\leq 2^k$, then there is a single edge
  connecting $(v,k)$ to $(w,k)$.
\item If $k\geq 0$ and $v\in \Gamma^{(0)}$, there is an edge  joining
  $(v,k)$ to $(v,k+1)$. 
\end{enumerate}
\item $\mc{H}^{(2)}$ contains three kinds of $2$-cells:
\begin{enumerate}
\item If $\gamma\subset \mc{H}^{(1)}$ is a circuit composed of three
  horizontal edges, then there is a $2$-cell (a \emph{horizontal
  triangle}) attached along $\gamma$.
\item If $\gamma\subset \mc{H}^{(1)}$ is a circuit composed of two
  horizontal edges and two vertical edges, then there is a $2$-cell (a
  \emph{vertical square}) attached along $\gamma$.
\item If $\gamma\subset \mc{H}^{(1)}$ is a circuit composed of three
  horizontal edges and two vertical ones, then there is a $2$-cell (a
  \emph{vertical pentagon}) attached along $\gamma$, unless $\gamma$
  is the boundary of the union of a vertical square and a horizontal
  triangle.
\end{enumerate}
\end{itemize}
\end{definition}

\begin{remark}\label{r:subset}
As the full subgraph of $\mc{H}(\Gamma)$ containing the vertices
$\Gamma^{(0)}\times\{0\}$ is isomorphic to $\Gamma$, we may think of
$\Gamma$ as a subset of $\mc{H}(\Gamma)$.
\end{remark}

\begin{remark}
Whenever $\mc{H}(\Gamma)$ is to be thought of as a metric space, we
will always implicitly ignore the $2$-cells, and regard
$\mc{H}(\Gamma)^{(1)}$ as a metric graph with all edges of length
one.
\end{remark}

\begin{definition} \label{d:Depth1}
Let $\Gamma$ be a graph and $\mc{H}(\Gamma)$ the associated
combinatorial horoball.  Define a \emph{depth} function
\[D: \mc{H}(\Gamma) \to [0,\infty)\]
which satisfies:
\begin{enumerate}
\item $D(x)=0$ if $x\in \Gamma$,
\item $D(x)=k$ if $x$ is a vertex $(v,k)$, and
\item $D$ restricts to an affine function on each $1$-cell and on each
  $2$-cell.
\end{enumerate}
\end{definition}

\begin{definition} \label{d:toppart}
Let $\Gamma$ be a graph and  $\mc{H} = \mc{H}(\Gamma)$ the associated
combinatorial horoball.  For $N \ge 1$, let $\mc{H}_N \subset \mc{H}$
be the full sub-graph with vertex set $\Gamma^{(0)} \times \{ 0, \ldots , N\}$.
\end{definition}

The following observation will be important in Section \ref{surgered}.

\begin{observation} \label{o:2step}
Let $\Gamma$ be a graph, $\mc{H}_N$ as in Definition \ref{d:toppart} above, and $\partial \mc{H}_N$ be the full sub-graph with vertex set
$\Gamma^{(0)} \times \{ N \}$.

Let $\mc{H}' = \mc{H}(\partial \mc{H}_N)$.  Identify the copies of
$\partial \mc{H}_N$ in $\mc{H}_N$ and $\mc{H}'$.  The resulting
complex is isomorphic to $\mc{H}$.
\end{observation}

\begin{proposition} \label{p:combhorolii}
Let $\Gamma$ be a connected
$1$-complex so that no edge joins a vertex to itself. Then $\mc{H}(\Gamma)$ is simply-connected and
satisfies a linear combinatorial isoperimetric inequality with
constant at most $3$.
\end{proposition}
\begin{proof}
Let $c$ be a combinatorial loop in $\mc{H}(\Gamma)$.  To
prove the proposition it suffices to show that $c$ can be filled by
a disk of area at most $3|c|_1$.

Let $j$ be minimal so that there is some vertex $(v,j)$ in $c$.
Since we may clearly suppose $c$ has no backtracking, there
is at least one horizontal edge at depth $j$ in $c$.
By gluing pentagons to edges in $c$ at depth $j$, and squares
when pentagons are not possible, we can reduce the length
of $c$, and increase the minimal depth.  Repeating this procedure,
we eventually end with a path of length $3$, or $4$, which is entirely
horizontal.  A path of length $3$ can be filled be a horizontal triangle,
whilst a path of length $4$ can be filled by two pentagons beneath it.
Being slightly careful about counting shows that the isoperimetric 
constant is at most $3$, as required.
\end{proof}

\begin{theorem} \label{t:combhorohyp}
Let $\Gamma$ be any $1$-complex.  Then 
$\mc{H}(\Gamma)^{(1)}$ is $\delta$-hyperbolic, where
$\delta$ is independent of $\Gamma$.
\end{theorem}
\begin{proof}
This follows from Proposition \ref{p:combhorolii} and Proposition 
\ref{p:combtohyp}.
\end{proof}

\begin{remark}
By studying the geometry of geodesics in combinatorial horoballs
as in the results below, it is possible to directly prove that any
combinatorial horoball is $20$-hyperbolic (and $20$ is not optimal).
\end{remark}

Geodesics in combinatorial horoballs are particularly easy to understand.

\begin{lemma}\label{l:gamma}
Let $\mc{H}(\Gamma)$ be a combinatorial horoball.
Suppose that $x, y \in \mc{H}(\Gamma)$ are distinct vertices.  
Then there is a geodesic
$\geod(x,y)=\geod(y,x)$
between $x$ and $y$ which consists of at most two vertical
segments and a single horizontal segment of length at most $3$. 

Moreover, any other geodesic between $x$ and $y$
is Hausdorff distance at most $4$ from this 
geodesic.
\end{lemma}
\begin{proof}
Let $\gamma'$ be any geodesic joining $x$ to $y$. 

We observe that if 
$h=[\gamma'(t_1),\gamma'(t_2)]$ is a maximal horizontal segment of
length greater than $1$, then 
$D(\gamma'(t_1-1))$ and $D(\gamma'(t_2+1))$ are both smaller than
$D(h)$ (see Figure \ref{f:shortcut2}).
\begin{figure}[htbp]
\begin{center}
\begin{picture}(0,0)%
\includegraphics{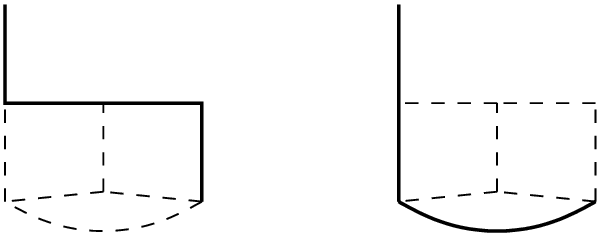}%
\end{picture}%
\setlength{\unitlength}{4144sp}%
\begingroup\makeatletter\ifx\SetFigFont\undefined%
\gdef\SetFigFont#1#2#3#4#5{%
  \reset@font\fontsize{#1}{#2pt}%
  \fontfamily{#3}\fontseries{#4}\fontshape{#5}%
  \selectfont}%
\fi\endgroup%
\begin{picture}(2744,1079)(2229,-1118)
\end{picture}%
\caption{In a geodesic, a horizontal segment not at maximal depth must
have length $1$.}
\label{f:shortcut2}
\end{center}
\end{figure}
It is easy to see that no geodesic in $A$ can contain a horizontal
segment of length longer than $5$ (Figure \ref{f:shortcut6}).
\begin{figure}[htbp]
\begin{center}
\begin{picture}(0,0)%
\includegraphics{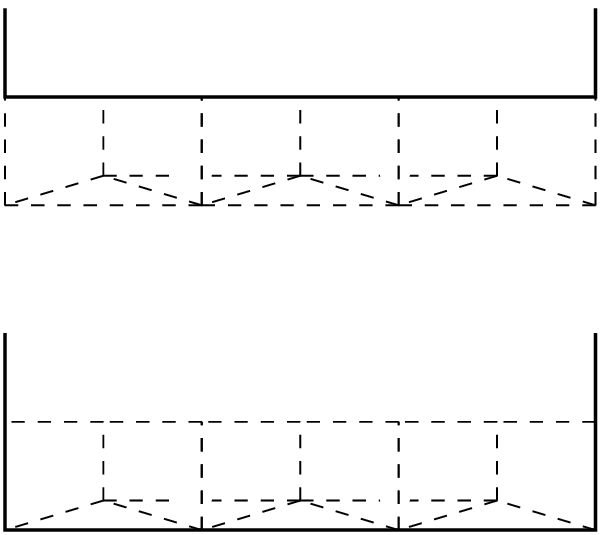}%
\end{picture}%
\setlength{\unitlength}{4144sp}%
\begingroup\makeatletter\ifx\SetFigFont\undefined%
\gdef\SetFigFont#1#2#3#4#5{%
  \reset@font\fontsize{#1}{#2pt}%
  \fontfamily{#3}\fontseries{#4}\fontshape{#5}%
  \selectfont}%
\fi\endgroup%
\begin{picture}(2744,2429)(2229,-2288)
\end{picture}%
\caption{A geodesic can contain no horizontal segment of length
greater than $5$.}
\label{f:shortcut6}
\end{center}
\end{figure}
Indeed, the geodesic $\gamma'$ can contain at most $5$ horizontal
edges in total.  

We next observe that there can be at most two horizontal segments in
$\gamma'$ other than the one at maximal depth, at most one on each side of the deepest one.  In fact, there
can be at most \emph{one} horizontal segment other than the one at
maximal depth.   Figure \ref{f:threehoriz}  
\begin{figure}[htbp]
\begin{center}
\begin{picture}(0,0)%
\includegraphics{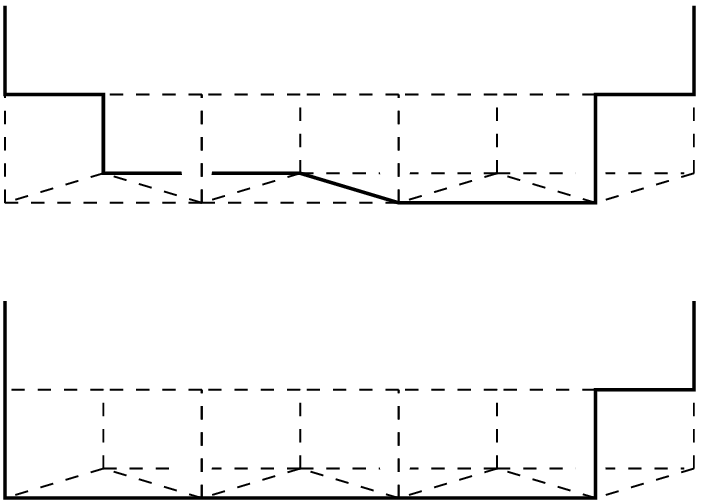}%
\end{picture}%
\setlength{\unitlength}{4144sp}%
\begingroup\makeatletter\ifx\SetFigFont\undefined%
\gdef\SetFigFont#1#2#3#4#5{%
  \reset@font\fontsize{#1}{#2pt}%
  \fontfamily{#3}\fontseries{#4}\fontshape{#5}%
  \selectfont}%
\fi\endgroup%
\begin{picture}(3194,2294)(2229,-2153)
\end{picture}%
\caption{A path with three horizontal segments and a way to shorten
the path.}
\label{f:threehoriz}
\end{center}
\end{figure}
shows representative paths with three horizontal segments together with
ways to shorten them (there are other possibilities, which are easy
to deal with).

Let $M'$ be the
maximum depth achieved by $\gamma'$.  There is a geodesic $\gamma''$
obtained from $\gamma'$ by pushing all horizontal segments of
$\gamma'$ down to $D^{-1}(M')$.  The Hausdorff distance between
$\gamma'$ and $\gamma''$ is at most $1\frac{1}{2}$.
The geodesic $\gamma''$ consists of at most two
vertical segments and one horizontal segment of length at most $5$.
If the horizontal segment $h \subset \gamma''$ has length $4$ 
(respectively $5$) then there is another path with the same endpoints as $h$, and the same length, consisting of two
vertical edges and a horizontal path of length $2$ (respectively $3$).
Replacing $h$ with this new path if necessary, we obtain a geodesic
of the form required by the lemma.

Now let $\gamma$ be {\em any} geodesic satisfying the conclusion
of the lemma.  We argue that the Hausdorff distance between
$\gamma$ and $\gamma''$ is at most $2\frac{1}{2}$ which, combined
with the earlier estimate on the Hausdorff distance between $\gamma'$
and $\gamma''$, completes the proof.

Let $M$ be the maximum depth of $\gamma$.  The reader
may readily verify that $0 \le M - M' \le 1$.  
If $M=M'$, then $\gamma''$ lies within Hausdorff
distance at most $1\frac{1}{2}$ of $\gamma$. 
If $M=M'+1$, then the Hausdorff distance between $\geod(x,y)$ and
$\gamma''$ is at most $2\frac{1}{2}$. 
\end{proof}

\begin{lemma}
If $A$ is a combinatorial horoball, then the Gromov boundary
consists of a single point, denoted $e_A$.  Moreover, for
any $x \in A$, there is a geodesic ray from $x$ to $e_A$ consisting
entirely of vertical edges.  Any geodesic ray from $x$ to $e_A$
is Hausdorff distance at most $1\frac{1}{2}$ from the vertical ray.
\end{lemma}
\begin{proof}
A geodesic ray has at most one horizontal edge.  Given this
observation, the proof is similar to that of Lemma \ref{l:gamma}.
\end{proof}

\subsection{The augmentation}

Let $G$ be a finitely generated group, with a finite
collection $\mc{P}$ of subgroups.  Let $S$ be a compatible
generating set for $G$.  We define an
augmentation of the relative Cayley graph complex of $G$
by combinatorial horoballs.  This
augmentation will be hyperbolic exactly when 
$G$ is hyperbolic relative to $\mc{P}$ (see Theorem
\ref{t:tfae} below).

\begin{definition}\label{d:X1}
Let $G$ be a finitely generated group, let
$\mc{P}=\{P_1,\ldots,P_n \}$ be a (finite) 
family of finitely generated subgroups of $G$, and let $S$ be a
generating set for $G$ so that $P_i\cap S$ generates $P_i$ for each
$i\in \{1,\ldots,n\}$.  
For each $i\in \{1,\ldots,n\}$, let $T_i$ be a left transversal for
$P_i$ (i.e. a collection of representatives for left cosets of $P_i$ in $G$
which contains exactly one element of each left coset). 

For each $i$, and each $t\in T_i$, let $\Gamma_{i,t}$ be the full
subgraph of the Cayley graph $\Gamma(G,S)$ which contains $tP_i$.
Each $\Gamma_{i,t}$ is isomorphic to the Cayley graph of $P_i$ with
respect to the generators $P_i\cap S$.
Then we define
\[ X = \Gamma(G,S) \cup (\cup \{\mc{H}(\Gamma_{i,t})^{(1)}\mid 1\leq i\leq n,
t\in T_i\}),\]
where the graphs $\Gamma_{i,t}\subset \Gamma(G,S)$ and
$\Gamma_{i,t}\subset \mc{H}(\Gamma_{i,t})$ are identified as suggested
in Remark \ref{r:subset}.
\end{definition}

\begin{remark}
The vertex set of $X$ can naturally be identified with the set of $4$-tuples $(i,t,p,k)$, where
$i\in\{1,\ldots,n\}$, $t\in T_i$, $p\in P_i$, and $k\in \N$.

We will use this identification without comment in the sequel.
\end{remark}

\begin{remark}
The group $G$ acts isometrically and properly on the graph $X(G,\mc{P},S)$.
\end{remark}

A path in $X$ starting at $1$ determines an element of $G$ in the following manner:  Each horizontal edge is naturally labelled by a
group element.  Take the product of the labels of the horizontal
edges in the path (and ignore the vertical edges) with the order
coming from the path. 

Paths in quotients $X/H$, where $H \le G$, and
also paths which start at a point $(i,1,1,k)$ which lies directly
`beneath' $1$ also naturally determine elements of $G$.

Supposing further that $G$ is finitely presented relative to $\mc{P}$,
there is a natural locally finite simply connected $2$-complex with
skeleton $X(G,\mc{P},S)$:

\begin{definition}\label{d:X2}
Let $G$ be a finitely generated group, let
$\mc{P}=\{P_1,\ldots,P_n \}$ be a (finite) 
family of finitely generated subgroups of $G$, and let $S$ be a
generating set for $G$ so that $P_i\cap S$ generates $P_i$ for each
$i\in \{1,\ldots,n\}$.  Let $S' = S\setminus\cup\mc{P}$, and
suppose that
\[G = \langle \mc{A},\mc{P}\mid \mc{R}\rangle\]
is a finite \emph{relative} presentation of $G$, in the sense of
Definition \ref{d:relpres} above (where $\mc{A} \subset S$).
Then we may form a locally finite $2$-complex $X(G,\mc{P},S,\mc{R})$,
whose one-skeleton is the space $X(G,\mc{P},S)$ from Definition
\ref{d:X1}, and which contains the following $2$-cells:
\begin{itemize}
\item $2$-cells from $\mc{R}$: Each relator in $\mc{R}$ determines
  a
  loop $\gamma$ in $\Gamma(G,S)$ beginning at $1$.  For each
  $g\in G$, there is a $2$-cell attached along $g\gamma$.
\item $2$-cells from combinatorial horoballs:  For each $i$ and each
  $t\in T_i$ we have an embedding of $\mc{H}(\Gamma_{i,t})^{(1)}$ into
  $X(G,\mc{P},S)$.  If $\gamma\in \mc{H}(\Gamma_{i,t})$ bounds a
  $2$-cell, then there is a $2$-cell attached along its image in
  $X(G,\mc{P},S,\mc{R})$.  
\end{itemize}
\end{definition}

\begin{remark} \label{r:alt}
The $2$-complex $X$ can also be obtained from the relative 
Cayley complex
$C(G,S,\mc{R})$ as in Definition \ref{d:relCayleycomplex} by attaching
combinatorial horoballs to the cosets $gP_i$ (including the 
$2$-cells in the horoballs).  In particular, and this will be important
a number of times throughout this paper, the relative Cayley complex
is canonically embedded in $X$.
\end{remark}

As in Definition \ref{d:Depth1}, we define a depth function for
the space $X(G,\mc{P},S,\mc{R})$:
\begin{definition}
Define a \emph{depth} function
\[D: X(G,\mc{P},S,\mc{R})\to [0,\infty)\]
which satisfies:
\begin{enumerate}
\item $D(x)=0$ if $x\in G$,
\item $D(x)=n$ if $x$ is a vertex $(i,t,p,n)$,
\item $D$ is equivariant, and 
\item $D$ restricts to an affine function on each $1$-cell and on each
  $2$-cell.
\end{enumerate}
\end{definition}

\begin{remark}
Because $D$ is $G$-equivariant, it induces a depth
function on the quotient space $X/H$, for any subgroup
$H$ of $G$.  We refer to the depth function on the 
quotient by $D$ also.
\end{remark}

\begin{remark}
We observe that $D^{-1}(0) = \CComp$, the relative Cayley complex of
$G$.
\end{remark}

\begin{remark}
In Section \ref{section:bicombing} below, we need to choose
$\mc{R}$ carefully.  We will need all sufficiently short loops in
$\CComp$ to be able to be filled in $\CComp$.  This can be ensured
by including all of the short relations in $G$ in $\mc{R}$.
\end{remark}

We remark that $\Xgraph$ naturally breaks up into 
$\CComp$ and an
equivariant family of combinatorial horoballs.

\begin{definition}\label{d:horoball}
Let $L>0$. An \emph{$L$-horoball} is a component of
$D^{-1}[L,\infty)$.  A \emph{$0$-horoball} is the maximal subcomplex of
  $X(G,\mc{P},S,\mc{R})$ with vertices 
\[tP_i\cup \{(i,t,p,k)\mid p\in P_i, k\in\N \} \]
for some $i\in\{1,\ldots,n\}$ and some $t\in T_i$.
\end{definition}

\begin{remark}
The space $X$ is given the path metric, where each edge has
length $1$.
We haven't specified what kind of metric to put on the $2$-cells.
Thus, whenever we are discussing the metric of $X$, we will
simply pretend that $X$ is a $1$-complex.  In particular, a {\em 
geodesic} in $X$ will always refer to a geodesic path in the 
$1$-skeleton.

The purpose of the $2$-cells is that we want a simply-connected 
space for which we will prove a linear isoperimetric 
inequality.  This will imply that $\Xgraph$ is a $\delta$-hyperbolic
space for some $\delta$, in case $G$ is hyperbolic relative to
$\P$.  See Theorem \ref{t:tfae} below.
\end{remark}

\subsection{Notions of relative hyperbolicity}
The main result of this subsection is Theorem \ref{t:tfae}, which
gives a collection of statements which are equivalent
to relative hyperbolicity.
In particular, $G$ is hyperbolic relative to $\P$ if and only
if $X(G,\P,S)$ is Gromov hyperbolic for any appropriate
choice of $S$ as in Definition \ref{d:X1} above.

The next result will allow us to translate the hyperbolicity
of the coned-off Cayley complex into hyperbolicity of the
space $X$.  In order to make the proof of Theorem \ref{t:tfae}
easier, we choose to phrase it in terms of linear homological
isoperimetric inequalities, though there is of course an analogous
version with combinatorial isoperimetric inequalities (and the proof
of this analogue is somewhat easier).

\begin{theorem} \label{t:Xlii}
Suppose that $G$ is finitely presented relative to
$\P$, that $S$ is a finite
generating set for $G$ so that $P = \langle S \cap P \rangle$
for each $P \in \P$, and suppose that
$\langle S , \P \mid \mc{R} \rangle$ is a finite
relative presentation for $G$.

If the coned-off Cayley complex $\hat{C}(G,\P,S,\mc{R})$
satisfies a linear homological isoperimetric inequality
then $X(G,\P,S,\mc{R})$ satisfies a linear homological
isoperimetric inequality.
\end{theorem}
\begin{proof}
Let $\Gamma = \Gamma(G,S)$ be the Cayley 
graph of $G$
with respect to $S$.  Let $\hat{\Gamma} = 
\hat{\Gamma}(G,\P,S)$ be the 
coned-off Cayley
graph and $\hat{C} =\hat{C}(G,\P,S,\mc{R})$ the coned-off Cayley
complex.

By assumption,
$\hat{C}$ satisfies a linear homological 
isoperimetric inequality. Let $K$ be the isoperimetric
constant for $\hat{C}$.

Take a $1$-cycle $c$ in $X^1$.  By Theorem 
\ref{t:coherent}, we may assume without loss of generality
that $c$ is a simple loop.  If the support of 
$c$ lies entirely within a single horoball, then $c$ may be filled
with a combinatorial disk of area at most $3|c|_1$, by
Proposition \ref{p:combhorolii}.  Therefore, we may suppose that
the support of $c$ does not lie entirely in a horoball.

Decompose
$c$ into pieces which lie in $\Gamma$,
and pieces which lie entirely in a single horoball.

For each maximal sub-path of $c$ which lies in a
horoball, there is a path of length $2$ in 
$\hat{\Gamma}$ with the same endpoints.  This 
gives a loop $\hat{c}$ in $\hat{\Gamma}$.
Clearly $|\hat c| \le |c|$.

There is a rational $2$-chain $\omega$ in $\hat{C}$ with
boundary $\hat c$ so that
\[	|\omega|_1 \le K |\hat c|_1 .		\]

Consider an infinite-valence vertex $v \in \hat\gamma$, and
the $2$-cells in $\text{supp}(\omega)$ which intersect $v$.
Each such $2$-cell is a triangle, and has a single edge in $\Gamma$.
Thus to each infinite valence vertex $v$ in $\hat c$ is associated (via
the $2$-chain $\omega$) a $1$-chain $p_v$ whose support is contained
in the link of $v$ and whose $1$-norm is exactly the $1$-norm of 
$\omega$ restricted to the star of $v$.  

The $1$-chains $p_v$ (of which there are finitely many for the given
path $\hat c$) have total length bounded by $|\omega|_1$, and 
give a method of `surgering' $c$.  We have the decomposition
of $c$ into paths $q_i$ in the Cayley graph and paths $r_v$, where
$r_v$ is a path through a horoball corresponding to the infinite valence
vertex $v$.  This induces a decomposition
of $c$ as 
\[	c = \sum_v (r_v + p_v) + \left( \sum_i q_i - \sum_v p_v \right),	\]
where the paths are oriented so that $r_v + p_v$ is a $1$-cycle,
and so is $\left( \sum_i q_i - \sum_v p_v \right)$.

Now,
\begin{eqnarray*}
|c|_1 & \le & \sum_v |r_v + p_v|_1 + \sum_i |q_i|_1 + \sum_v |p_v|_1 \\
& \le & 2K |\hat c|_1 + |c|_1\\
& \le & (2K+1) |c|_1.
\end{eqnarray*}
(The second inequality holds because $\sum_v |p_v|_1 \le |\omega|_1$,
and because $\sum_v |r_v|_1 + \sum_i |q_i|_1 = |c|_1$.)

The $1$-cycles $r_v+p_v$ lie entirely in a horoball, and can be filled
efficiently, as described above (this follows from Theorem 
\ref{t:coherent} and Proposition \ref{p:combhorolii}).

Therefore, it suffices to fill the $1$-cycle $c_1 =  \left( \sum_i 
q_i - \sum_v p_v \right)$ efficiently.  The $1$-cycle $c_1$ has 
support entirely in the Cayley graph $\Gamma$.
Therefore, $c_1$ can be interpreted as a $1$-cycle in $\hat\Gamma$.
Hence there is a $2$-chain $\omega_1$ in $\hat\Gamma$ so that
$\partial \omega_1 = c_1$ and $|\omega_1|_1 \le K |c_1|_1$.

Let $\omega_1'$ be the $2$-chain which is equal to $\omega_1$
on the Cayley complex, and equal to zero on $2$-cells adjacent
to infinite valence vertices.  Then $\partial \omega_1'$ decomposes
into a sum of $c$, and a collection of $1$-cycles $d_i$ each of which
lies entirely on the link of an infinite valence vertex.  Also,
$|\omega_1'|_1 \le |\omega_1|_1$, so $|\partial \omega_1'|_1 \le
M |\omega_1'|_1 \le M |\omega_1|_1$, where $M$ is the
maximum length of an attaching map of a $2$-cell in $\hat\Gamma$.

Since the support of $\omega_1'$ is entirely contained in the Cayley
complex, it can be considered as a $2$-chain in $X$, whose boundary
is the sum of $c_1$ and a collection of $1$-cycles, $d_i$ each of
which lives entirely in a single parabolic coset.  Each $d_i$ can
be filled by a $2$-chain $\nu_i$ whose support is entirely contained
in the appropriate combinatorial horoball, so that 
$|\nu_i|_1 \le 3 |d_i|_1$.  We also have that 
$\sum_i |d_i|_1 \le |\partial \omega_1'|_1$, by choosing the $d_i$ to
have distinct supports (choose one $d_i$ for each parabolic coset).

Now, by choosing appropriate orientations, we have
$\partial (\omega_1' + \sum_i \nu_i) = c_1$.  We also have
\begin{eqnarray*}
| \omega_1' + \sum_i \nu_i |_1 & \le & |\omega_1'|_1 + \sum_i |\nu_i|_1\\
& \le & |\omega_1|_1 + 3\sum_i |d_i|_1 \\
& \le & K |c_1|_1 + 3|\partial \omega_1'|_1 \\
& \le & K|c_1|_1  + M|\omega_1|_1 \\
& \le & (KM+K)|c_1|_1.
\end{eqnarray*}
This finishes the proof that $X$ satisfies a linear homological
isoperimetric inequality.
\end{proof}

We state the combinatorial version of Theorem \ref{t:Xlii} below 
for completeness.
The proof is entirely analogous to that of Theorem \ref{t:Xlii}
with loops playing the part of $1$-cycles and disks the part of
filling $2$-chains.
\begin{theorem}
Let $G, \P, S$ and $\mc{R}$ be as in the statement of 
Theorem \ref{t:Xlii} above.

If $\hat{C}(G,\P,S,\mc{R})$ satisfies a linear combinatorial isoperimetric
inequality with constant $K$ then
$X(G,\P,S,\mc{R})$ satisfies a linear combinatorial isoperimetric
inequality, and we can take the constant to be $K_1  = 3K(2K+1)$.
\end{theorem}

The following is the main result of this section, and 
gathers together a few notions of relative hyperbolicity,
including some that are new in this paper.

\begin{theorem} \label{t:tfae}
Suppose that $G$ is a finitely generated group,
that $\P = \{ P_1, \ldots , P_n \}$ is a finite collection of
finitely generated subgroups of $G$, that
$G = \langle \mc{A}, \P \mid \mathcal R \rangle$ is a finite
relative presentation for $G$, and $S$ is a compatible
generating set containing $\mc{A}$.

Let $\hat{\Gamma}$ be the coned-off Cayley graph for 
$\Gamma$ with respect to $S$ and $\P$, and let
$\hat{C}$ be the coned-off Cayley complex.  Let
$X(G,\P,S,\mc{R})$ be as defined in Definition \ref{d:X2} above.
The following are equivalent:
\begin{enumerate}
\item \label{tfae1} $G$ is hyperbolic relative to $\P$ in the sense of
Gromov;
\item \label{tfae2} $G$ is hyperbolic relative to $\P$ (i.e. $\hat\Gamma$ is Gromov hyperbolic and fine);
\item \label{tfae3}$\hat{C}$ satisfies a linear combinatorial isoperimetric
inequality;
\item \label{tfae4} $\hat{C}$ satisfies a linear homological isoperimetric
inequality;
\item \label{tfae5} $X^{(1)}$ is Gromov hyperbolic;
\item \label{tfae6} $X$ satisfies a linear combinatorial isoperimetric
inequality;
\item \label{tfae7} $X$ satisfies a linear homological isoperimetric
inequality.
\end{enumerate}
\end{theorem}
\begin{proof}
\cd{ & \eqref{tfae5}\ar @{=>}[r]\ar @{=>}[dl] & \eqref{tfae1}\ar @{=>}
[r] & \eqref{tfae2}\ar @{=>}[d]\\
\eqref{tfae6} \ar @{=>}[r]  & \eqref{tfae7}\ar @{=>}[u] &
\eqref{tfae4} \ar @{=>}[l] & \eqref{tfae3} \ar @{=>}[l]}
The main result of \cite{bumagin:rh} is that \eqref{tfae1} implies
\eqref{tfae2} (cf. Remark \ref{r:bcp}).
By Proposition 
\ref{p:conelii}, \eqref{tfae2} implies \eqref{tfae3}.

By Lemma \ref{l:combtohom}, \eqref{tfae3} implies 
\eqref{tfae4} and \eqref{tfae6} implies \eqref{tfae7}.

By Theorem \ref{t:Xlii}, \eqref{tfae4} implies 
\eqref{tfae7}.

Now, \eqref{tfae7} implies \eqref{tfae5}, by Theorem \ref{t:homtohyp}.

Proposition
\ref{p:hyptocomb} gives that \eqref{tfae5} implies \eqref{tfae6},
provided we can find a bound on the area of fillings of short loops in
$X$ (where ``short'' means at most $16\delta$).   Any combinatorial
loop in $X$ of length at most
$16\delta$ which lies in a horoball can be filled by a combinatorial
disk of length $48\delta$, by Proposition \ref{p:combhorolii}.  Up to
the $G$-action, there are only finitely many loops of length less than
$16\delta$ which do not lie in a horoball.  Since $X$ is simply
connected, these can all be filled, and so there is some universal
constant $C(X)$ so that any loop in $X$ of length less than $16\delta$
can be filled by a disk of area at most $C(X)$.

It remains to observe that the space $X$ satisfies the conditions of
Definition \ref{d:gromovrh}, and so \eqref{tfae5} implies \eqref{tfae1}.

Gathering together these implications, the theorem
is proved.
\end{proof}

\subsection{Metric properties of $X$}

We now suppose that $G$ is hyperbolic relative to $\P$,
and that $X = X(G,\P,S)$ as in Definition \ref{d:X1} (for the
moment we are only concerned with the metric properties
of $X$, and we restrict our attention to the $1$-complex version
of $X$).

By Theorem \ref{t:tfae}, the graph $X$ is $\delta$-hyperbolic for some
$\delta$.  We assume that this $\delta$ satisfies the conditions
of Remark \ref{r:deltaprops}

\begin{lemma}\label{lemma:convex}
If $L>\delta$, the $L$-horoballs are convex in $\Xgraph$. 
\end{lemma}
\begin{proof}
Let $H_1$ be a $1$-horoball, and let $H_L$ be the $L$-horoball
contained in $H_1$.
We
observe that $H_{L}$ is convex in $H_1$ (where $H_1$ is endowed with
its path metric).  Thus if $H_{L}$ fails to
be convex in $\Xgraph$, then two points in $H_L$ are connected by a
geodesic which passes through $D^{-1}(0)$.  Let $p$ and $q$ be two such points
in $H_{L}$, chosen to have minimal distance from one another in
the path metric on $H_1$.  (Note that this distance must be at least
$2 L$ so that the geodesic between them actually leaves $H_1$. They
must also satisfy $D(p)=D(q)=L$.)
Choose another point $r\in H_{L}$ so that 
$\max\{d_{H_1}(p,r),d_{H_1}(q,r)\}< d_{H_1}(p,q)$ and $D(r)=L$.
There are then $X$-geodesics $[p,r]$ and
$[q,r]$ which lie entirely in $H_{L}$, whereas $p$ and $q$ are
joined by a geodesic which includes points in $D^{-1}(0)$, i.e., points
which are of distance at least $L$ from $H_{L}$.  Since the
triangle formed by these geodesics is $\delta$-slim, this implies that
$L\leq\delta$, a contradiction.
\end{proof}

We extend the choice of geodesics in Lemma \ref{l:gamma} 
to a choice of
geodesics between any two points in $X$:
\begin{lemma}\label{l:geodesicbicombing}
If $G$ is torsion-free then there is an antisymmetric, $G$-equivariant geodesic bicombing $\geod$ on $X$, so
that if $x$ and $y$ lie in the same $L$-horoball for $L>2\delta$, then
$\geod(x,y)$ is as described in Lemma \ref{l:gamma}.
\end{lemma}
\begin{proof}
Choose a complete set $\{ o_i \}$ of representatives for the orbits of
vertices of $X$ under the $G$-action.  For each $o_i$ and each
$x \in X$, choose a geodesic $\geod(o_i,x)$ as in Lemma \ref{l:gamma}.
Extend equivariantly and antisymmetrically.
Because $G$ is torsion free, there is no element
of $G$ which exchanges two points of $X$.  Thus equivariance and
antisymmetry may happily coexist.
\end{proof}

\begin{remark}
We say the bicombing in \ref{l:geodesicbicombing} is
`antisymmetric' because we consider paths to be maps.  If
paths are considered as subsets, it would be symmetric.  We
will sometimes blur this distinction, but it will not introduce
any confusion.  For more about the parametrizations of 
the bicombing, see Paragraph \ref{ss:param}.

When we define the homological bicombing in Section 
\ref{section:bicombing}, it will be {\em antisymmetric}, of
course.
\end{remark}

\begin{torsionremark} \label{r:tf1}
If $G$ is not torsion-free, then there may be no bicombing
as in Lemma \ref{l:geodesicbicombing}, because of the possible
presence of $2$-torsion.  The way to proceed in the
presence of torsion is to let $\gamma(x,y)$ be the collection of
all
of the geodesics between $x$ and $y$.  This makes a number
of the arguments later in this paper more awkward.  At times it
is also useful to consider the ``average" of all of the geodesics
between $x$ and $y$.
\end{torsionremark}

The following is a slight generalization of the usual notion of
cone types.
\begin{definition}\label{d:conetype} [Cone types]
Suppose that $G$ is a group, that $\Xi$ is a graph equipped
with a free $G$-action and that $x \in \Xi$ is a vertex.  

For a combinatorial path $\gamma \co I \to \Xi$, let $[\gamma]$ denote
the $G$-orbit of $\gamma$.  If $v$ is a vertex in $\Xi$, the {\em cone type of $v$ viewed from $x$} is the collection of classes $[\gamma]$ of paths
for which 
\begin{enumerate}
\item there exists $g \in G$ so that $g \cdot \gamma(0) = v$; and
\item $d(x, g \cdot \gamma(1)) = d(x,v) + |\gamma|$.
\end{enumerate}
\end{definition}
Note that if $v_1$ and $v_2$ have the same cone type then, in 
particular, $v_1$ and $v_2$ lie in the same $G$-orbit.

In case $\Xi$ is the Cayley graph of a finitely generated group $G$, and
$x = 1$, the above definition is equivalent to the usual notion of 
cone types (see \cite[III.$\Gamma$.2.16]{bridhaef:book}, for example).

The following result follows directly from the proof of \cite[Theorem
III.$\Gamma$.2.18]{bridhaef:book}.

\begin{lemma} \label{l:finmanycones}
Let $G$ be a finitely generated group acting freely on the locally
compact $\delta$-hyperbolic graph $\Xi$, and suppose that $x \in \Xi$
is a vertex.  
Each orbit of vertices in $\Xi$ contains only finitely many cone types
viewed from $x$.
\end{lemma}

\begin{lemma} \label{l:thicklocal}
Suppose that $G$ is hyperbolic relative to $\mc{P} = \{ P_1, \ldots , P_n \}$, where no $P_i$ is equal to $G$.  Suppose further that $S$ is a 
compatible generating set for $G$ with respect to $\mc{P}$, and
that $X = X(G,\mc{P},S)$ is $\delta$-hyperbolic.  Finally, suppose
that $\mathrm{diam}(\Gamma(P_i,S\cap P_i)) \geq 2^{15\delta +1}$ for each $i$.

For each $L > 0$ there exists a $10\delta$-local geodesic $\gamma 
\subset X$ of length at least $L$ so that $\gamma$ does not intersect
any $(15\delta+1)$-horoball in $X$.
\end{lemma}
\begin{proof}
Proceed as follows.  Let $H_1$ and $H_2$ be distinct  
$10\delta$-horoballs in $X$.  Let $\nu$ be a shortest path
between $H_1$ and $H_2$.  Note that the length of
$\nu$ is at least $20\delta$.

Let $x = \nu \cap H_1$.  For some $i$, some $t \in T_i$ and some
$p \in P_i$, we have $x = (i,t,p,10\delta)$.  By hypothesis, there exists
$q \in P_i$ so that the distance in $\Gamma(P_i,S \cap P_i)$ between
$p$ and $q$ is exactly $2^{15\delta +1}$.  This implies that the
geodesic $\geod(x,y)$ between $x$ and $y = (i,t,q,10\delta)$ intersects
$D^{-1}(15\delta)$ but not $D^{-1}(15\delta +1)$.  This geodesic
consists of two vertical segments, and a single horizontal segment
of length $2$.

Let $\nu' = qp^{-1} \cdot \nu$.  Our $10\delta$-local geodesic segment
is then constructed as follows:  Start with $\nu$, followed by 
$\geod(x,y)$.  Then take $\nu'$ to the $10\delta$-horoball 
$H_3 = qp^{-1} H_2$.  From the endpoint $x'$ of $\nu'$ in $H_3$, 
construct a a path $\geod(x',y')$ exactly as above (See Figure
\ref{f:localgeodesic}).  
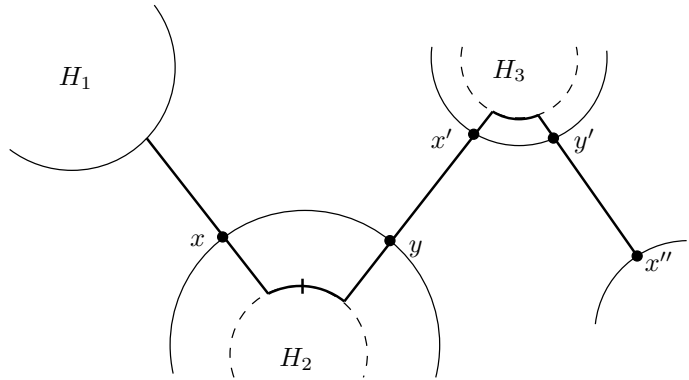
\begin{figure}[htbp]
\begin{center}
\begin{picture}(0,0)%
\includegraphics{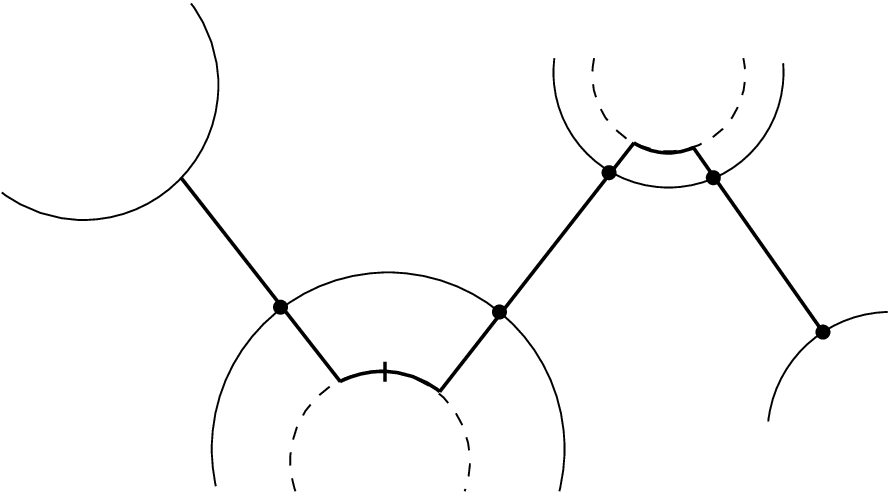}%
\end{picture}%
\setlength{\unitlength}{4144sp}%
\begingroup\makeatletter\ifx\SetFigFont\undefined%
\gdef\SetFigFont#1#2#3#4#5{%
  \reset@font\fontsize{#1}{#2pt}%
  \fontfamily{#3}\fontseries{#4}\fontshape{#5}%
  \selectfont}%
\fi\endgroup%
\begin{picture}(4066,2245)(1613,-3288)
\put(3236,-3212){\makebox(0,0)[lb]{\smash{{\SetFigFont{10}{12.0}{\familydefault}{\mddefault}{\updefault}{\color[rgb]{0,0,0}$H_2$}%
}}}}
\put(4010,-2530){\makebox(0,0)[lb]{\smash{{\SetFigFont{10}{12.0}{\familydefault}{\mddefault}{\updefault}{\color[rgb]{0,0,0}$y$}%
}}}}
\put(4511,-1483){\makebox(0,0)[lb]{\smash{{\SetFigFont{10}{12.0}{\familydefault}{\mddefault}{\updefault}{\color[rgb]{0,0,0}$H_3$}%
}}}}
\put(5421,-2643){\makebox(0,0)[lb]{\smash{{\SetFigFont{10}{12.0}{\familydefault}{\mddefault}{\updefault}{\color[rgb]{0,0,0}$x''$}%
}}}}
\put(1917,-1529){\makebox(0,0)[lb]{\smash{{\SetFigFont{10}{12.0}{\familydefault}{\mddefault}{\updefault}{\color[rgb]{0,0,0}$H_1$}%
}}}}
\put(2701,-2491){\makebox(0,0)[lb]{\smash{{\SetFigFont{10}{12.0}{\familydefault}{\mddefault}{\updefault}{\color[rgb]{0,0,0}$x$}%
}}}}
\put(4141,-1906){\makebox(0,0)[lb]{\smash{{\SetFigFont{10}{12.0}{\familydefault}{\mddefault}{\updefault}{\color[rgb]{0,0,0}$x'$}%
}}}}
\put(4996,-1906){\makebox(0,0)[lb]{\smash{{\SetFigFont{10}{12.0}{\familydefault}{\mddefault}{\updefault}{\color[rgb]{0,0,0}$y'$}%
}}}}
\end{picture}%
\caption{How to make an arbitrarily long $10\delta$-local geodesic
  which stays in the ``thick'' part of $X$.}
\label{f:localgeodesic}
\end{center}
\end{figure}
Continue in this
manner.  This construction can be iterated as many times as necessary
to ensure the path has at least length $L$, and it is not difficult
to see that the ensuing path is a $10\delta$-local geodesic (note that the path $\nu$, and its translates, must begin and end with a vertical path of
length $10\delta$, since it is a shortest path between horoballs).
\end{proof}

\begin{theorem} \label{t:onethickZ}
Suppose that $G$ is hyperbolic relative to $\mc{P} = \{ P_1, \ldots , P_n \}$, where no $P_i$ is equal to $G$.  Suppose further that $S$ is a 
compatible generating set for $G$ with respect to $\mc{P}$, and
that $X = X(G,\mc{P},S)$ is $\delta$-hyperbolic.  Finally, suppose
that $\mathrm{diam}(\Gamma(P_i,S\cap P_i)) \geq 2^{15\delta +1}$ for each $i$.

Then there exists a hyperbolic element $g \in G$ which has an axis
$\gamma \in X$ so that $\gamma \subset D^{-1}[0,19\delta]$.
\end{theorem}
\begin{proof}
Let
$K_{19\delta}$ be the total number of cone types of all points in
$D^{-1}[0,19\delta]$ as viewed from $1 \in X$.

Let $\sigma$ be a $10\delta$-local geodesic of length longer than
$\frac{7}{3}(K_{19\delta}) + 2\delta$ as in Lemma \ref{l:thicklocal}.  In
particular, $\sigma$ does not penetrate any $(15\delta+1)$-horoball.
By translating $\sigma$ by an element of $G$, we suppose that
$\sigma$ begins at $1 \in X$.

By \cite[Theorem III.H.1.13(1)]{bridhaef:book}, 
any $10\delta$-local geodesic is contained in the $2\delta$-neighborhood of
any geodesic joining the endpoints.  A simple argument then shows
that in fact the geodesic is contained in the $4\delta$-neighborhood
of the $k$-local geodesic.  Also, any $10\delta$-local 
geodesic is a $(\frac{7}{3},2\delta)$-quasi-geodesic (by
\cite[Theorem III.H.1.13(3)]{bridhaef:book}).

Let $\rho$ be a geodesic segment joining the endpoints of $\sigma$.
We have ensured that the length of $\rho$ is greater than 
$K_{19\delta}$, that $\rho$ does not intersect any 
$(19\delta +1)$-horoball, and that $\rho$ starts at $1$.
\footnote{The following argument is very similar to that of
\cite[Proposition III.$\Gamma$.2.22]{bridhaef:book}.}

Therefore, there exist vertices $v_1, v_2 \in \rho$ which have the
same cone type as viewed from $1$.  Suppose that $v_1$ occurs before $v_2$ on
$\rho$.  Let $\rho_1$ be the subpath of $\rho$ from the
$1$ to $v_1$, let $\rho_2$ be that part between $v_1$ and
$v_2$, and let $\rho_3$ be the remainder of $\rho$.

Since $v_1$ and $v_2$ have the same cone type, they are in the 
same $G$-orbit.  Let $g \in G$ be so that $g.v_1 = v_2$.  Since
$\rho_1\rho_2\rho_3$ is a geodesic, and $v_1$ and $v_2$ have
the same cone type as viewed from $1$, the path
$\rho_1 \rho_2 (g\cdot \rho_2) (g \cdot \rho_3)$ is also a geodesic.

In turn, this implies that $\rho_1\rho_2(g \cdot \rho_2)(g^2\cdot \rho_2)(g^2 \cdot \rho_3)$ is a geodesic.

Iterating this argument, we see that the path
\[	\rho_2 (g \cdot \rho_2) (g^2 \cdot \rho_2) \cdots	,	\]
is a geodesic ray, starting at $v_1$.  Denote this ray by $r$.
Note that $r \subset g^{-1} \cdot r$.  Let $\gamma$ be the union of
the paths $g^{-i} \cdot r$ as $i \to \infty$, parametrized in the obvious
way.  This is a bi-infinite geodesic line, contained in $D^{-1}[0,19\delta]$,
upon which $g$ acts by translation.
\end{proof}

\begin{theorem} \label{t:independentZ}
Suppose that $G$ is hyperbolic relative to $\mc{P} = \{ P_1, \ldots , P_n \}$, that $G \neq P_i$ for each $i$ and that $P_1$ is infinite.  Then
$G$ is non-elementary relatively hyperbolic.
\end{theorem}
\begin{proof}
Suppose first that some parabolics are finite.  Let $X'$ be the space 
obtained from Definition \ref{d:X1} using {\em all} of the parabolics, and 
$X$ the space obtained using only the infinite parabolics.  Up to
quasi-isometry, $X'$ is $X$ with a locally finite collection of rays 
attached to the cosets of the finite parabolics. In particular, $X$ is
Gromov hyperbolic if and only if $X'$ is.  Moreover, an element $g \in G$ 
acts hyperbolically on $X'$ if and only if it acts hyperbolically on $X$.

Therefore, we are free to assume that all of the parabolic subgroups
of $G$ are infinite.

Let $\gamma$ be the geodesic from Theorem \ref{t:onethickZ},
which is an axis of a hyperbolic element $g \in G$.  Let 
$\gamma^+$ and $\gamma^-$ be the points in $\partial X$ at
either end of $\gamma$.

Let $H$ be a $25\delta$-horoball in $X$, and let $\e_H$ be
the point in $\partial X$ coming from $H$.  Consider an
ideal geodesic triangle, $T$, with vertices $\gamma^-, \gamma^+$ and
$\e_H$, and edge $\gamma$ between $\gamma^-$ and $\gamma^+$.
Suppose furthermore, that the geodesics with endpoint $\e_H$
are vertical after depth $2\delta$.
By Lemma \ref{l:idealthin}, this triangle is $3\delta$-slim.

The triangle $T$ intersects $H \cap D^{-1}(25\delta)$ in a pair of 
points $\{ x^+, x^- \}$ which are at most $3\delta$ apart.  Since the parabolic $P$
which stabilizes $H$ is infinite, there exists $p \in P$ so that
$d_X(\{ p \cdot x^-, p \cdot x^+ \}, \{ x^-, x^+ \}) \ge 10 \delta$.

Now, $pT$ is another ideal triangle, with one vertex $\e_H$ and the
opposite side an axis for $pgp^{-1}$.  We claim that
$\mbox{Fix}_{\partial X}(g) \cap \mbox{Fix}_{\partial X}(pgp^{-1}) =
\emptyset$.

An easy
argument shows that the Hausdorff distance between two geodesics
with the same endpoints at infinity is at most $2\delta$.

Suppose, for instance that $p \cdot \gamma^+ = \gamma^-$.  Then
the fact that the geodesics from $p \cdot \gamma^+$ to $\e_H$
and from $\gamma^-$ to $\e_H$ are vertical below the $2\delta$
level of $H$, easily implies that $d_X(p \cdot x^+, x^- )
\le 4\delta$, in contradiction to the choice of $p$.

Thus, the elements $g$ and $pgp^{-1}$ are hyperbolic elements with
disjoint fixed sets in $\partial X$, and $G$ is non-elementary
relatively hyperbolic, as required.
\end{proof}

\subsection{Constants}
By virtue of Theorem \ref{t:tfae}, we may assume that
$X(G,\mc{P},S)$ is $\delta$-hyperbolic for some $\delta$; we are free
to assume that $\delta$ is an integer and that $\delta\geq 100$.  It is
useful to have a few different scales to work at, and so we choose the
following constants.  The choices made will be justified by the
results of the subsequent sections.

\begin{definition}\label{d:constants}
We set $K=10\delta$, $L_1=100 K$, and $L_2=3 L_1$.
\end{definition}

We have made no attempt to make this choice of constants
optimal.

\def\GHH {\mathcal C}  
\def\CHH {{\mathcal D}^\infty}  
\def\NHH {\mathcal C^0}  
\def\PHH {{\mathcal C}^K} 

\newcommand{\Prop}[4]{\mathcal{B}^{#1}({#2},{#3};{#4})}
\newcommand{\Sp}[3]{\geod^{#3}({#1},{#2})}

\def\Lone {100K} 

 \section{Convex sets and betweenness}
\label{s:horoballs}

In this section we prove a theorem about 
collections of convex sets in an arbitrary proper
Gromov hyperbolic space.  The construction in this section is crucial to 
the construction
of {\em preferred paths} in the next section.  In turn, preferred
paths are
the key to our construction of the bicombing in Section 
\ref{section:bicombing}.

\begin{definition} \label{d:globule}
Let $\Upsilon$ be a geodesic metric space.
A collection $\mc{G}$ of convex sets is
{\em $N$-separated} if
for all $A,B \in\mc{G}$ the distance between $A$ and
$B$ is at least $N$.

The elements of $\mc{G}$ will be called {\em globules}.

Let $\text{Isom}(\Upsilon;\mc{G})$ be the collection of all isometries
$g$ of $\Upsilon$ so that for all $A\in \mc{G}$ we have $gA \in \mc{G}$.
\end{definition}

Throughout the remainder of this section, we will suppose that
$\mc{G}$ is a $50\delta$-separated collection of convex subsets
of a $\delta$-hyperbolic space $\Upsilon$.

\begin{remark}
In subsequent sections of this paper, we will apply the results
of this section in case $\Upsilon = X(G,\mc{P},S,\mc{R})$
and $\mc{G}$ is the collection of all $L_1$-horoballs in $X$.

Another interesting example is $\mathbb{H}^n$ with some collection
of (sufficiently separated) horoballs.  It is worth remarking that the existence of a function $\CHH$ satisfying
the properties (A1)--(A7) below is not obvious even in the case
of horoballs in $\H^n$.
\end{remark}

The main purpose of this section is to define, for any
pair of points $a,b \in \Upsilon$, a family of globules
$\CHH_{a,b}$, which will be the collection of globules
which are `between' $a$ and $b$ (in case $a$ is contained 
in a globule $A$ we will have $A \in \CHH_{a,b}$ for all $b$).
We want our collections $\CHH_{a,b}$
to satisfy various conditions, listed as Axioms (A1)-(A7)
below.  Most of these are quite straightforward to ensure
but Axioms (A5)--(A7), the most important for the applications,
are much more difficult to guarantee.

\begin{remark}
The construction performed in this section can be done
if the globules are only quasi-convex.  How far the globules must be 
separated depends on the constants of quasi-convexity.  
\end{remark}

The following two lemmas will be useful:
\begin{lemma}
If $C$ is a convex subset of $\Upsilon$ then for any $R\ge 0$,
the set $N_R(C)$ is $\delta$-quasi-convex.
\end{lemma}
\begin{lemma}\label{l:tube}
Let $A$, $B$ be disjoint closed $\delta$-quasi-convex sets in the proper $\delta$-hyperbolic space
$\Upsilon$.  Let $a$, $a'\in A$, and $b$, $b'\in B$.  If $q$ is a geodesic
between $a$ and $b$, and $q'$ a geodesic between $a'$ and $b'$, then
$q'$ lies in the $3\delta$-neighborhood of $A\cup B\cup q$.

Moreover, for any $4\delta<R<\frac{1}{2}d(A,B)$, the Hausdorff
distance between the (nonempty) sets $q\smallsetminus 
N_R(A\cup B)$ and $q'\smallsetminus N_R(A\cup B)$ is
at most $5\delta$.
\end{lemma}
\begin{proof}
Let $[a,a']$ be any geodesic between $a$ and $a'$, and let $[b,b']$ be
any geodesic between $b$ and $b'$;  since $A$ and $B$ are $\delta$-quasi-convex,
$[a,a']\subset N_\delta(A)$ and $[b,b']\subset N_\delta(B)$.
Consider the quadrilateral $[a,a']\cup q\cup[b,b']\cup[b,a']$, 
pictured in Figure \ref{f:tubelemma}.  
\begin{figure}[htbp]
\begin{center}
\begin{picture}(0,0)%
\includegraphics{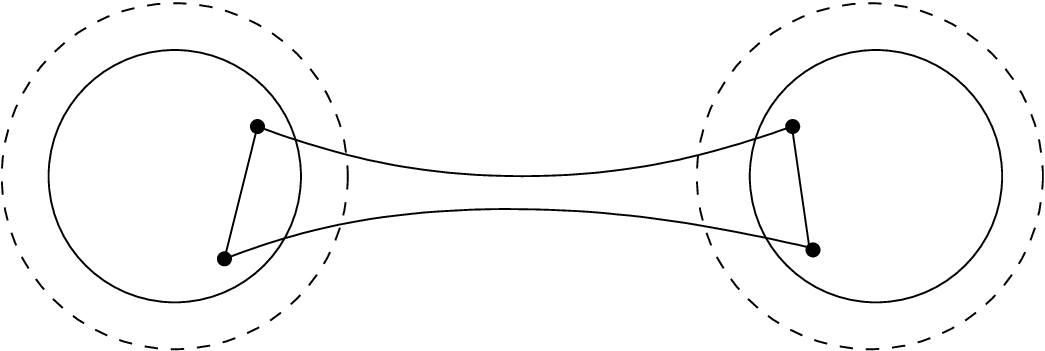}%
\end{picture}%
\setlength{\unitlength}{4144sp}%
\begingroup\makeatletter\ifx\SetFigFont\undefined%
\gdef\SetFigFont#1#2#3#4#5{%
  \reset@font\fontsize{#1}{#2pt}%
  \fontfamily{#3}\fontseries{#4}\fontshape{#5}%
  \selectfont}%
\fi\endgroup%
\begin{picture}(4776,1598)(2284,-3898)
\put(3260,-2869){\makebox(0,0)[lb]{\smash{{\SetFigFont{12}{14.4}{\familydefault}{\mddefault}{\updefault}{\color[rgb]{0,0,0}$a$}%
}}}}
\put(6005,-2869){\makebox(0,0)[lb]{\smash{{\SetFigFont{12}{14.4}{\familydefault}{\mddefault}{\updefault}{\color[rgb]{0,0,0}$b$}%
}}}}
\put(3120,-3479){\makebox(0,0)[lb]{\smash{{\SetFigFont{12}{14.4}{\familydefault}{\mddefault}{\updefault}{\color[rgb]{0,0,0}$a'$}%
}}}}
\put(6094,-3441){\makebox(0,0)[lb]{\smash{{\SetFigFont{12}{14.4}{\familydefault}{\mddefault}{\updefault}{\color[rgb]{0,0,0}$b'$}%
}}}}
\put(2624,-3136){\makebox(0,0)[lb]{\smash{{\SetFigFont{12}{14.4}{\familydefault}{\mddefault}{\updefault}{\color[rgb]{0,0,0}$A$}%
}}}}
\put(6450,-3085){\makebox(0,0)[lb]{\smash{{\SetFigFont{12}{14.4}{\familydefault}{\mddefault}{\updefault}{\color[rgb]{0,0,0}$B$}%
}}}}
\put(4588,-3009){\makebox(0,0)[lb]{\smash{{\SetFigFont{12}{14.4}{\familydefault}{\mddefault}{\updefault}{\color[rgb]{0,0,0}$q$}%
}}}}
\put(4556,-3504){\makebox(0,0)[lb]{\smash{{\SetFigFont{12}{14.4}{\familydefault}{\mddefault}{\updefault}{\color[rgb]{0,0,0}$q'$}%
}}}}
\end{picture}%
\caption{Tube lemma.}
\label{f:tubelemma}
\end{center}
\end{figure}
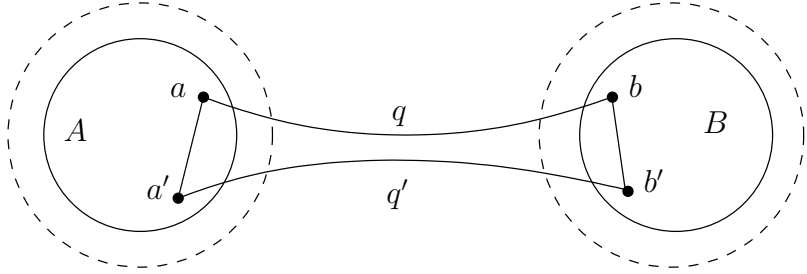
Since $\Upsilon$ is $\delta$-hyperbolic, this quadrilateral is
$2\delta$-slim.  In particular, if $x$ is any point on $q'$, then 
\[x\in N_{2\delta}\left( [a,a']\cup q\cup[b,b'] \right) \subset 
N_{3\delta} \left( A \cup q\cup B \right).\]

For the second assertion of the lemma,
let $R$ be as in the statement;  since $R<\frac{1}{2}d(A,B)$,
the set $q\smallsetminus N_R(A\cup B)$ is non-empty.
Let $x\in q\smallsetminus N_R(A\cup B)$.  We must show that 
\[d(x,q'\smallsetminus N_R(A\cup B))\leq 5\delta.\]

As above, there is a point
$x'\in[a,a']\cup q'\cup[b,b']$ so that $d(x,x')\leq 2\delta$.  Since
$R>3\delta$, and $x$ is at least $R$ from $A\cup B$, the
point $x'$ must be on $q'$.  If $x'$ is outside the $R$-neighborhood
of $A\cup B$, then we are done, so suppose without loss of generality
that $x'$ lies in $N_R(A)$.

Let $y$ be the point on $q'$ which is
distance exactly $R$ from $A$.
Let $t$ be a point
in $A$ with $d(y,t)=d(y,A)=R$.  Consider the 
geodesic triangle pictured in
Figure \ref{f:tubelemma2}
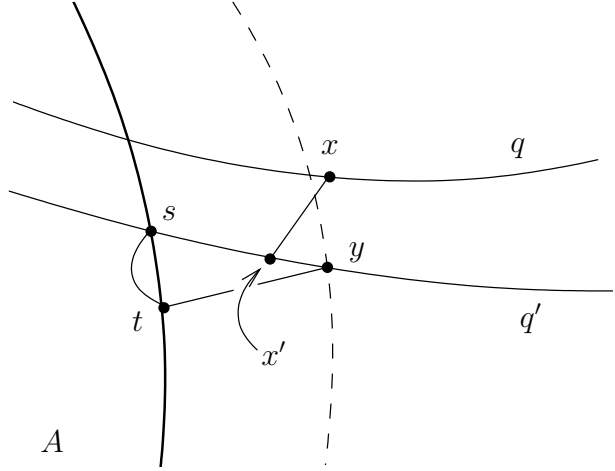
\begin{figure}[htbp]
\begin{center}
\begin{picture}(0,0)%
\includegraphics{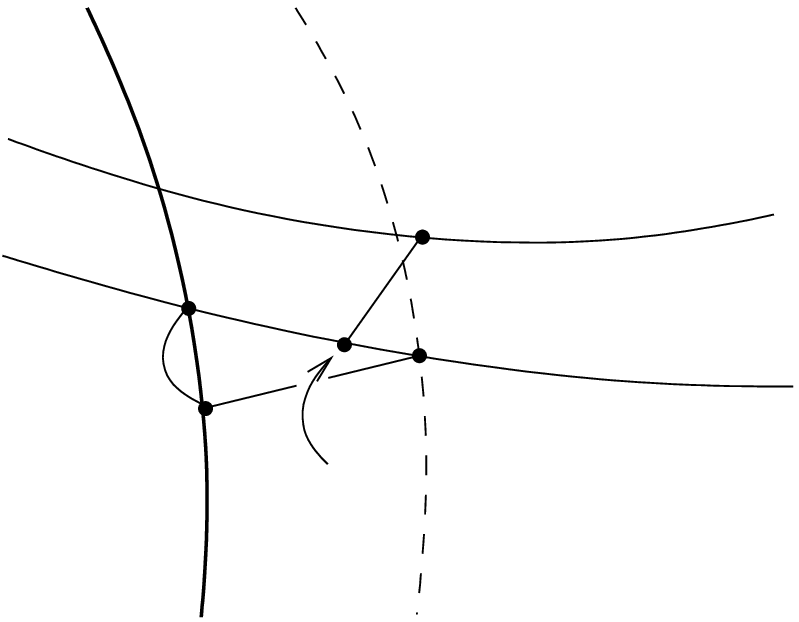}%
\end{picture}%
\setlength{\unitlength}{4144sp}%
\begingroup\makeatletter\ifx\SetFigFont\undefined%
\gdef\SetFigFont#1#2#3#4#5{%
  \reset@font\fontsize{#1}{#2pt}%
  \fontfamily{#3}\fontseries{#4}\fontshape{#5}%
  \selectfont}%
\fi\endgroup%
\begin{picture}(3638,2828)(2390,-4885)
\put(5458,-4037){\makebox(0,0)[lb]{\smash{{\SetFigFont{12}{14.4}{\familydefault}{\mddefault}{\updefault}{\color[rgb]{0,0,0}$q'$}%
}}}}
\put(5401,-2996){\makebox(0,0)[lb]{\smash{{\SetFigFont{12}{14.4}{\familydefault}{\mddefault}{\updefault}{\color[rgb]{0,0,0}$q$}%
}}}}
\put(2592,-4787){\makebox(0,0)[lb]{\smash{{\SetFigFont{12}{14.4}{\familydefault}{\mddefault}{\updefault}{\color[rgb]{0,0,0}$A$}%
}}}}
\put(3317,-3401){\makebox(0,0)[lb]{\smash{{\SetFigFont{12}{14.4}{\familydefault}{\mddefault}{\updefault}{\color[rgb]{0,0,0}$s$}%
}}}}
\put(3139,-4062){\makebox(0,0)[lb]{\smash{{\SetFigFont{12}{14.4}{\familydefault}{\mddefault}{\updefault}{\color[rgb]{0,0,0}$t$}%
}}}}
\put(4270,-3007){\makebox(0,0)[lb]{\smash{{\SetFigFont{12}{14.4}{\familydefault}{\mddefault}{\updefault}{\color[rgb]{0,0,0}$x$}%
}}}}
\put(4435,-3618){\makebox(0,0)[lb]{\smash{{\SetFigFont{12}{14.4}{\familydefault}{\mddefault}{\updefault}{\color[rgb]{0,0,0}$y$}%
}}}}
\put(3914,-4266){\makebox(0,0)[lb]{\smash{{\SetFigFont{12}{14.4}{\familydefault}{\mddefault}{\updefault}{\color[rgb]{0,0,0}$x'$}%
}}}}
\end{picture}%
\caption{Tube lemma.}
\label{f:tubelemma2}
\end{center}
\end{figure}
with vertices $s$, $t$, and $y$.  This triangle is $\delta$-thin and
the geodesic between $s$ and $t$ is contained in $N_\delta(A)$.
Since $d(x,A)\geq R$, we have
$d(x',A)\geq R-2\delta >2\delta$, so $x'$ cannot be as close as
$\delta$ to the geodesic between $s$ and $t$; thus 
there is a point $z$ on the geodesic between $y$ and $t$ so that
$d(x',z)\leq \delta$ and $d(z,y)=d(x',y)$.  We have
\[d(z,A)\geq d(x',A)-\delta \geq R-3\delta,\]  
and so $d(z,y)\leq 3\delta$; this implies $d(x',y)=d(z,y)\leq 3\delta$.
Finally, $d(x,y)\leq d(x,x')+d(x',y)\leq 5\delta$, and the lemma is proved.
\end{proof}

\begin{remark}
The assumptions of Lemma \ref{l:tube} can be weakened.  At the cost of a slightly more complicated proof, the assumptions of `closed' and 'proper' can be removed.

However we will only ever use Lemma \ref{l:tube} exactly as it is stated.
\end{remark}

\begin{definition} \label{d:geods}
For $a,b \in \Upsilon$, let $\alpha(a,b)$ be the set of all
geodesics between $a$ and $b$.
\end{definition}

\begin{definition}
Suppose that $a,b \in\Upsilon$.  Let $\NHH (a,b)$ denote the set
of globules which intersect some element of $\alpha(a,b)$ nontrivially.

For any $R$ let $\mc{C}^R(a,b)$ denote the set of globules
$P$ so that for every $\gamma\in \alpha(a,b)$
\[	N_R(P) \cap \gamma \neq \emptyset.	\]
\end{definition}
Since bigons are $\delta$-thin, $\NHH(a,b) \subseteq \PHH(a,b)$
for all $a,b\in \Upsilon$. (Recall $K=10\delta$.)
\begin{remark}\label{r:orientation}
For a geodesic $\gamma$ in $\Upsilon$, the set of globules
$P \in \mc{G}$ so that $N_{2K}(P) \cap \gamma \neq \emptyset$
inherits a natural
order by projection to the geodesic.  This is because globules
lie at least $50\delta=5K$ apart from each other.  This induces an
order on $\mc{C}^{2K}(a,b)$ for each $a,b \in \Upsilon$ and
each $\gamma \in \alpha(a,b)$;
this order is independent of $\gamma$.

Note that if $\mc{C}^{2K}(a,b)$
intersects $\mc{C}^{2K}(c,d)$ nontrivially, then the order on the
intersection inherited from $\mc{C}^{2K}(c,d)$ will either coincide
with the order inherited from $\mc{C}^{2K}(a,b)$ or with its reverse.
\end{remark}

For any $a,b \in \Upsilon$, there are only finitely many elements of $\PHH (a,b)$:
the size is bounded by $\frac{1}{K} d(a,b) + 1$.  

Let $O_{\mc{G}}$ be the set of totally ordered finite subsets of $\mc{G}$

Thus we have functions
\[\NHH\co \Upsilon\times\Upsilon\to O_{\mc{G}} , \]
and 
\[\PHH\co \Upsilon\times\Upsilon\to O_{\mc{G}},\]
Note that the action of $\text{Isom}(\Upsilon;\mc{G})$ on $\mc{G}$ 
induces an action on $O_{\mc{G}}$ and that $\NHH$ and $\PHH$ are
equivariant with respect to this action.
(Note also that the same subset appears many times in
$O_{\mc{G}}$, once 
for each possible total order on the subset.  Thanks to Remark
\ref{r:orientation}, only two of these orders actually appear in the
image of $\NHH$ or $\PHH$.)

For a function
\[\GHH\co \Upsilon\times\Upsilon\to O_{\mc{G}},\]
we will denote $\GHH(a,b)$ by $\GHH_{a,b}$.
In the properties (A1)--(A6) defined below, the symbol "$\subset$" is 
used to denote 
``ordered subset'', not just ``subset''.  

We will use the following notation for subintervals: Let $a,b \in 
\Upsilon$.  If $A$ and $B$ are contained in $\GHH_{a,b}$ then
\begin{eqnarray*}
\GHH_{a,b}[A,B] &=& \{	Y \in \GHH_{a,b} \mid A\le Y \le B \},\\
\GHH_{a,b}(a,A] &=&     \{ Y \in \GHH_{a,b} \mid Y \le A \},\mbox{ and }\\
\GHH_{a,b}[A,b) &=& \{ Y \in \GHH_{a,b} \mid A \le Y \}.
\end{eqnarray*}

Here are some useful conditions $\GHH$ might satisfy:

\begin{enumerate}
 \item[(A1)]  For all $a,b \in \Upsilon$, $\mc{C}^0_{a,b} \subseteq \mc{C}_{a,b}$;
 \item[(A2)]For all $a,b \in \Upsilon$, $\mc{C}_{a,b} \subseteq \mc{C}^K_{a,b}$; 
 \item[(A3)] For all $a,b \in \Upsilon$, $\mc{C}_{a,b}  = (\mc{C}_{b,a})^{op}$;
 \item[(A4)] $\GHH$ is $\text{Isom}(\Upsilon;\mc{G})$-equivariant. 
 \item[(A5)]If $A,B \in \mc{C}_{a,b} \cap \mc{C}_{c,d}$ for some $a,b,c,d \in \Upsilon$
and $A,B \in \mc{G}$ then $\mc{C}_{a,b}[A,B] = \mc{C}_{c,d}[A,B]$; 
\item[(A6)] If $A \in \mc{C}_{a,b} \cap \mc{C}_{a,c}$ for some $a,b,c \in \Upsilon$ and
$A \in \mc{G}$ then $\mc{C}_{a,b}(a,A] = \mc{C}_{a,c}(a,A]$.
\item[(A7)] If $A \in \mc{C}_{a,b} \cap \mc{C}_{c,b}$ for some 
$a,b,c \in\Upsilon$ and $A \in \mc{G}$ then 
$\mc{C}_{a,b}[A,b) = \mc{C}_{c,b}[A,b)$.
\end{enumerate}

Note that Axiom (A7) follows from (A6) and (A3).

The purpose of this section is to find a function $\CHH$ which
satisfies all seven of the axioms (A1)--(A7).  For the application in the
next section, (A5)--(A7) are the most important.  They are also the most difficult
to ensure, although our definition of $\CHH$ in Definition 
\ref{d:filtration} below is designed to make (A5)--(A7)
as apparent as possible.

The approach to constructing $\CHH$ is as follows.
We start with $\NHH$ with the order as in Remark \ref{r:orientation}
and observe:
\begin{lemma} 
The function $\NHH$ satisfies (A1) -- (A4)
\end{lemma}

For a pair $a,b \in \Upsilon$, the final $\CHH_{a,b}$ will come from 
the preliminary $\NHH_{a,b}$ by adding new elements in order 
to make (A5)--(A7) hold.  It is not at all obvious that enforcing 
(A5)--(A7) whilst
retaining (A1)--(A4) is possible.  This is the content of the proof of 
Theorem \ref{t:ppw} below.

By (A2), all of the elements of $\CHH_{a,b}$ must lie in $\PHH_{a,b}$.  

\begin{figure}[htbp]
\begin{center}
\begin{picture}(0,0)%
\includegraphics{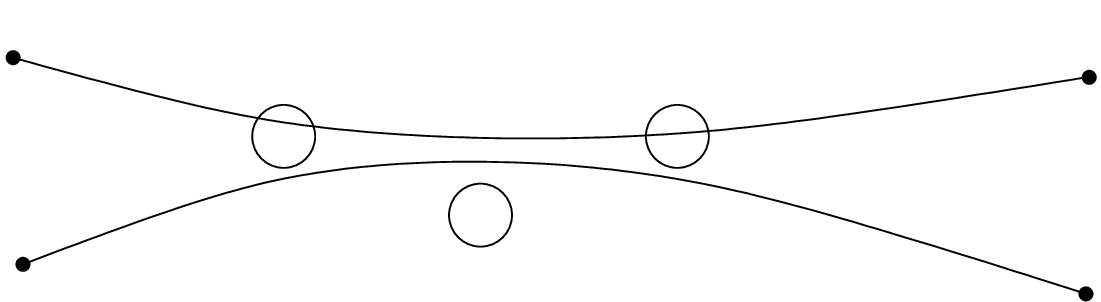}%
\end{picture}%
\setlength{\unitlength}{4144sp}%
\begingroup\makeatletter\ifx\SetFigFont\undefined%
\gdef\SetFigFont#1#2#3#4#5{%
  \reset@font\fontsize{#1}{#2pt}%
  \fontfamily{#3}\fontseries{#4}\fontshape{#5}%
  \selectfont}%
\fi\endgroup%
\begin{picture}(5018,1366)(1741,-1493)
\put(3826,-1141){\makebox(0,0)[lb]{\smash{{\SetFigFont{12}{14.4}{\familydefault}{\mddefault}{\updefault}{\color[rgb]{0,0,0}$H$}%
}}}}
\put(2926,-781){\makebox(0,0)[lb]{\smash{{\SetFigFont{12}{14.4}{\familydefault}{\mddefault}{\updefault}{\color[rgb]{0,0,0}$A$}%
}}}}
\put(4726,-781){\makebox(0,0)[lb]{\smash{{\SetFigFont{12}{14.4}{\familydefault}{\mddefault}{\updefault}{\color[rgb]{0,0,0}$B$}%
}}}}
\put(6481,-1321){\makebox(0,0)[lb]{\smash{{\SetFigFont{12}{14.4}{\familydefault}{\mddefault}{\updefault}{\color[rgb]{0,0,0}$d$}%
}}}}
\put(1846,-1141){\makebox(0,0)[lb]{\smash{{\SetFigFont{12}{14.4}{\familydefault}{\mddefault}{\updefault}{\color[rgb]{0,0,0}$c$}%
}}}}
\put(6661,-376){\makebox(0,0)[lb]{\smash{{\SetFigFont{12}{14.4}{\familydefault}{\mddefault}{\updefault}{\color[rgb]{0,0,0}$b$}%
}}}}
\put(1756,-286){\makebox(0,0)[lb]{\smash{{\SetFigFont{12}{14.4}{\familydefault}{\mddefault}{\updefault}{\color[rgb]{0,0,0}$a$}%
}}}}
\end{picture}%
\caption{Failure of (A5).}
\label{f:typeM}
\end{center}
\end{figure}
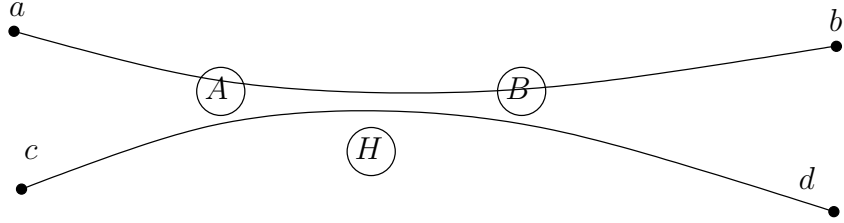

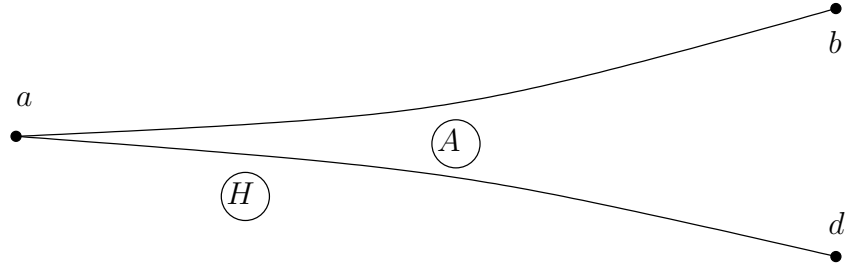
\begin{figure}[htbp]
\begin{center}
\begin{picture}(0,0)%
\includegraphics{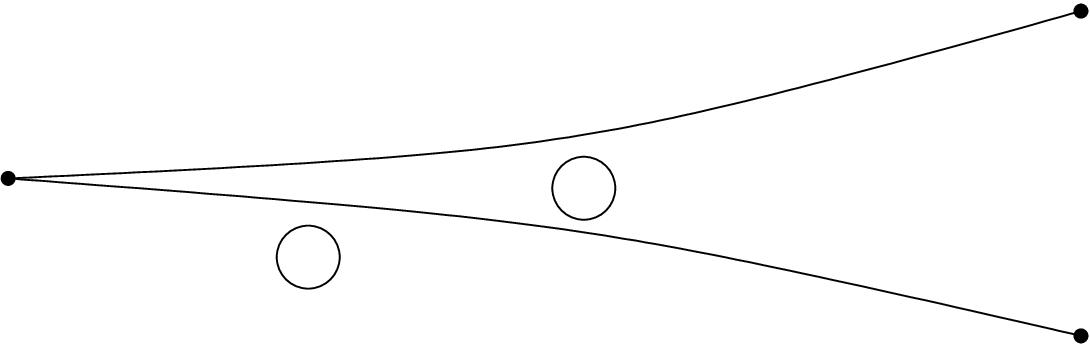}%
\end{picture}%
\setlength{\unitlength}{4144sp}%
\begingroup\makeatletter\ifx\SetFigFont\undefined%
\gdef\SetFigFont#1#2#3#4#5{%
  \reset@font\fontsize{#1}{#2pt}%
  \fontfamily{#3}\fontseries{#4}\fontshape{#5}%
  \selectfont}%
\fi\endgroup%
\begin{picture}(4981,1559)(1718,-1493)
\put(1756,-556){\makebox(0,0)[lb]{\smash{{\SetFigFont{12}{14.4}{\familydefault}{\mddefault}{\updefault}{\color[rgb]{0,0,0}$a$}%
}}}}
\put(3016,-1141){\makebox(0,0)[lb]{\smash{{\SetFigFont{12}{14.4}{\familydefault}{\mddefault}{\updefault}{\color[rgb]{0,0,0}$H$}%
}}}}
\put(4276,-826){\makebox(0,0)[lb]{\smash{{\SetFigFont{12}{14.4}{\familydefault}{\mddefault}{\updefault}{\color[rgb]{0,0,0}$A$}%
}}}}
\put(6616,-241){\makebox(0,0)[lb]{\smash{{\SetFigFont{12}{14.4}{\familydefault}{\mddefault}{\updefault}{\color[rgb]{0,0,0}$b$}%
}}}}
\put(6616,-1321){\makebox(0,0)[lb]{\smash{{\SetFigFont{12}{14.4}{\familydefault}{\mddefault}{\updefault}{\color[rgb]{0,0,0}$d$}%
}}}}
\end{picture}%
\caption{Failure of (A6).}
\label{f:typeL}
\end{center}
\end{figure}

Below we define a filtration of $\CHH$, which will be the minimal (in
an obvious sense of the word minimal) possible function satisfying
(A1)--(A7).  Before giving the definition, consider the possible ways
that (A5)--(A7) might fail to hold for some $\GHH_{a,b}$.
Either (A5) fails (in the ``middle'' of $\GHH_{a,b}$; see Figure \ref{f:typeM}), (A6) fails (on the ``left''; see Figure \ref{f:typeL}) or 
(A7) fails.

Definition
\ref{d:filtration} below can be thought of in the following way:
Start with $\NHH$, and ``fix'' every failure of axiom (A5)--(A7); after doing these repairs you will have obtained $\mc{D}^1$,
which still doesn't satisfy (A5)--(A7).  Repairing $\mc{D}^1$ yields
$\mc{D}^2$, and so on.
\begin{definition}\label{d:filtration}
For any pair $a,b \in \Upsilon$, we make the following definitions:
Let $\mc{D}^0_{a,b} = \NHH_{a,b}$, and define inductively, 
for $i \ge 0$:
\begin{eqnarray*}
\mc{M}^i_{a,b} &=& \bigcup_{c,d\in \Upsilon} 
\left(  \bigcup_{A<B \in \mc{D}^i_{a,b} \cap \mc{D}^i_{c,d}} \mc{D}^i_{c,d}[A,B] \right),	\\
\mc{L}^i_{a,b} &=& 
\bigcup_{d \in \Upsilon} \left( \bigcup_{A \in \mc{D}^i_{a,b} \cap \mc{D}^i_{a,d}}
\mc{D}^i_{a,d}(a,A] \right),\\
\mc{R}^i_{a,b} &=& 
\bigcup_{c\in \Upsilon} \left( \bigcup_{A\in \mc{D}^i_{c,b}\cap \mc{D}^i_{a,b}}
\mc{D}^i_{c,b}[A,b) \right) ,\mbox{ and} \\
\mc{D}^{i+1}_{a,b}  &=& \mc{D}^i_{a,b} \cup
\mc{M}^i_{a,b} \cup 
\mc{L}^i_{a,b} \cup
\mc{R}^i_{a,b}.
\end{eqnarray*}
Finally, we define
\[\CHH_{a,b} = \bigcup_{i=1}^\infty \mc{D}^i_{a,b}.\]
\end{definition}

Notice that the set $\mc{M}^i_{a,b}$ fixes the failure of (A5) for
$\mc{D}^i_{a,b}$, and so on.

The following is the main result of this section, which will allow us
to define preferred paths in the next section.

\begin{theorem}\label{t:ppw}
Suppose that $\mc{G}$ is a $50\delta$-separated collection
of convex subsets of a $\delta$-hyperbolic space $\Upsilon$.
Then the function $\CHH$ as defined in Definition
\ref{d:filtration} satisfies the axioms (A1)--(A7).
\end{theorem}
Definition \ref{d:filtration} is tailored so as to make (A5)--(A7) as apparent as
possible (once it is known that $\CHH$ satisfies (A2), (A5)--(A7) are 
immediate).
As explained below, the hard part of proving Theorem \ref{t:ppw} is
Axiom (A2).  We will proceed by induction.  On their own, Axioms
(A1)--(A4) for $\mc{D}^i$ do not seem strong enough to imply Axioms
(A1)--(A4) for $\mc{D}^{i+1}$.  
In order for the inductive proof to work, we need to
impose further conditions, which are encapsulated in the following
definitions.

\begin{definition} \label{d:guards}
Let $a,b \in \Upsilon$ and $C \in  \mc{D}^n_{a,b}$ with $a,b\not\in C$.
A {\em pair of $(n,a,b)$-guards} of $C$ is a pair $(Z,W)$ (each of
which may be either a point or a globule) so that
there exists $x,y \in \Upsilon$ for which:
\begin{enumerate}
\item $C \in\NHH_{x,y}$;
\item Either $Z=a=x$ or $Z \in \PHH_{x,y} \cap \mc{D}^n_{a,b}$ and $Z <C$; and
\item Either $W = b=y$ or $W \in \PHH_{x,y} \cap \mc{D}^n_{a,b}$ and
$W> C$.
\end{enumerate}

We inductively define {\em $(n,a,b)$-sentinels} of $C$ by saying that
$(n,a,b)$-guards of $C$ are $(n,a,b)$-sentinels of $C$ and that
if $D$ is an $(n,a,b)$-sentinel of $C$ then any $(n,a,b)$-guards
of $D$ are $(n,a,b)$-sentinels of $C$.
\end{definition}

\begin{remark}
For any $n, a,b$ and $C$ as in Definition \ref{d:guards}, $(n,a,b)$-guards
of $C$ come in pairs (one to the left of $C$ and one to the right).
A single $(n,a,b)$-guard may occur in many pairs.

We also remark that for any $n\ge 0$, any $a,b\in \Upsilon$ and any $C \in \mc{G}$, a pair
of $(n,a,b)$-guards of $C$ is also a pair of $(n+1,a,b)$-guards of $C$.
\end{remark}

We now introduce the property which will form the inductive hypothesis
in the proof of Theorem \ref{t:ppw}.

\begin{definition} \label{d:prop}
Given an integer $n \ge 0$, a pair $a,b\in \Upsilon$ and a globule $C \in \mc{G}$ we let 
$\Prop{n}{a}{b}{C}$ be the conjunction of the following three statements:
\begin{enumerate}
\item $C \in \mc{D}^n_{a,b}$;
\item $C \in \PHH_{a,b}$; and 
\item either (i) $C$ has a pair of $(n,a,b)$-guards; (ii) $a \in C$; or
(iii) $b \in C$.
\end{enumerate}
\end{definition}

\begin{lemma} \label{l:ifitheni+1}
For any $n \ge 0 $, any $a,b \in \Upsilon$ and any $C \in \mc{G}$, if $\Prop{n}{a}{b}{C}$ holds
then $\Prop{n+1}{a}{b}{C}$ also holds.
\end{lemma}

\begin{definition} \label{d:Sp}
Let $a,b \in \Upsilon$ and $C,D \in \mc{G}$.
Suppose that $\gamma \in \alpha(a,b)$ and that 
$\gamma$ intersects both $N_\beta(C)$ and
$N_\beta(D)$ nontrivially for some $\beta \le 2K$.
Suppose also that $C<D$ in the order on
$\mc{C}^0_{a,b}$ as in Remark \ref{r:orientation}.

Let $y$ be the last point on $\gamma$ in $N_\beta(C)$
and $z$ the first point on $\gamma$ in $N_\beta(D)$.
Define \[\Sp{C}{D}{\beta}\] 
to be the subsegment of $\gamma$ between $y$ and $z$.  
\end{definition}

\begin{proof}[Proof (Theorem {\ref{t:ppw}})]
Axiom (A4) is obvious.  The (unordered) set theoretic parts of (A1)
and (A3) are also
obvious.  The set theoretic part of (A2) is the key:  By Remark
\ref{r:orientation}, this gives sense to (and implies) all the
statements about ordered sets.  Axioms (A5)--(A7) will then follow from the
construction.  Axiom (A2) follows from the following inductive
statement:

\begin{claim} \label{c:induction}
For $i \ge 0$, all $a,b\in \Upsilon$ and all $C \in \mc{G}$ the following are
equivalent:
\begin{enumerate}
\item\label{c:Di} $C \in \mc{D}^i_{a,b}$;
\item\label{c:Pi} $\Prop{i}{a}{b}{C}$ holds, and $\Prop{i}{a}{b}{D}$ holds 
whenever $D$ is an $(i,a,b)$-sentinel of $C$.
\end{enumerate}
\end{claim}

Note that for any $i \ge 0$, any $a,b \in \Upsilon$ and any $C \in \mc{G}$,
the statement ``(\ref{c:Pi}) implies (\ref{c:Di})'' holds
tautologically (by Statement \ref{d:prop}.(1)).

We will prove Claim \ref{c:induction} by induction on $i$.

{\noindent \bf Base Case:} \label{l:basecase}
For any $a,b \in \Upsilon$ and $C \in \mc{G}$, if $C \in \mc{D}^0_{a,b}$ then
$\Prop{0}{a}{b}{C}$.  Moreover, whenever $D$ is an
$(0,a,b)$-sentinel of $C$, $\Prop{0}{a}{b}{D}$ holds.

\begin{proof}[Proof (Base case).]
Suppose $C \in \mc{D}^0_{a,b}$.  If $a \in C$ or $b\in C$ then
the result is immediate. Thus we may suppose that 
$a,b \not\in C$.

Statement \ref{d:prop}.(1) holds by assumption.

Statement \ref{d:prop}.(2) holds because $\mc{D}^0_{a,b} =
\NHH_{a,b} \subseteq \PHH_{a,b}$.

Statement \ref{d:prop}.(3) holds because $(a,b)$ form a pair of
$(0,a,b)$-guards of $C$.

Now suppose that $D$ is a $(0,a,b)$-sentinel of $C$. This
implies that $D \in \mc{D}^0_{a,b}$ and the same argument applies.
\end{proof}

{\bf Inductive Hypothesis:}
Fix $n \ge 1$. For any $0 \le i < n$, any $a,b\in \Upsilon$ and any
$C \in \mc{G}$ if $C \in \mc{D}^i_{a,b}$ then $\Prop{i}{a}{b}{C}$ is
true and $\Prop{i}{a}{b}{D}$ holds for any $(i,a,b)$-sentinel 
$D$ of $C$.

{\bf Inductive step:}
Consider $a,b\in \Upsilon$ and $C \in \mc{G}$ and suppose that $C \in \mc{D}^n_{a,b}$.  We
wish to prove that $\Prop{n}{a}{b}{C}$ holds and that 
$\Prop{n}{a}{b}{D}$ holds for any $(n,a,b)$-sentinel $D$ of $C$.  

Observe that if $D$ is an $(n,a,b)$-sentinel of $C$ then it is
in particular in $\mc{D}^n_{a,b}$.  Therefore since $C \in 
\mc{D}^n_{a,b}$ was arbitrary, it suffices to prove that 
$\Prop{n}{a}{b}{C}$ holds.

By Definition \ref{d:filtration}, one of four situations must occur:
 \begin{enumerate}
 \item \label{Case1} $C \in \mc{D}^{n-1}_{a,b}$;
 \item \label{Case2} $C \in \mc{M}^{n-1}_{a,b}$;
\item \label{Case3} $C  \in \mc{L}^{n-1}_{a,b}$; or
\item \label{Case4} $C \in \mc{R}^{n-1}_{a,b}$.
\end{enumerate}
We deal with each of these situations in turn.

{\bf Case \ref{Case1}:} In this case $\Prop{n-1}{a}{b}{C}$ holds
by induction, and by Lemma \ref{l:ifitheni+1} $\Prop{n}{a}{b}{C}$
also holds.

{\bf Case \ref{Case2}:} 
Suppose that $C \not\in \mc{D}^{n-1}_{a,b}$ but that
$C \in \mc{M}^{n-1}_{a,b}$.

Therefore, there are $c,d \in \Upsilon$ and $A,B \in \mc{G}$ so 
that
\begin{enumerate}
\item $A,B \in \mc{D}^{n-1}_{a,b}\cap \mc{D}^{n-1}_{c,d}$, and
\item $C \in \mc{D}^{n-1}_{c,d}[A,B]$.
\end{enumerate}

We wish to show that $\Prop{n}{a}{b}{C}$ holds. Statement 
\ref{d:prop}.(1) is clear.

We now prove Statement \ref{d:prop}.(2) of $\Prop{n}{a}{b}{C}$.

Since $C \in \mc{D}^{n-1}_{c,d}$, the statement $\Prop{n-1}{c}{d}{C}$
holds by the inductive hypothesis, as does $\Prop{n-1}{c}{d}{D}$ for
any $(n-1,c,d)$-sentinel $D$ of $C$.  Note that $A < C < B$, so
$c,d \not\in C$ and $C$ has a pair of $(n-1,c,d)$ guards $Z$ and $W$,
by property $\Prop{n-1}{c}{d}{C}.(3)$.  Suppose that
$x, y \in \Upsilon$ are the points associated to 
$C,Z$ and $W$ from Definition \ref{d:guards}.
\begin{figure}[htbp]
\begin{center}
\begin{picture}(0,0)%
\includegraphics{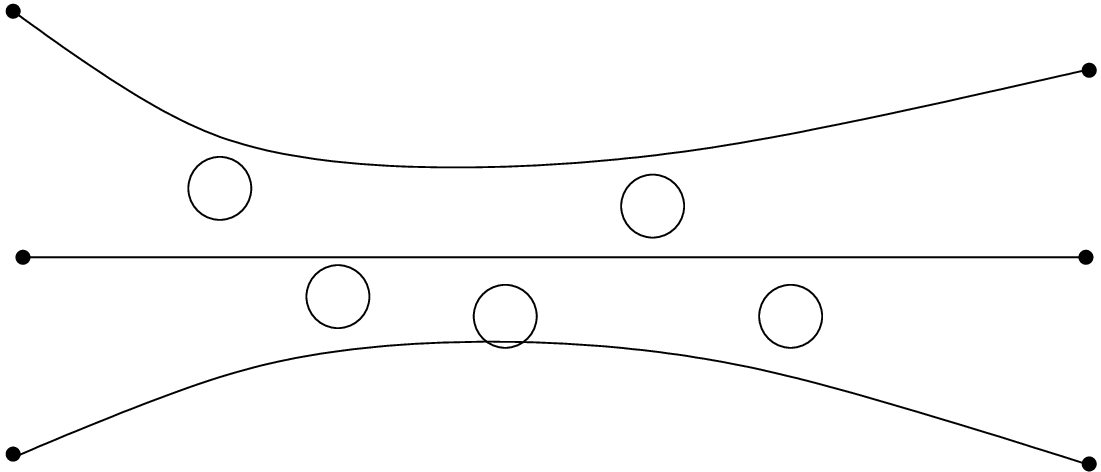}%
\end{picture}%
\setlength{\unitlength}{4144sp}%
\begingroup\makeatletter\ifx\SetFigFont\undefined%
\gdef\SetFigFont#1#2#3#4#5{%
  \reset@font\fontsize{#1}{#2pt}%
  \fontfamily{#3}\fontseries{#4}\fontshape{#5}%
  \selectfont}%
\fi\endgroup%
\begin{picture}(5018,2144)(1741,-2303)
\put(3196,-1546){\makebox(0,0)[lb]{\smash{{\SetFigFont{12}{14.4}{\familydefault}{\mddefault}{\updefault}{\color[rgb]{0,0,0}$Z$}%
}}}}
\put(5266,-1636){\makebox(0,0)[lb]{\smash{{\SetFigFont{12}{14.4}{\familydefault}{\mddefault}{\updefault}{\color[rgb]{0,0,0}$W$}%
}}}}
\put(2656,-1051){\makebox(0,0)[lb]{\smash{{\SetFigFont{12}{14.4}{\familydefault}{\mddefault}{\updefault}{\color[rgb]{0,0,0}$A$}%
}}}}
\put(4636,-1141){\makebox(0,0)[lb]{\smash{{\SetFigFont{12}{14.4}{\familydefault}{\mddefault}{\updefault}{\color[rgb]{0,0,0}$B$}%
}}}}
\put(3961,-1636){\makebox(0,0)[lb]{\smash{{\SetFigFont{12}{14.4}{\familydefault}{\mddefault}{\updefault}{\color[rgb]{0,0,0}$C$}%
}}}}
\put(1756,-2086){\makebox(0,0)[lb]{\smash{{\SetFigFont{12}{14.4}{\familydefault}{\mddefault}{\updefault}{\color[rgb]{0,0,0}$x$}%
}}}}
\put(1801,-1141){\makebox(0,0)[lb]{\smash{{\SetFigFont{12}{14.4}{\familydefault}{\mddefault}{\updefault}{\color[rgb]{0,0,0}$c$}%
}}}}
\put(1756,-466){\makebox(0,0)[lb]{\smash{{\SetFigFont{12}{14.4}{\familydefault}{\mddefault}{\updefault}{\color[rgb]{0,0,0}$a$}%
}}}}
\put(6661,-691){\makebox(0,0)[lb]{\smash{{\SetFigFont{12}{14.4}{\familydefault}{\mddefault}{\updefault}{\color[rgb]{0,0,0}$b$}%
}}}}
\put(6661,-1231){\makebox(0,0)[lb]{\smash{{\SetFigFont{12}{14.4}{\familydefault}{\mddefault}{\updefault}{\color[rgb]{0,0,0}$d$}%
}}}}
\put(6661,-2131){\makebox(0,0)[lb]{\smash{{\SetFigFont{12}{14.4}{\familydefault}{\mddefault}{\updefault}{\color[rgb]{0,0,0}$y$}%
}}}}
\end{picture}%
\caption{Case (2).}
\label{f:XYZW}
\end{center}
\end{figure}

We now define globules $P$ and $Q$.  In case $Z=x=c$, let $P=A$.
Otherwise, $Z$ and $A$ are both globules in $\mc{D}^{n-1}_{c,d}$
and we let $P = \max\{ A,Z \}$ with respect to the order on
$\mc{D}^{n-1}_{c,d}$.  Similarly, if $W=y=d$ then $Q = B$ and otherwise
$Q = \min\{ B,W \}$.  See Figure \ref{f:XYZW} for an illustrative picture.
In the figure, $P=Z$ and $B=Q$, and they are all globules.

{\bf Claim 1:} $P,Q \in \mc{C}^{K+5\delta}_{x,y}$.

We only consider $P$, as the argument for $Q$ is identical.
If $P = Z$ then $Z \in \mc{C}^K_{x,y} \subseteq \mc{C}^{K+5\delta}_{x,y}$.

Thus suppose that $P=A$ and that $A \neq Z$.  There are now two
cases, depending on whether $Z= x =c $ or not.

If $Z = x = c$, let  $\beta_{x,y}^{5\delta,K+5\delta}(Z,C)$ denote that
portion of $\beta_{x,y}$ between the last point in $N_{5\delta}(Z)$ and
the first point in $N_{K+5\delta}(C)$, and define 
$\beta_{x,y}^{5\delta,K+5\delta}(Z,C)$ similarly.
Note that if $N_{K+5\delta}(A)$ and $N_{5\delta}(x)$ are not disjoint,
then $A=P$ is certainly contained in $\mc{C}^{K+5\delta}_{x,y}$.
Otherwise, since the singleton
$\{ Z \} = \{ x\}$ is convex, Lemma \ref{l:tube} implies that the Hausdorff distance
between $\beta_{x,y}^{5\delta,K+5\delta}(Z,C)$ and
$\beta_{x,y}^{5\delta,K+5\delta}(Z,C)$ is at most $5\delta$.
But $A \in \mc{D}^{n-1}_{c,d}$, so by the inductive hypothesis
$A\in \mc{C}^K_{c,d}$. Therefore $\beta_{x,y}$ passes within
$K + 5\delta$ of $A =P$, as required.  See Figure \ref{f:claim1point}.

\begin{figure}[htbp]
\begin{center}
\begin{picture}(0,0)%
\includegraphics{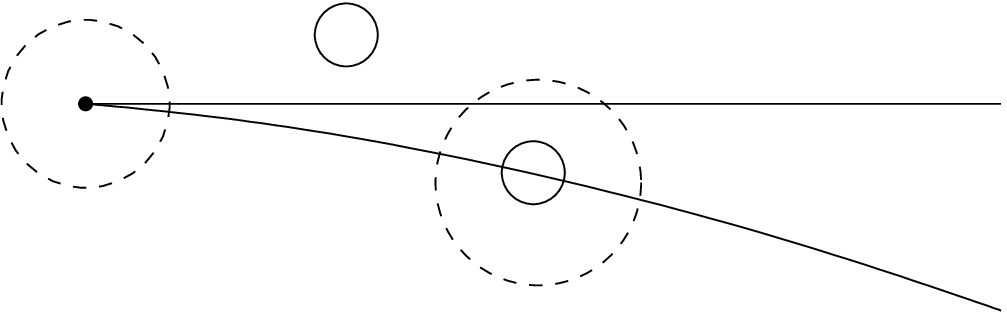}%
\end{picture}%
\setlength{\unitlength}{4144sp}%
\begingroup\makeatletter\ifx\SetFigFont\undefined%
\gdef\SetFigFont#1#2#3#4#5{%
  \reset@font\fontsize{#1}{#2pt}%
  \fontfamily{#3}\fontseries{#4}\fontshape{#5}%
  \selectfont}%
\fi\endgroup%
\begin{picture}(4589,1423)(1364,-3448)
\put(2881,-2266){\makebox(0,0)[lb]{\smash{{\SetFigFont{12}{14.4}{\familydefault}{\mddefault}{\updefault}{\color[rgb]{0,0,0}$A$}%
}}}}
\put(1486,-2401){\makebox(0,0)[lb]{\smash{{\SetFigFont{12}{14.4}{\familydefault}{\mddefault}{\updefault}{\color[rgb]{0,0,0}$Z=c=x$}%
}}}}
\put(3691,-2851){\makebox(0,0)[lb]{\smash{{\SetFigFont{12}{14.4}{\familydefault}{\mddefault}{\updefault}{\color[rgb]{0,0,0}$C$}%
}}}}
\put(4951,-3031){\makebox(0,0)[lb]{\smash{{\SetFigFont{12}{14.4}{\familydefault}{\mddefault}{\updefault}{\color[rgb]{0,0,0}$\beta_{x,y}$}%
}}}}
\put(4906,-2401){\makebox(0,0)[lb]{\smash{{\SetFigFont{12}{14.4}{\familydefault}{\mddefault}{\updefault}{\color[rgb]{0,0,0}$\beta_{c,d}$}%
}}}}
\end{picture}%
\caption{Claim 1, in case $Z=x=c$. Dotted lines indicate
$N_K(Z)$ and $N_K(C)$.}
\label{f:claim1point}
\end{center}
\end{figure}
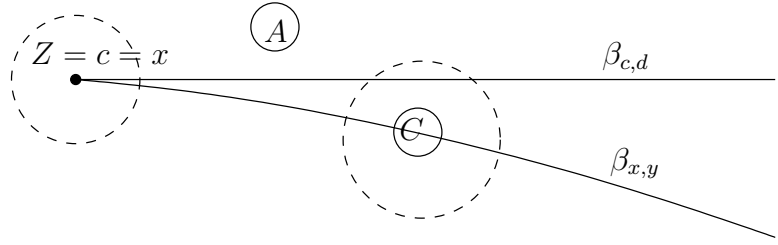

Suppose then that $Z \neq x$ or $Z\neq c$.  Then $Z$ is a globule
and is separated from $A$ by at least $L_1$.  By Lemma \ref{l:tube}
again, the Hausdorff distance between
$\beta_{c,d}^{K+5\delta}(Z,C)$ and $\beta_{x,y}^{K+5\delta}(Z,C)$
is at most $5\delta$.  Once again we know by induction that
$A \in \mc{C}^K_{c,d}$, which proves that $\beta_{x,y}$
passes within $K+5\delta$ of $A$.  See Figure \ref{f:claim1}.  This proves Claim 1.

\begin{figure}[htbp]
\begin{center}
\begin{picture}(0,0)%
\includegraphics{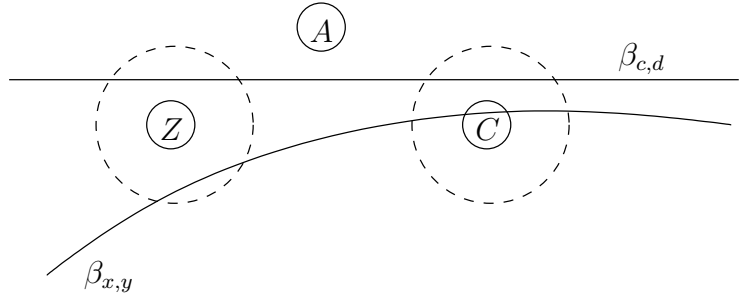}%
\end{picture}%
\setlength{\unitlength}{4144sp}%
\begingroup\makeatletter\ifx\SetFigFont\undefined%
\gdef\SetFigFont#1#2#3#4#5{%
  \reset@font\fontsize{#1}{#2pt}%
  \fontfamily{#3}\fontseries{#4}\fontshape{#5}%
  \selectfont}%
\fi\endgroup%
\begin{picture}(4389,1754)(1744,-3779)
\put(2656,-2851){\makebox(0,0)[lb]{\smash{{\SetFigFont{12}{14.4}{\familydefault}{\mddefault}{\updefault}{\color[rgb]{0,0,0}$Z$}%
}}}}
\put(3556,-2266){\makebox(0,0)[lb]{\smash{{\SetFigFont{12}{14.4}{\familydefault}{\mddefault}{\updefault}{\color[rgb]{0,0,0}$A$}%
}}}}
\put(4546,-2851){\makebox(0,0)[lb]{\smash{{\SetFigFont{12}{14.4}{\familydefault}{\mddefault}{\updefault}{\color[rgb]{0,0,0}$C$}%
}}}}
\put(2206,-3706){\makebox(0,0)[lb]{\smash{{\SetFigFont{12}{14.4}{\familydefault}{\mddefault}{\updefault}{\color[rgb]{0,0,0}$\beta_{x,y}$}%
}}}}
\put(5401,-2401){\makebox(0,0)[lb]{\smash{{\SetFigFont{12}{14.4}{\familydefault}{\mddefault}{\updefault}{\color[rgb]{0,0,0}$\beta_{c,d}$}%
}}}}
\end{picture}%
\caption{Claim 1, in case $Z$ is not $x$ and $c$. Dotted lines indicate
$N_K(Z)$ and $N_K(C)$.}
\label{f:claim1}
\end{center}
\end{figure}

{\bf Claim 2:}  $P,Q \in \mc{C}^{K+5\delta}_{a,b}$.

Once again, we only consider $P$, as the argument for $Q$ is identical.

If $P=A$ then $P \in \mc{D}^{n-1}_{a,b}$ so by induction we
have $\Prop{n-1}{a}{b}{P}$ and $P\in \mc{C}^K_{a,b} \subseteq
\mc{C}^{K+5\delta}_{a,b}$.

Suppose then that $P=Z$ and that
$Z \neq A$ (this is depicted in Figure \ref{f:XYZW}).
Certainly we also have $Z \neq B$.

We know that $A,B \in \mc{D}^{n-1}_{a,b} \cap \mc{D}^{n-1}_{c,d}$
so by the inductive hypothesis
$A,B \in \mc{C}^K_{a,b} \cap \mc{C}^K_{c,d}$.

Let $\beta_{a,b} \in \alpha(a,b)$ and $\beta_{c,d} \in
\alpha(c,d)$ be arbitrary.  Then $\beta_{a,b}$ and $\beta_{c,d}$
pass within $K$ of both $A$ and $B$.  Since $P = Z$ is
an $(n-1,c,d)$-guard of $C$, and $c \neq Z$ we know
$P\in \mc{D}^{n-1}_{c,d}$ and so by another application of
the inductive hypothesis $P\in\mc{C}^K_{c,d}$.
Therefore $\beta_{c,d}$ passes within $K$ of $P$.

By Lemma \ref{l:tube}, the Hausdorff distance between
$\beta_{c,d}^{K+5\delta}(A,B)$ and $\beta_{a,b}^{K+5\delta}(A,B)$
is at most $5\delta$.  Thus, $\beta_{a,b}$ passes within $K+5\delta$
of $P$, proving Claim 2 (see Figure \ref{f:claim2}).

\begin{figure}[htbp]
\begin{center}
\begin{picture}(0,0)%
\includegraphics{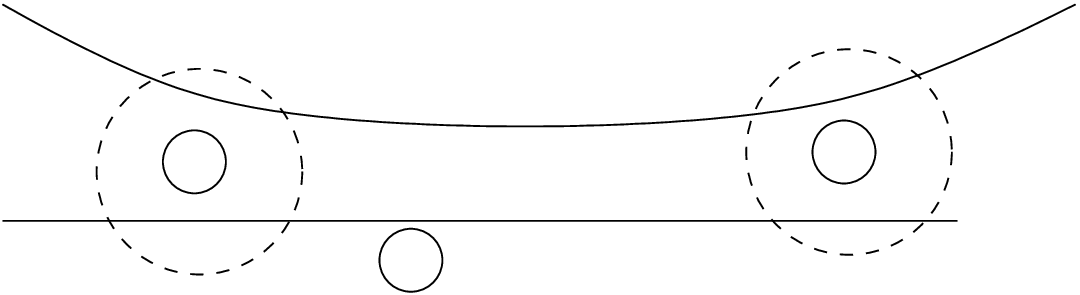}%
\end{picture}%
\setlength{\unitlength}{4144sp}%
\begingroup\makeatletter\ifx\SetFigFont\undefined%
\gdef\SetFigFont#1#2#3#4#5{%
  \reset@font\fontsize{#1}{#2pt}%
  \fontfamily{#3}\fontseries{#4}\fontshape{#5}%
  \selectfont}%
\fi\endgroup%
\begin{picture}(4929,1341)(1744,-2830)
\put(3736,-1951){\makebox(0,0)[lb]{\smash{{\SetFigFont{12}{14.4}{\familydefault}{\mddefault}{\updefault}{\color[rgb]{0,0,0}$\beta_{a,b}$}%
}}}}
\put(6211,-2536){\makebox(0,0)[lb]{\smash{{\SetFigFont{12}{14.4}{\familydefault}{\mddefault}{\updefault}{\color[rgb]{0,0,0}$\beta_{c,d}$}%
}}}}
\put(2566,-2311){\makebox(0,0)[lb]{\smash{{\SetFigFont{12}{14.4}{\familydefault}{\mddefault}{\updefault}{\color[rgb]{0,0,0}$A$}%
}}}}
\put(5536,-2266){\makebox(0,0)[lb]{\smash{{\SetFigFont{12}{14.4}{\familydefault}{\mddefault}{\updefault}{\color[rgb]{0,0,0}$B$}%
}}}}
\put(3556,-2761){\makebox(0,0)[lb]{\smash{{\SetFigFont{12}{14.4}{\familydefault}{\mddefault}{\updefault}{\color[rgb]{0,0,0}$Z$}%
}}}}
\end{picture}%
\caption{Claim 2, in case $A<Z$. Dotted lines indicate
$N_K(A)$ and $N_K(B)$.}
\label{f:claim2}
\end{center}
\end{figure}
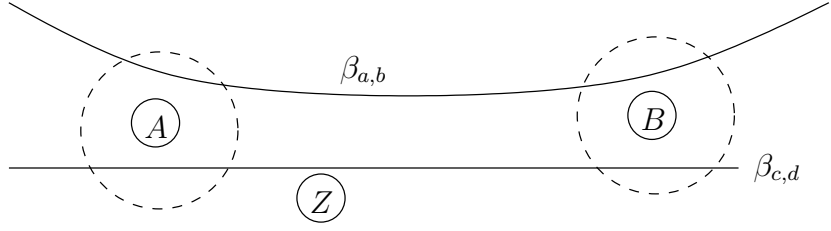

Note that whatever $P$ and $Q$ are, we always have $P < C < Q$.
Since $P,Q\in \mc{C}^{K+5\delta}_{a,b} \cap\mc{C}^{K+5\delta}_{x,y}$,
Lemma \ref{l:tube} implies that the Hausdorff distance between
$\beta_{a,b}^{K+10\delta}(P,Q)$ and $\beta_{x,y}^{K+10\delta}(P,Q)$
is at most $5\delta$.  Since $\beta_{x,y}^{K+10\delta}(P,Q)$ intersects
$C$, the geodesic $\beta_{a,b}$ passes within $5\delta$ of $C$ (Figure
\ref{f:PQ}). 
This proves Statement \ref{d:prop}.(2) of $\Prop{n}{a}{b}{C}$.
\begin{figure}[htbp]
\begin{center}
\begin{picture}(0,0)%
\includegraphics{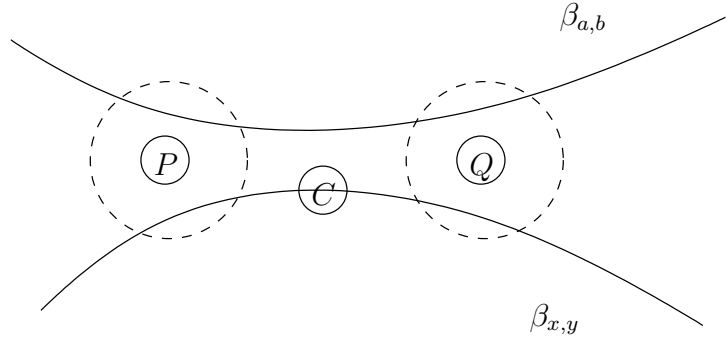}%
\end{picture}%
\setlength{\unitlength}{4144sp}%
\begingroup\makeatletter\ifx\SetFigFont\undefined%
\gdef\SetFigFont#1#2#3#4#5{%
  \reset@font\fontsize{#1}{#2pt}%
  \fontfamily{#3}\fontseries{#4}\fontshape{#5}%
  \selectfont}%
\fi\endgroup%
\begin{picture}(4299,2044)(1789,-3824)
\put(5086,-1951){\makebox(0,0)[lb]{\smash{{\SetFigFont{12}{14.4}{\familydefault}{\mddefault}{\updefault}{\color[rgb]{0,0,0}$\beta_{a,b}$}%
}}}}
\put(2656,-2851){\makebox(0,0)[lb]{\smash{{\SetFigFont{12}{14.4}{\familydefault}{\mddefault}{\updefault}{\color[rgb]{0,0,0}$P$}%
}}}}
\put(3601,-3031){\makebox(0,0)[lb]{\smash{{\SetFigFont{12}{14.4}{\familydefault}{\mddefault}{\updefault}{\color[rgb]{0,0,0}$C$}%
}}}}
\put(4906,-3751){\makebox(0,0)[lb]{\smash{{\SetFigFont{12}{14.4}{\familydefault}{\mddefault}{\updefault}{\color[rgb]{0,0,0}$\beta_{x,y}$}%
}}}}
\put(4546,-2851){\makebox(0,0)[lb]{\smash{{\SetFigFont{12}{14.4}{\familydefault}{\mddefault}{\updefault}{\color[rgb]{0,0,0}$Q$}%
}}}}
\end{picture}%
\caption{Case (2), Statement \ref{d:prop}.(2). Dotted lines indicate
$N_{K+5\delta}(P)$ and $N_{K+5\delta}(Q)$.}
\label{f:PQ}
\end{center}
\end{figure}

We now prove Statement \ref{d:prop}.(3) of $\Prop{n}{a}{b}{C}$.
We claim that $P$ and $Q$ form a pair of $(n,a,b)$-guards
of $C$.  It is certainly true that $P < C< Q$, and that keeping
the same $x$ and $y$ we still have $C \in \NHH_{x,y}$.
Note that $P$ and $Q$ are globules.
Thus we have to prove 
\begin{enumerate}
\item \label{subCase(i)} $P,Q \in \mc{D}^n_{a,b}$; and
\item \label{subCase(ii)} $P,Q \in \PHH_{x,y}$.
\end{enumerate}

Once again, we prove these statements for $P$, the 
arguments for $Q$ being analogous.

Suppose that $P = A$.  Then $P \in \mc{D}^{n-1}(a,b)\subseteq
\mc{D}^n_{a,b}$.
Otherwise, $P = Z$ and $A < P < B$.  Therefore
$P \in \mc{D}^{n-1}_{c,d}[A,B]$, which implies,
by the definition of $\mc{D}^n$, that $P \in \mc{M}^{n-1}_{a,b}
\subseteq \mc{D}^n_{a,b}$.
This proves Statement \ref{subCase(i)} above.

We now prove Statement \ref{subCase(ii)}.  If $P = Z$ then
$P \in \PHH_{x,y}$, so we are done.  Thus suppose that
$P = A$, and $Z < A$.  

We know that $A \in \mc{D}^{n-1}(c,d)$.  Certainly, $A<C$
so $d \not\in A$.  Suppose that $c \in A$.  This forces
$Z=x=c$, by Definition \ref{d:guards}, so $x \in A$ and
$A \in \mc{C}^0_{x,y} \subseteq \mc{C}^K_{x,y}$.

Suppose now that $c\not\in A$.  Then the inductive
hypothesis and property $\Prop{n-1}{c}{d}{A}$
imply that $A$ has a pair of
$(n-1,c,d)$-guards, which we denote $Z_A$ and $W_A$, and
an associated pair of points $x_A, y_A \in \Upsilon$, as in
Definition \ref{d:guards}.

We know that $A \in \mc{C}^0_{x_A,y_A}$ and that 
\begin{enumerate}
\item either $c=x_A=Z_A$ or $Z_A \in \mc{C}^K_{x_A,y_A}
\cap \mc{D}^{n-1}_{c,d}$; and
\item either $d=y_A=W_A$ or $W_A \in \mc{C}^K_{x_A,y_A}
\cap   \mc{D}^{n-1}_{c,d}$.
\end{enumerate}

In case $Z_A$ is a globule, we have $Z_A< A$, and in case
$W_A$ is a globule we have $A < W_A$.
(See Figure \ref{f:c3c3s2}).
\begin{figure}[htbp]
\begin{center}
\begin{picture}(0,0)%
\includegraphics{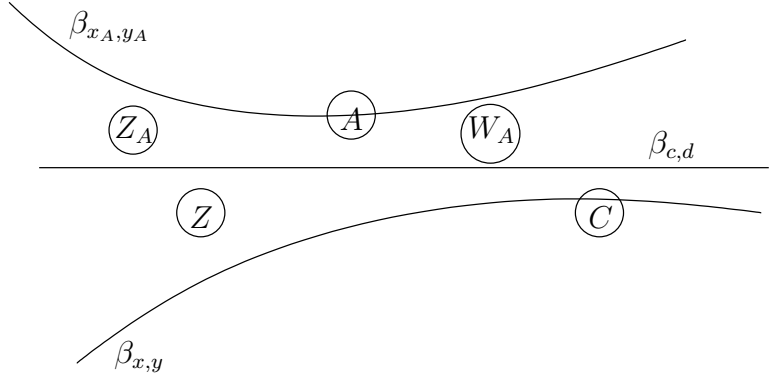}%
\end{picture}%
\setlength{\unitlength}{4144sp}%
\begingroup\makeatletter\ifx\SetFigFont\undefined%
\gdef\SetFigFont#1#2#3#4#5{%
  \reset@font\fontsize{#1}{#2pt}%
  \fontfamily{#3}\fontseries{#4}\fontshape{#5}%
  \selectfont}%
\fi\endgroup%
\begin{picture}(4569,2245)(1564,-3734)
\put(3556,-2266){\makebox(0,0)[lb]{\smash{{\SetFigFont{12}{14.4}{\familydefault}{\mddefault}{\updefault}{\color[rgb]{0,0,0}$A$}%
}}}}
\put(5401,-2401){\makebox(0,0)[lb]{\smash{{\SetFigFont{12}{14.4}{\familydefault}{\mddefault}{\updefault}{\color[rgb]{0,0,0}$\beta_{c,d}$}%
}}}}
\put(2656,-2851){\makebox(0,0)[lb]{\smash{{\SetFigFont{12}{14.4}{\familydefault}{\mddefault}{\updefault}{\color[rgb]{0,0,0}$Z$}%
}}}}
\put(5041,-2851){\makebox(0,0)[lb]{\smash{{\SetFigFont{12}{14.4}{\familydefault}{\mddefault}{\updefault}{\color[rgb]{0,0,0}$C$}%
}}}}
\put(4321,-2311){\makebox(0,0)[lb]{\smash{{\SetFigFont{12}{14.4}{\familydefault}{\mddefault}{\updefault}{\color[rgb]{0,0,0}$W_A$}%
}}}}
\put(2206,-2311){\makebox(0,0)[lb]{\smash{{\SetFigFont{12}{14.4}{\familydefault}{\mddefault}{\updefault}{\color[rgb]{0,0,0}$Z_A$}%
}}}}
\put(1936,-1681){\makebox(0,0)[lb]{\smash{{\SetFigFont{12}{14.4}{\familydefault}{\mddefault}{\updefault}{\color[rgb]{0,0,0}$\beta_{x_A,y_A}$}%
}}}}
\put(2206,-3661){\makebox(0,0)[lb]{\smash{{\SetFigFont{12}{14.4}{\familydefault}{\mddefault}{\updefault}{\color[rgb]{0,0,0}$\beta_{x,y}$}%
}}}}
\end{picture}%
\caption{Case (2), Statement \ref{d:prop}.(3) (one of many possible
arrangements).}
\label{f:c3c3s2}
\end{center}
\end{figure}

We now define $P_A$, which may be either a point or a globule.
In case $Z=x=c$, define $P_A = Z_A$ (which may be either a point
or a globule). In case $Z$ is a globule, we consider whether $Z_A$
is a point or a globule.  If $Z_A$ is a point then $P_A = Z$.  Otherwise
$P_A= \max\{ Z,Z_A\}$.  Note that if $Z$ is a globule, so is $P_A$.

We also define a globule $Q_A$.  If $W_A= y_A = d$ then
$Q_A= C$.  If $W_A$ is a globule then $Q_A= \min\{ W_A,C \}$.

{\bf Claim 3:} If $P_A$ is a globule then $P_A\in \mc{C}^{K+5\delta}_{x_A,y_A} \cap\mc{C}^{K+5\delta}_{x,y}$. In any case
$Q_A \in \mc{C}^{K+5\delta}_{x_A,y_A}\cap \mc{C}^{K+5\delta}_{x,y}$.

The proof of this is similar to those of Claims 1 and 2.

Let $\beta_{x_A,y_A}$ be any geodesic between $x_A$ and $y_A$.
In case $P_A$ is a globule, Claim 3 and Lemma \ref{l:tube}
imply that the Hausdorff
distance between $\beta_{x_A,y_A}^{K+10\delta}(P_A,Q_A)$
and $\beta_{x,y}^{K+10\delta}(P_A,Q_A)$ is at most $5\delta$.
Since $P_A <A<Q_A$ and $\beta_{x_A,y_A}$ intersects
$A$, the path $\beta_{x,y}$ passes within $5\delta$ of $A$.
Thus in this case $P= A\in   \mc{C}^K_{x,y}$,as required.

Suppose that $P_A$ is a point (see Figure \ref{f:c3c3s2corner}).  
In this case, $Z=x=c$ and $Z_A=x=x_A$, so $P_A = x=x_A$.  In case $N_{5\delta}(P_A)$ and
$A$ are not disjoint, $x$ lies within $5\delta$ of $A$,
and certainly $A\in \mc{C}^K_{x,y}$.     Otherwise, the singleton
$\{ x\}$ is convex so Lemma \ref{l:tube} implies that the
Hausdorff distance between those parts of $\beta_{x,y}$ and
$\beta_{x_A,y_A}$ between $N_{5\delta}(x)$ and $N_{K+5\delta}(C)$
are at Hausdorff distance at most $5\delta$ from each other.
Since $A \in \mc{C}^0_{x_A,y_A}$ this implies once again that
$P=A\in\mc{C}^K_{x,y}$ as required.

\begin{figure}[htbp]
\begin{center}
\begin{picture}(0,0)%
\includegraphics{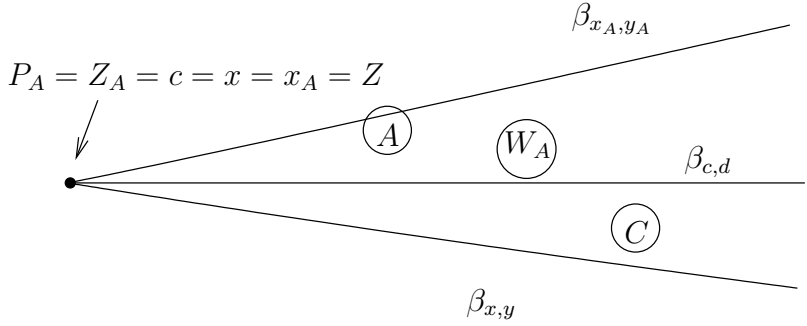}%
\end{picture}%
\setlength{\unitlength}{4144sp}%
\begingroup\makeatletter\ifx\SetFigFont\undefined%
\gdef\SetFigFont#1#2#3#4#5{%
  \reset@font\fontsize{#1}{#2pt}%
  \fontfamily{#3}\fontseries{#4}\fontshape{#5}%
  \selectfont}%
\fi\endgroup%
\begin{picture}(4797,1954)(1336,-3329)
\put(3556,-2266){\makebox(0,0)[lb]{\smash{{\SetFigFont{12}{14.4}{\familydefault}{\mddefault}{\updefault}{\color[rgb]{0,0,0}$A$}%
}}}}
\put(5401,-2401){\makebox(0,0)[lb]{\smash{{\SetFigFont{12}{14.4}{\familydefault}{\mddefault}{\updefault}{\color[rgb]{0,0,0}$\beta_{c,d}$}%
}}}}
\put(5041,-2851){\makebox(0,0)[lb]{\smash{{\SetFigFont{12}{14.4}{\familydefault}{\mddefault}{\updefault}{\color[rgb]{0,0,0}$C$}%
}}}}
\put(4321,-2311){\makebox(0,0)[lb]{\smash{{\SetFigFont{12}{14.4}{\familydefault}{\mddefault}{\updefault}{\color[rgb]{0,0,0}$W_A$}%
}}}}
\put(4096,-3256){\makebox(0,0)[lb]{\smash{{\SetFigFont{12}{14.4}{\familydefault}{\mddefault}{\updefault}{\color[rgb]{0,0,0}$\beta_{x,y}$}%
}}}}
\put(1351,-1906){\makebox(0,0)[lb]{\smash{{\SetFigFont{12}{14.4}{\familydefault}{\mddefault}{\updefault}{\color[rgb]{0,0,0}$P_A=Z_A=c=x=x_A=Z$}%
}}}}
\put(4726,-1546){\makebox(0,0)[lb]{\smash{{\SetFigFont{12}{14.4}{\familydefault}{\mddefault}{\updefault}{\color[rgb]{0,0,0}$\beta_{x_A,y_A}$}%
}}}}
\end{picture}%
\caption{Case (2), Statement \ref{d:prop}.(3) (another of many possible
arrangements).}
\label{f:c3c3s2corner}
\end{center}
\end{figure}

This finally proves that $P$ and $Q$ form a pair of $(n,a,b)$-guards
of $C$, and finishes the proof of Case (2).

{\bf Case 3:}  Suppose that $C \not\in \mc{D}^{n-1}_{a,b} \cup
\mc{M}^{n-1}_{a,b}$, but that $C \in \mc{L}^{n-1}_{a,b}$.

Therefore, there is $d\in \Upsilon$ and $A \in \mc{D}^{n-1}_{a,b}
\cap\mc{D}^{n-1}_{a,d}$ so that
$C \in \mc{D}^{n-1}_{a,d}(a,A]$.

We are required to show that $\Prop{n}{a}{b}{C}$ holds.
As usual, statement \ref{d:prop}.(1) is immediate.

The inductive hypothesis
implies that $\Prop{n-1}{a}{d}{C}$ holds.

Since $C \not\in \mc{D}^{n-1}_{a,b}$ and $A \in \mc{D}^{n-1}_{a,b}$ are globules, they are different and $C< A$. Therefore, $d \not\in C$, and
so either $a \in C$ or $C$ has a pair of $(n-1,a,d)$-guards.

If $a \in C$ then statements \ref{d:prop}.(2) and \ref{d:prop}.(3)
are clear.

Therefore we suppose that $a \not\in C$.
We now prove \ref{d:prop}.(2).

Let $Z$ and $W$ be a pair of $(n-1,a,d)$-guards for $C$, with
associated points $x,y \in\Upsilon$.
Note that either $Z=x=a$ or $Z \in \mc{C}^K_{x,y}
\cap \mc{D}^{n-1}_{a,d}$.  Also, either $W=y=d$ or
$W\in \mc{C}^K_{x,y} \cap\mc{D}^{n-1}_{a,d}$.

Let $P=Z$ and $Q  =\min\{  A,W \}$.  Note that $P$ may be
a point or a globule whereas $Q$ is definitely a globule.

Then following is entirely analogous to Claims 1, 2 and 3:

{\bf Claim 4:}  If $P$ is a globule then $P\in \mc{C}^{K+5\delta}_{a,b}
\cap \mc{C}^{K+5\delta}_{x,y}$. In any case $Q\in 
\mc{C}^{K+5\delta}_{a,b}\cap \mc{C}^{K+5\delta}_{x,y}$.

An identical argument to that which followed Claim 3 implies
that $C \in \mc{C}^K_{a,b}$, which proves \ref{d:prop}.(2).

We now prove \ref{d:prop}.(3).  Note that $b \not\in C$ and
we are assuming that $a\not\in C$. Therefore we are required
to show that $C$ has a pair of $(n,a,b$)-guards.  We will
show that $P$ and $Q$ form such a pair (with associated
points $x$ and $y$).

If $P=a$ then $P=x$ also, so \ref{d:guards}.(2) holds in this
case.  Otherwise, $P=Z$ is a globule and $P<C<A$ (See Figure
\ref{f:L1}).  Since 
$P,A \in \mc{D}^{n-1}_{a,d}$, and $A \in \mc{D}^{n-1}_{a,b}$
 this implies that $P \in \mc{L}^{n-1}_{a,b} \subseteq \mc{D}^n_{a,b}$.
 Thus we have proved \ref{d:guards}.(2).
 
 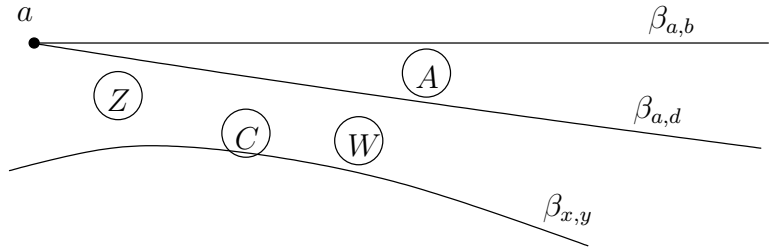
\begin{figure}[htbp]
\begin{center}
\begin{picture}(0,0)%
\includegraphics{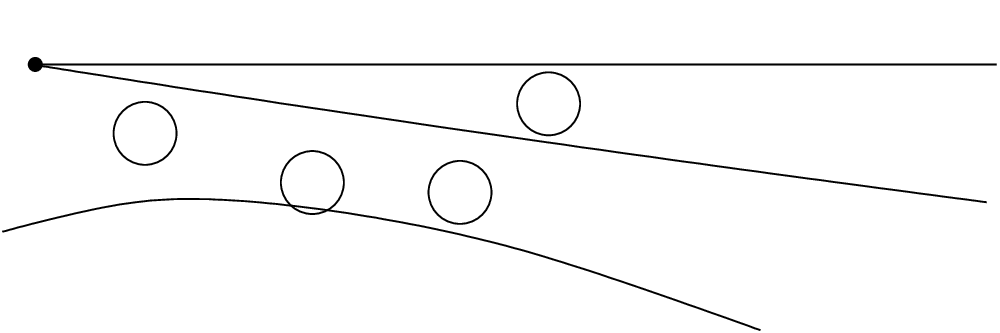}%
\end{picture}%
\setlength{\unitlength}{4144sp}%
\begingroup\makeatletter\ifx\SetFigFont\undefined%
\gdef\SetFigFont#1#2#3#4#5{%
  \reset@font\fontsize{#1}{#2pt}%
  \fontfamily{#3}\fontseries{#4}\fontshape{#5}%
  \selectfont}%
\fi\endgroup%
\begin{picture}(4569,1521)(1564,-3718)
\put(2926,-3121){\makebox(0,0)[lb]{\smash{{\SetFigFont{12}{14.4}{\familydefault}{\mddefault}{\updefault}{\color[rgb]{0,0,0}$C$}%
}}}}
\put(4006,-2761){\makebox(0,0)[lb]{\smash{{\SetFigFont{12}{14.4}{\familydefault}{\mddefault}{\updefault}{\color[rgb]{0,0,0}$A$}%
}}}}
\put(1621,-2356){\makebox(0,0)[lb]{\smash{{\SetFigFont{12}{14.4}{\familydefault}{\mddefault}{\updefault}{\color[rgb]{0,0,0}$a$}%
}}}}
\put(5401,-2401){\makebox(0,0)[lb]{\smash{{\SetFigFont{12}{14.4}{\familydefault}{\mddefault}{\updefault}{\color[rgb]{0,0,0}$\beta_{a,b}$}%
}}}}
\put(2161,-2896){\makebox(0,0)[lb]{\smash{{\SetFigFont{12}{14.4}{\familydefault}{\mddefault}{\updefault}{\color[rgb]{0,0,0}$Z$}%
}}}}
\put(3601,-3166){\makebox(0,0)[lb]{\smash{{\SetFigFont{12}{14.4}{\familydefault}{\mddefault}{\updefault}{\color[rgb]{0,0,0}$W$}%
}}}}
\put(5311,-2941){\makebox(0,0)[lb]{\smash{{\SetFigFont{12}{14.4}{\familydefault}{\mddefault}{\updefault}{\color[rgb]{0,0,0}$\beta_{a,d}$}%
}}}}
\put(4771,-3526){\makebox(0,0)[lb]{\smash{{\SetFigFont{12}{14.4}{\familydefault}{\mddefault}{\updefault}{\color[rgb]{0,0,0}$\beta_{x,y}$}%
}}}}
\end{picture}%
\caption{The proof of Case (3), Statement $\Prop{n}{a}{b}{C}$, in case
$C$ has $(n-1,a,d)$-guards.}
\label{f:L1}
\end{center}
\end{figure}

 We now show \ref{d:guards}.(3).  Note that $Q$ is a globule
 and that $C < Q$.  If $Q=W$ then $W$ is a globule and
 $W \in \mc{C}^K_{x,y} \cap \mc{D}^{n-1}_{a,d}$.  In this case
 also $W\in \mc{D}^{n-1}_{a,d}(a,A]\subseteq \mc{L}^{n-1}_{a,b}
 \subseteq \mc{D}^n_{a,b}$.  This proves \ref{d:guards}.(3)
 in this case.
 
 Suppose then that $Q=A\neq W$.  Then $Q \in \mc{D}^n_{a,b}$
 and we need only prove that $Q \in \mc{C}^K_{x,y}$.
 
 By the inductive hypothesis, since $A$ is an $(n-1,a,d)$-sentinel
 of $C$ and $C \in \mc{D}^{n-1}_{a,d}$, property
 $\Prop{n-1}{a}{d}{A}$ holds.
 Certainly $C<A$, and both are globules, so $a \not\in A$.
 If $d \in A$ then since $Q=A$ and $A \neq W$, it cannot
 be that $W$ is a globule.  Therefore in this case $d=W=y$,
 so $y \in A$ and $A\in \mc{C}^0_{x,y} \subseteq \mc{C}^K_{x,y}$.
 
 Suppose then that $d\not\in A$. In this case, \ref{d:prop}.(3)
 implies that $A$ has a pair of $(n-1,a,d)$-guards, which we
 denote by $Z_A$ and $W_A$, with associated points
 $x_A,y_A\in \Upsilon$ (See Figure \ref{f:L2}).  As before, we define
 $P_A=  \max \{ Z_A,C\}$ and $Q_A = \min\{ W_A,W \}$.  Note
 that $P_A$ is always a globule but $Q_A$ might be a point.
 
 \begin{figure}[htbp]
\begin{center}
\begin{picture}(0,0)%
\includegraphics{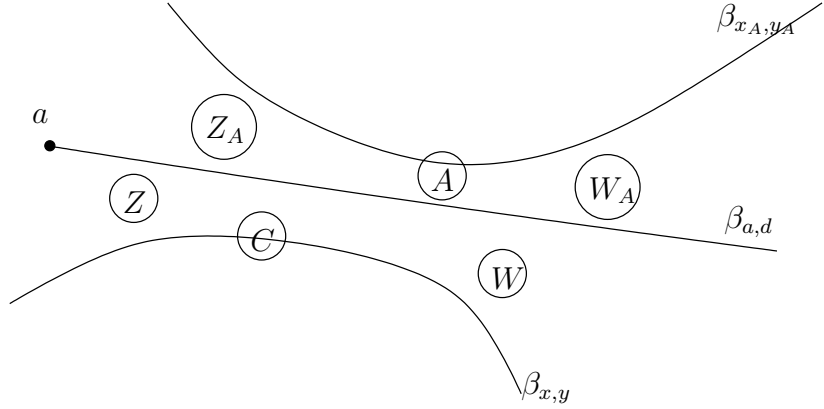}%
\end{picture}%
\setlength{\unitlength}{4144sp}%
\begingroup\makeatletter\ifx\SetFigFont\undefined%
\gdef\SetFigFont#1#2#3#4#5{%
  \reset@font\fontsize{#1}{#2pt}%
  \fontfamily{#3}\fontseries{#4}\fontshape{#5}%
  \selectfont}%
\fi\endgroup%
\begin{picture}(4884,2449)(1474,-4049)
\put(2926,-3121){\makebox(0,0)[lb]{\smash{{\SetFigFont{12}{14.4}{\familydefault}{\mddefault}{\updefault}{\color[rgb]{0,0,0}$C$}%
}}}}
\put(4006,-2761){\makebox(0,0)[lb]{\smash{{\SetFigFont{12}{14.4}{\familydefault}{\mddefault}{\updefault}{\color[rgb]{0,0,0}$A$}%
}}}}
\put(4366,-3346){\makebox(0,0)[lb]{\smash{{\SetFigFont{12}{14.4}{\familydefault}{\mddefault}{\updefault}{\color[rgb]{0,0,0}$W$}%
}}}}
\put(1621,-2356){\makebox(0,0)[lb]{\smash{{\SetFigFont{12}{14.4}{\familydefault}{\mddefault}{\updefault}{\color[rgb]{0,0,0}$a$}%
}}}}
\put(2161,-2896){\makebox(0,0)[lb]{\smash{{\SetFigFont{12}{14.4}{\familydefault}{\mddefault}{\updefault}{\color[rgb]{0,0,0}$Z$}%
}}}}
\put(4546,-3976){\makebox(0,0)[lb]{\smash{{\SetFigFont{12}{14.4}{\familydefault}{\mddefault}{\updefault}{\color[rgb]{0,0,0}$\beta_{x,y}$}%
}}}}
\put(2656,-2446){\makebox(0,0)[lb]{\smash{{\SetFigFont{12}{14.4}{\familydefault}{\mddefault}{\updefault}{\color[rgb]{0,0,0}$Z_A$}%
}}}}
\put(4951,-2806){\makebox(0,0)[lb]{\smash{{\SetFigFont{12}{14.4}{\familydefault}{\mddefault}{\updefault}{\color[rgb]{0,0,0}$W_A$}%
}}}}
\put(5761,-2986){\makebox(0,0)[lb]{\smash{{\SetFigFont{12}{14.4}{\familydefault}{\mddefault}{\updefault}{\color[rgb]{0,0,0}$\beta_{a,d}$}%
}}}}
\put(5716,-1771){\makebox(0,0)[lb]{\smash{{\SetFigFont{12}{14.4}{\familydefault}{\mddefault}{\updefault}{\color[rgb]{0,0,0}$\beta_{x_A,y_A}$}%
}}}}
\end{picture}%
\caption{Case (3), proving that $P$ and $Q$ form a pair
of $(n,a,b)$-guards for $C$ (when $Q=A$).}
\label{f:L2}
\end{center}
\end{figure}
 
 We now make the customary
 
 {\bf Claim 5:}  If $Q_A$ is a globule then $Q_A \in  \mc{C}^{K+5\delta}_{x_A,y_A} \cap \mc{C}^{K+5\delta}_{x,y}$. In any case
$P_A\in  \mc{C}^{K+5\delta}_{x_A,y_A} \cap \mc{C}^{K+5\delta}_{x,y}$.

Since $A \in \mc{C}^0_{x_A,y_A}$, the argument following Claim 3 now implies that  $A \in \mc{C}^K_{x,y}$ as required.
This establishes \ref{d:guards}.(3).  This proves that property
$\Prop{n}{a}{b}{C}$, which finally completes the proof of
Case 3.

{\bf Case 4:} Finally, suppose that $C \not\in \mc{D}^{n-1}_{a,b}$
but $C \in \mc{R}^{n-1}_{a,b}$.

The proof that $\Prop{n}{a}{b}{C}$ holds is identical (with a 
left-right reflection) to that of Case 3.

This completes the proof of Claim \ref{c:induction}, which as
we noted above implies that axiom (A2) holds for the collections
$\CHH_{A,B}$.  In turn, by Remark \ref{r:orientation}, this gives a
coherent order on each of the sets.  The function $\CHH$ was
designed to make axioms (A5)--(A7) immediately apparent, once
axiom (A2) holds, and axioms (A1), (A3) and (A4) are now clear
also.

This completes the proof of Theorem \ref{t:ppw}.
\end{proof}

\begin{corollary}
Suppose that $A,B \in \mc{G}$, and that
$a_1,a_2 \in A$ and $b_1,b_2\in B$.  Then
$\CHH_{a_1,b_1} = \CHH_{a_2,b_2}$.

Also, if $x \in \Upsilon$ then $\CHH_{a_1,x} = \CHH_{a_2,x}$.
\end{corollary}
\begin{proof}
Axiom (A1) implies that $A,B \in \CHH_{a_1,b_1}
\cap \CHH_{a_2,b_2}$.  An application of (A5)
implies that $\CHH_{a_1,b_1} \subset \CHH_{a_2,b_2}$,
and, symmetrically, $\CHH_{a_2,b_2} \subset \CHH_{a_1,b_1}$
(note that we are implicitly using the order and separation properties
of globules).

The proof of the second assertion is similar (using axiom (A7)
instead).
\end{proof}

\begin{notation} \label{n:horofamilies}
Suppose that $A,B\in \mc{G}$.  Then we denote the set
$\CHH_{a,b}$ for $a\in A$ and $b \in B$ by $\CHH_{A,B}$.
This is well-defined by the above corollary.

Similarly, if $x\in \Upsilon$ and $A \in \mc{G}$, then the sets
$\CHH_{A,x}$ and $\CHH_{x,A}$ are well-defined.
\end{notation}

\def\Hpairs {\mc{H}}

\def\pp {p} 
\def\ppbar {\bar{\pp}}
\def\Ltwo {3 L_1}  
\newcommand{\FirstH}[3] {M^{#1}_{{#2}{#3}}}
\newcommand{\LastH}[3] {N^{#1}_{{#2}{#3}}}
\def\Pairs {\mathfrak P}
\def\Hdist{3K+21\delta+14} 
\def\Displace{K+12\delta+9}
\def\Hexradius{K+2\delta}
\def\Jumpbound{4K+18\delta+9}
\def\Waitbound{4K+18\delta+9}
\def\Pepsilon{20 K+120\delta+72} 
\def\Deltaprime{6K+48\delta+28}

\def\corner {bite}
\def\bite {bite}
\def\Corner {Bite}
\def\nibble {nibble}
\def\Nibble {Nibble}
\def\dip {dip}
\def\Dip {Dip}
\def\plunge {plunge}
\def\Plunge {Plunge}

\section{Preferred paths, preferred triangles, and skeletal fillings}\label{section:pp}

Throughout this section, and the remainder of the paper,
$G$ will denote a torsion-free
group which is hyperbolic relative to a (finite) collection of
(finitely generated) subgroups $P$.  We suppose that
$S$ is a finite compatible generating set for $G$, and that $\delta$
is the constant of hyperbolicity for $X(G,P,S)$.

In this section we define {\em preferred paths}.  
Preferred paths give an equivariant, symmetric
quasi-geodesic bicombing of the space $X=X(G,P,S)$,
and will be a key to our construction of the bicombing
in Section \ref{section:bicombing} below.  Having defined preferred
paths, we then use them
to understand the combinatorial types of 
triangles whose sides are preferred paths.

We fix some notation and terminology:
Recall that $\delta$ gives rise to the constants $K = 10 \delta$,
$L_1=\Lone$, and $L_2=\Ltwo$.  
Denote by $\mc{H}$ the collection of $L_1$-horoballs in
$X$ (Definition \ref{d:horoball}).  In this section, the word
`horoball' (without a prefix) will refer to an element of $\mc{H}$.

Suppose that $P$ is a horoball and $N>0$.  There is a unique
$N$-horoball which intersects $P$ nontrivially.  We denote this
$N$-horoball by $P^N$.

Any horoball $P$ has a single accumulation point in $\partial X$,
which we denote $e_P$.  Denote by $\Hbound$ the collection of such
accumulation points.

\begin{definition} \label{d:geodAB}
Recall that in Lemma \ref{l:geodesicbicombing} we chose a
$G$-equivariant antisymmetric geodesic bicombing $\gamma$.  Let
$\Omega(X)$ denote the set of unit-speed geodesic paths in $X$, up to
orientation preserving reparametrization.
We will extend $\gamma$ to an antisymmetric map from most of
$(X\cup\mc{H}\cup\Hbound)^2$ to $\Omega(X)$.  (If $A$ is a horoball,
we leave $\gamma(A,A)$, $\gamma(e_A,A)$, $\gamma(A,e_A)$ and
$\gamma(e_A,e_A)$ undefined.)

First, for each pair of horoballs $A,B$, we choose (in an antisymmetric 
and $G$-equivariant way), a geodesic $\geod(A,B)$, which realizes the
distance between $A$ and $B$.   Second, for each point $a \in X$ and
each horoball $A\in\mc{H}$, we (equivariantly) choose a geodesic path
$\gamma(a,A)$ which realizes the distance from $a$ to $A$.
The path $\gamma(A,a)$ is the time-reverse of $\gamma(a,A)$.

Third, we extend $\gamma$ to points in $\Hbound$.  
If $x \in X$ and $\e_A \in \Hbound$ then we 
define $\geod(x,\e_A)$ to be the concatenation
of $\geod(x,A)$ with the vertical geodesic ray in $A$
from the endpoint of $\geod(x,A)$ to $\e_A$.

Finally, if $\e_A, \e_B \in \Hbound$ correspond
to distinct horoballs $A$ and $B$ in $\mc{H}$ then 
we define $\geod(\e_A, \e_B)$ to be the path 
$\geod(A,B)$, together with vertical paths on either end.
Note that $\geod(\e_A,\e_B)$ is always a 
geodesic line.
\end{definition}

\begin{torsionremark} \label{r:tf2}
As noted in Remark \ref{r:tf1}, $2$-torsion may prevent a choice
of paths as in Definition \ref{d:geodAB} above.

Also, if there is torsion, then parabolic subgroups of $G$ could intersect
nontrivially, which makes the $G$-equivariance problematic.

For this section, these problems can be solved by considering the
union of all shortest paths between $A$ and $B$ (of which there
are only finitely many, for a given $A$ and $B$).  In future
sections, this is a more subtle problem: in fact we should also use the ``average" as well as the union.

Similarly, we should take the union (or average) of all geodesics
from a point $a$ to a horoball $A$ which realize the distance from
$a$ to $A$.
\end{torsionremark}

\subsection{Definition and basic properties of preferred
paths}

In this subsection we apply the construction of Section
\ref{s:horoballs}.  The family $\mc{H}$ satisfies
the hypotheses of Theorem \ref{t:ppw}, therefore
there is a function
$\CHH \co X \times X \to \mc{O}_{\mc{H}}$ satisfying 
Axioms (A1)--(A7).

Suppose that $a,b \in X \cup \Hbound$.  We now associate
a collection of horoballs $\Hpairs_{a,b}$ to the pair
$\{ a,b \}$.  If $a,b \in X$ then $\Hpairs_{a,b} = \CHH_{a,b}$.
If $a\in X$ and $b =e_B \in \Hbound$ then $\Hpairs_{a,b}
= \CHH_{a,B}$ as defined in Notation \ref{n:horofamilies}.
Similarly, if $a = e_A \in \Hbound$ and $b \in X$ then
$\Hpairs_{a,b} = \CHH_{A,b}$ and if $a = e_A, b = e_B \in \Hbound$
then $\Hpairs_{a,b} = \CHH_{A,B}$.

\begin{remark}
For any $a,b \in X \cup \Hbound$,  the set $\Hpairs_{a,b}$
is a linearly ordered collection of horoballs, and this
order is compatible with the order obtained from
projection to $\geod(a,b)$.
\end{remark}

\begin{definition}\label{d:sigma}
For each $A \in \mc{H}$, and each pair $x,y \in A
\cup \{ \e_A \}$, we define $\sigma (x,y) = \sigma(y,x)$
in the following manner:  If $x = \e_A$, then 
$\sigma (x,y)$ is the vertical ray from $y$ to $\e_A$.
Thus suppose that $x,y \in A$. Then for some $t \in T$
we have $x = (t,p_1,k_1)$ and $y = (t,p_2,k_2)$ .
In case $p_1 = p_2$, then $\sigma (x,y)$ is the
vertical geodesic between $x$ and $y$.
Otherwise, suppose that $d_P(p_1,p_2)$ satisfies
\[	2^{N-1} < d_P(p_1,p_2) \leq 2^N	,	\]
for $N \in \N$, where $d_P(p_1,p_2)$ is the distance
between $p_1$ and $p_2$ in $P$, with respect to the generating
set $S \cap P$.  Then define
$R' = \max \{ N, k_1, k_2 \}$ and $R = R'$, unless
$R = L_2$ in which case $R = L_2+1$.
Then define $\sigma(x,y)$ to be the path which
consists of vertical paths from $x$ and $y$ to
depth $R$, and then joins the endpoints of
these vertical paths with the (unique) edge of length
$1$ at depth $R$.  
\end{definition}

Let $m$ be the midpoint of the unique horizontal edge of
$\sigma(x,y)$.  The point $m$ cuts $\sigma(x,y)$ into two geodesic
segments. 
The path $\sigma(x,y)$ is not usually itself
a geodesic, but is very close to
one.

\begin{remark}
We insist that the paths $\sigma(x,y)$ do not have a horizontal
edge at depth $L_2$ in order to simplify some future arguments.
In particular, this means that our preferred paths are `transverse'
to $D^{-1}(L_2)$.
\end{remark}

\begin{lemma}\label{l:sigma}
Let $x$ and $y$ lie in the same horoball $A$.
The path $\sigma(x,y)$ is Hausdorff distance at
most $5$ from $\geod (x,y)$.
\end{lemma}
\begin{proof}
The geodesic $\geod(x,y)$ 
coincides with
$\sigma(x,y)$ to depth $S$, where $S$ is the maximum depth of
$\geod(x,y)$.
Moreover, the sub-path $\sigma(x,y)\cap A^S$ has
length at most $7$,
while the (horizontal) sub-path $\geod(x,y)\cap A^S$ has length
at most $3$ (see Lemma \ref{l:gamma}).  Thus $\sigma(x,y)$ and
$\geod(x,y)$ are Hausdorff distance at most $5$ from
one another.
\end{proof}

We now want to use the collections of horoballs $\Hpairs$
to define {\em preferred paths} between any two points in
$X \cup \Hbound$.

Let $\alpha, \beta \in X \cup \Hbound$, and let 
$\Hpairs_{\alpha,\beta} = \{ A_1, \ldots , A_k \}$.
Note that $\Hpairs_{\alpha,\beta}$ is a finite 
linearly ordered set.  

\begin{definition}\label{d:preferredpaths}
{\em Preferred paths.}
If $\Hpairs_{\alpha,\beta}$ is empty, then we define
\[\pp_{\alpha,\beta}=\geod(\alpha,\beta).\]  
If $\Hpairs_{\alpha,\beta}$
contains a single horoball $A$, then we set
\[\pp_{\alpha,\beta} = \geod(\alpha,A)\cup \sigma(a_1,a_2)\cup\geod(A,\beta),\]
where $a_1$ is the terminal point of $\geod(\alpha,A)$, and $a_2$ is
the terminal point of $\geod(\beta,A)$.
Otherwise, we define the preferred path $\pp_{\alpha,\beta}$
to be 
\[ \geod (\alpha,A_1) \cup r_1 \cup 
\geod (A_1,A_2) \cup  r_2 \cup \cdots
\cup \geod (A_{k-1},A_k) \cup r_k
\cup \geod (A_k,\beta) ,	\]
where the $r_i$ are paths 
$\sigma(a_{i,1},a_{i,2})$, where $a_{1,1}$ is the terminal point of 
$\geod(\alpha, A_1)$ and otherwise
$a_{i,1}$ is the terminal point of
$\geod(A_{i-1}, A_i)$;  the point $a_{2,k}$ is the terminal point of
$\geod(A_k,\beta)$ and otherwise $a_{i,2}$ is the initial point of
$\geod(A_i, A_{i+1})$.

This describes $\pp_{\alpha,\beta}$ as a set, but we will often
consider it as a map, and parametrize by arc length.
\end{definition}
\begin{figure}[htbp]
\begin{center}
\begin{picture}(0,0)%
\includegraphics{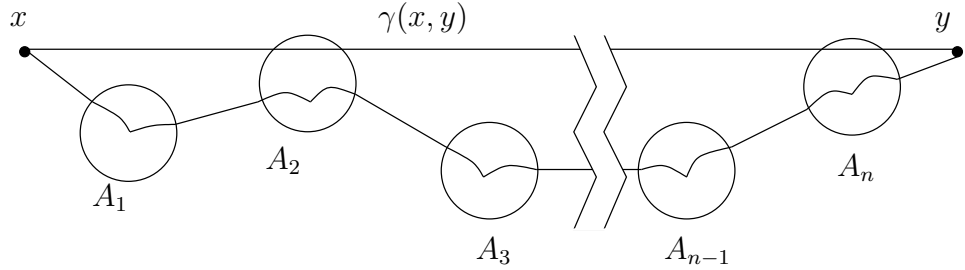}%
\end{picture}%
\setlength{\unitlength}{4144sp}%
\begingroup\makeatletter\ifx\SetFigFont\undefined%
\gdef\SetFigFont#1#2#3#4#5{%
  \reset@font\fontsize{#1}{#2pt}%
  \fontfamily{#3}\fontseries{#4}\fontshape{#5}%
  \selectfont}%
\fi\endgroup%
\begin{picture}(5723,1639)(1426,-4274)
\put(3646,-2806){\makebox(0,0)[lb]{\smash{{\SetFigFont{12}{14.4}{\familydefault}{\mddefault}{\updefault}{\color[rgb]{0,0,0}$\geod(x,y)$}%
}}}}
\put(1441,-2806){\makebox(0,0)[lb]{\smash{{\SetFigFont{12}{14.4}{\familydefault}{\mddefault}{\updefault}{\color[rgb]{0,0,0}$x$}%
}}}}
\put(6976,-2806){\makebox(0,0)[lb]{\smash{{\SetFigFont{12}{14.4}{\familydefault}{\mddefault}{\updefault}{\color[rgb]{0,0,0}$y$}%
}}}}
\put(4231,-4201){\makebox(0,0)[lb]{\smash{{\SetFigFont{12}{14.4}{\familydefault}{\mddefault}{\updefault}{\color[rgb]{0,0,0}$A_3$}%
}}}}
\put(5356,-4201){\makebox(0,0)[lb]{\smash{{\SetFigFont{12}{14.4}{\familydefault}{\mddefault}{\updefault}{\color[rgb]{0,0,0}$A_{n-1}$}%
}}}}
\put(6391,-3706){\makebox(0,0)[lb]{\smash{{\SetFigFont{12}{14.4}{\familydefault}{\mddefault}{\updefault}{\color[rgb]{0,0,0}$A_n$}%
}}}}
\put(1936,-3886){\makebox(0,0)[lb]{\smash{{\SetFigFont{12}{14.4}{\familydefault}{\mddefault}{\updefault}{\color[rgb]{0,0,0}$A_1$}%
}}}}
\put(2971,-3661){\makebox(0,0)[lb]{\smash{{\SetFigFont{12}{14.4}{\familydefault}{\mddefault}{\updefault}{\color[rgb]{0,0,0}$A_2$}%
}}}}
\end{picture}%
\caption{A preferred path.}
\label{f:prefschem}
\end{center}
\end{figure}

\begin{torsionremark} \label{r:tf3}
We have indicated in Remarks \ref{r:tf1} and \ref{r:tf2} how
the paths $\geod(x,y)$ and $\geod(A,B)$ might be defined
using averages in order to ensure antisymmetry and
$G$-equivariance.

In this framework, preferred paths will consist of all
of the ``possible" preferred paths,
superimposed like the states of a quantum system.

These more complicated preferred paths will continue to have many of the properties that honest preferred paths do, though these properties
become a deal more cumbersome to state, let alone prove.  
\end{torsionremark}

\begin{definition}
Let $\ppbar_{x,y}$ be the path obtained from $\pp_{x,y}$ by replacing
each subsegment of the form $\sigma(a_{i,1},a_{i,2})$ by the corresponding
geodesic $\geod(a_{i,1},a_{i,2})$.
\end{definition}

\begin{lemma}\label{l:localgeodesic}
If $A_i\in \Hpairs_{x,y}$, then either $\ppbar_{x,y}\cap A^{L}$ is geodesic
for any $L\geq 2\delta$, or $d(a_{i,1},a_{i,2})<3$, where $a_{i,1}$
and $a_{i,2}$ are the first and last points on $\pp_{x,y}\cap A_i$.
\end{lemma}
\begin{proof}
Let $\ppbar_{x,y} \cap A \cap D^{-1}(2\delta) = \{ b_1, b_2 \}$ and
let $\gamma = \geod_{b_1,b_2}$.  This consists of two vertical
segments, and a single horizontal segment of length at most
$3$.  

If $\gamma$ intersects $D^{-1}(L_1)$ nontrivially, then the vertical
segments in $\gamma$ and the vertical segments in $\ppbar_{x,y}
\cap A^L$ coincide, and the horizontal segments have the
same length.  Therefore, in this case, $\ppbar_{x,y} \cap A^L$ is
a geodesic.

Otherwise let $c_1$ and $c_2$ be the deepest points of the 
vertical subsegments of $\gamma$, with $c_1$ lying directly
above $a_{i,1}$ and $c_2$ directly above $a_{i,2}$.  Then
the distance between $c_1$ and $c_2$ is at most $3$, which
implies that the distance between $a_{i,1}$ and $a_{i,2}$ is
less than $3$.
\end{proof}

\subsubsection{Parametrizations} \label{ss:param}
We pause briefly to discuss the parametrizations of
the paths $\geod(x,y)$, $\pp_{x,y}$ and $\ppbar_{x,y}$.
These parametrizations are always by arc length. For
each $x, y \in X \cup \Hbound$ choose
$\mc{I}_{x,y}$, $\bar{\mc{I}}_{x,y}$ and $\mc{J}_{x,y}$ as follows:
If $x,y \in X$ then
$\mc{I}_{x,y} = [0,\mbox{length}(\pp_{x,y})]$, and similarly
$\bar{\mc{I}}_{x,y} = [0 ,\mbox{length}(\ppbar_{x,y})]$ and
$\mc{J}_{x,y} = [0,d(x,y)]$.

If $x \in X$ and $y \in \Hbound$, then all of the intervals
are $[0,\infty)$.

If $x \in \Hbound$, then each of the paths $\geod(x,y)$,
$\pp_{x,y}$ and $\ppbar_{x,y}$ begin with an infinite
vertical ray.  If $y$ is contained in the same $L_1$-horoball
as $x$ then the three paths coincide and we define
$\mc{I}_{x,y} = \bar{\mc{I}}_{x,y} = \mc{J}_{x,y} = 
(-\infty,L_1 - D(y)]$.

Suppose then that $x \in \Hbound$ and that $y$ is not contained
in the same $L_1$-horoball as $x$.  Let $x'$ be the first point
on $\pp_{x,y}$ contained in $D^{-1}(L_1)$.  Then
the interval $\mc{I}_{x',y}$ is already defined, and we define
$\mc{I}_{x,y} = (\infty,0] \cup \mc{I}_{x',y}$, and suppose that
$\pp_{x,y}(0) = x'$.  

We make similar definitions for
$\bar{\mc{I}}_{x,y}$ and $\mc{J}_{x,y}$.

\begin{proposition}\label{p:async}
For any two points $x$ and $y$ in $X \cup \Hbound$, 
there are nondecreasing functions
\[f=f_{xy}\co \mc{I}_{x,y} \to \mc{J}_{x,y} \]
and 
\[g=g_{xy}\co \mc{J}_{x,y} \to \mc{I}_{x,y} \]
satisfying:
\begin{enumerate}
\item For all but finitely many points
  $t\in \mc{I}_{x,y}$ and $s\in \mc{J}_{x,y}$, $f'(t)$
  and $g'(s)$ exist and are either $0$ or $1$.
\item\label{displace} For any $t\in \mc{I}_{x,y}$, and
$s\in \mc{J}_{x,y}$ we have
\[ d(\pp_{x,y}(t), \geod(x,y)(f(t))\leq \Displace  \mbox{, and} \] 
\[ d(\pp_{x,y}(g(s)), \geod(x,y)(s)\leq \Displace.\]
\item If $\pp_{x,y}(t)$ is outside of $A^{L_1-(\Hexradius)}$ for every
  $A\in\Hpairs_{x,y}$, then $f'(t)=1$.

  Similarly, if $\geod(x,y)(t)$ is outside of $A^{L_1-(\Hexradius)}$ for every
  $A\in\Hpairs_{x,y}$, then $g'(t)=1$.
\item\label{waitjump} Let $A\in\Hpairs_{x,y}$ and  
  $\hat{A}=A^{L_1-(\Hexradius)}$.  Let
  $[t_1,t_2]=\pp_{x,y}^{-1} (\hat{A})$, and 
  $[s_1,s_2]=\geod(x,y)^{-1}(\hat{A})$.   
  The total
  waiting times in $A$, i.e. the maximum of the measures of
  \[\{t_1\leq t\leq t_2\mid f'(t)=0 \}\mbox{ and } 
  \{s_1\leq s\leq s_2\mid g'(s)=0 \}\] 
  is at most $\Waitbound$.
  Moreover, the  total jumping times in $A$, i.e. the maximum of the
  measures of
  \[\{s_1\leq s\leq s_2\mid \not\exists t\mbox{, } f(t)=s\}\mbox{ and }
  \{t_1\leq t\leq t_2\mid \not\exists s\mbox{, } g(s)=t\}\] 
  is at most $\Jumpbound$.
\end{enumerate}
\end{proposition}
\begin{proof}
We suppose for the moment that $x,y \in D^{-1}[0,L_1]$.  We will
deal with the other cases at the end of the proof.  In fact, we only
deal with the cases that $x,y \in D^{-1}[0,L_1] \cup \Hbound$.  The
other cases are similar, but are not required for later applications.

We use the same symbol for a path and its image throughout.
Note that a map from the domain of one path to the domain of
another induces a map from the image of one path to the image of
the second.  Conversely, if both paths are embeddings, then a map
between their images induces a map between their domains.

Recall that the path $\pp_{x,y}$ is defined to be a concatenation
of geodesic segments outside of horoballs, and paths of the form
$\sigma(u,v)$ through horoballs.

Consider the difference between the paths $\sigma(u,v)$ and 
$\geod(u,v)$.  This consists of a path of length at most 
$9$ in $\sigma(u,v)$ and a path of length at most $5$ in $\geod(u,v)$.
Define functions 
\begin{eqnarray*}
{\mc{I}_{x,y}} & \xrightarrow{f_1} & \bar{\mc{I}}_{x,y} , \mbox{ and}\\
\bar{\mc{I}}_{x,y}  &\xrightarrow{g_1} & \mc{I}_{x,y},
\end{eqnarray*}
so that the induced maps on the intersection of the images are
the identity, and the induced maps on the differences send a
component to its initial point.

The functions $f_1$ and $g_1$ have derivatives at all but finitely 
many points,
the derivatives are always $0$ or $1$, and the waiting times and
jumping distances are at most $9$ per horoball in $\mc{H}_{x,y}$.

It is not difficult to see that it suffices now to define appropriate
functions
\begin{eqnarray*}
{\bar{\mc{I}}_{x,y}}  & \xrightarrow{f_2}& \mc{J}_{x,y} \mbox{ and } \\
{\mc{J}_{x,y}} & \xrightarrow{g_2} & \bar{\mc{I}}_{x,y}  .
\end{eqnarray*}

There are three cases, depending on whether $\Hpairs_{x,y}$ is empty,
contains a single horoball, or contains at least two horoballs.

{\bf Case 1:}
If $\Hpairs_{x,y}$ is empty, then $\pp_{x,y}=\bar\pp_{x,y} = \geod(x,y)$, and we can
set both $f$ and $g$ equal to the identity map.

{\bf Case 2:}
Suppose $\Hpairs_{x,y}$ is a single horoball.  This case is 
strictly easier than Case 3, and we leave it as an exercise.

{\bf Case 3:} There is more than one horoball in $\Hpairs_{x,y}$.
By Axiom (A2), the intersection $\geod(x,y)\cap A_i^{L_1-K}$
is non-empty for each $1 \le i \le n$.

Let $u$ be the first point on $\geod(x,y)$ which intersects $A_1^{L_1-K}$,
and let $v$ be the first point on $\bar\pp_{x,y}$ which intersects
$A_1^{L_1-K}$.
We define $f_2$ and $g_2$ to be the identity on
$[0, (u,v)_x]$.  Note that the segments 
$\bar\pp_{x,y}[0, (u,v)_x]$ and 
$\geod(x,y)[0,(u,v)_x]$ have the same length and are 
identified under the map to the comparison tripod $Y_{u,v,x}$.
Thus, for example, the distance between
$\bar\pp_{x,y}(t)$ and $\geod(x,y)(f_2(t))$ is at most $\delta$ for
all $t \in [0, (u,v)_x]$.  There is clearly no waiting or
jumping in this interval.  

Similarly, let $w$ be the last point on $\geod(x,y)$ intersecting $A_n$,
and let $z$ be the last point on $\bar\pp_{x,y}$ intersecting $A_n$.
We can define the map 
\[	f_2 \co \Big[ |\bar\pp_{x,y}| - (y_1,y_2)_y, |\bar\pp_{x,y}| \Big] 
\to \Big[ d(x,y) - (w,z)_y, d(x,y) \Big]	,\]
to be the unique (orientation preserving) isometry, and 
$g_2 = (f_2)^{-1}$ on this interval.  Once again, the images
of these intervals are identified under the map to the comparison
tripod $Y_{w,z,y}$.  We also define 
$y_{1,n} = \geod(x,y)(d(x,y) - (w,z)_y)$ and 
$y_{2,n} = \bar\pp_{x,y}(d(x,y) - (w,z)_y)$.

We now define the maps $f_2$ and $g_2$ on the region between
adjacent horoballs $A_i$ and $A_{i+1}$.

Let $1 \le i \le n-1$.  Let $a$ be the last point on $\geod(x,y) \cap 
A_i^{L_1-K}$,
and let $b$ be the first point on $\geod(x,y) \cap A_{i+1}^{L_1-K}$.  
Similarly, let $c$ be the last point on $\bar\pp_{x,y} \cap A_i^{L_1-K}$
and let $d$ be the first point on $\bar\pp_{x,y} \cap A_{i+1}^{L_1-K}$. 
Suppose that
$a = \geod(x,y)(t_a)$, and define times $t_b, t_c$ and $t_d$ 
analogously.  By drawing two comparison tripods (one for
the triangle $\Delta(a,b,c)$ and one for $\Delta(b,c,d)$) we see that 
\[	I_{a,b} = \Big[ t_a + (b,c)_a, t_b - (c,d)_b \Big], \mbox{ and } \]
\[  	I_{c,d} = \Big[ t_c + (a,b)_c, t_d - (b,c)_d \Big]	,	\]
have the same length (note that between $A_i^{L_1-K}$ and 
$A_{i+1}^{L_1-K}$,
both $\bar\pp_{x,y}$ and $\geod(x,y)$ are geodesics 
parametrized by arc length).

Define $f_2 \co I_{c,d} \to I_{a,b}$ and $g_2 \co I_{a,b} \to I_{c,d}$
to be orientation preserving isometries.  Note that the comparison
tripods imply that for $t \in I_{a,b}$, 
\[	d(\geod(x,y)(t), \bar\pp_{x,y}(g_2(t))) \le 2\delta,	\] 
and for $s\in I_{c,d}$,
\[	d(\bar\pp_{x,y}(s), \geod(x,y)(f_2(s))) \le 2\delta.	\]

We claim now that those parts of 
$\bar{\mc{I}}_{x,y}$ on which $f_2$ is not
yet defined are segments $I_1, \ldots , I_n$ so that
$\bar\pp_{x,y}(I_j) \subseteq N_{2\delta}(A_j^{L_1-K})$, and
similarly that those parts of
$\mc{J}_{x,y}$ on which $g_2$ is not yet defined are
segments $I_1', \ldots , I_n'$ so that
$\geod(x,y)(I_j') \subseteq N_{2\delta}(A_j^{L_1-K})$.
To see this consider, for example, the points
$\geod(x,y)(t_a+(b,c)_a)$ and $\bar\pp_{x,y}(t_c+(a,b)_c)$, for the
points $a,b,c$ and $d$ as above.  These are
typical endpoints of the intervals on which $f_2$ and $g_2$
are not yet defined.
Any geodesic between $a$ and $c$ lies entirely within 
the convex set $A_i^{L_1-K}$.
Therefore, the point $\geod(x,y)(t_a + (b,c)_a)$ lies within
$\delta$ of $A_i^{L_1-K}$, because it is in the preimage of the central point
of the comparison tripod $Y_{a,b,c}$.  Now, $\bar\pp_{x,y}(t_c +
(a,b)_c)$ lies within $\delta$ of the point on the geodesic $[c,b]$ 
in the preimage of the central point of $Y_{a,b,c}$, and in particular,
$\bar\pp_{x,y}(t_c + (a,b)_c)$ lies within $2\delta$ of $A_i^{L_1-K}$.
Similar arguments for all other points prove the claim.

We now define $f_2$ and $g_2$ on the remaining subsegments,
$I_j$ and $I_j'$.  We fix $1 \le j \le n$.  Let $x_1$ and $y_1$ be the
images of the endpoints of $I_j'$ under $\geod(x,y)$, and let
$x_2$ and $y_2$ be the images of the endpoints of $I_j$ under
$\bar\pp_{x,y}$.  The path $\bar\pp_{x,y}(I_j)$ is a concatenation
of three geodesics, one traveling towards $A_j$ as quickly as possible,
one through $A_j$ and one traveling away from $A_j$ as quickly
as possible.  Let $x_3$ be the first point on $\bar\pp_{x,y}(I_j) \cap
A_j$ and let $y_3$ be the last point on $\bar\pp_{x,y}(I_j) \cap A_j$.
(See Figure \ref{f:waitjump1}.)

It is clear from the previous arguments that $d(x_1,x_2), d(y_1,y_2)
\le 2\delta$.  Also, $x_2$ is within $K+2\delta$ of $A_j$.  However, between $x_2$ and $x_3$, the 
path $\bar\pp_{x,y}$ is traveling as quickly as possible towards
$A_j$, so $d(x_2,x_3) \le K+2\delta$.  Similarly, $d(y_2,y_3) \le
K + 2\delta$.

Suppose that $x_1 = \geod(x,y)(s_1), y_1 = \geod(x,y)(s_2),
x_2 = \bar\pp_{x,y}(t_1), x_3 = \bar\pp_{x,y}(t_2), 
y_3 = \bar\pp_{x,y}(t_3)$ and $y_2 = \bar\pp_{x,y}(t_4)$.

For $t \in [t_1,t_2]$, define $f_2(t) = s_1$ , and for $t \in [t_3,t_4]$
define $f_2(t) = s_2$.

Drawing comparison tripods as in Figure \ref{f:waitjump2} allows us
to see that
\[	J_1 = \Big[ s_1 + (x_3,y_1)_{x_1}, s_2 - (x_3,y_3)_{y_1} \Big],
\mbox{ and }	\]
\[ J_2 = \Big[ t_2 + (x_1,y_1)_{x_3}, t_3 - (x_3,y_1)_{y_3} \Big]	\]
have the same length.  

Assume first that $J_1$ and $J_2$ are nonempty,
and define $f_2 \co J_2 \to J_1$ and
$g_2 \co J_1 \to J_2$ to be orientation preserving isometries.

\begin{figure}[htbp]
\begin{center}
\begin{picture}(0,0)%
\includegraphics{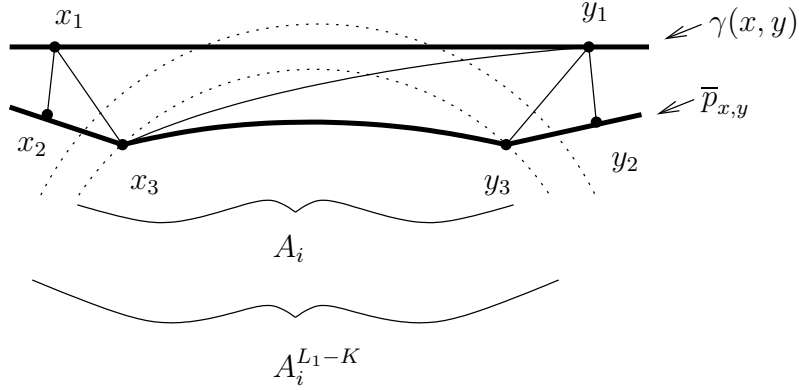}%
\end{picture}%
\setlength{\unitlength}{4144sp}%
\begingroup\makeatletter\ifx\SetFigFont\undefined%
\gdef\SetFigFont#1#2#3#4#5{%
  \reset@font\fontsize{#1}{#2pt}%
  \fontfamily{#3}\fontseries{#4}\fontshape{#5}%
  \selectfont}%
\fi\endgroup%
\begin{picture}(4233,2404)(1858,-5444)
\put(3466,-4651){\makebox(0,0)[lb]{\smash{{\SetFigFont{12}{14.4}{\familydefault}{\mddefault}{\updefault}{\color[rgb]{0,0,0}$A_i$}%
}}}}
\put(3466,-5371){\makebox(0,0)[lb]{\smash{{\SetFigFont{12}{14.4}{\familydefault}{\mddefault}{\updefault}{\color[rgb]{0,0,0}$A_i^{L_1-K}$}%
}}}}
\put(6031,-3751){\makebox(0,0)[lb]{\smash{{\SetFigFont{12}{14.4}{\familydefault}{\mddefault}{\updefault}{\color[rgb]{0,0,0}$\bar\pp_{x,y}$}%
}}}}
\put(6076,-3301){\makebox(0,0)[lb]{\smash{{\SetFigFont{12}{14.4}{\familydefault}{\mddefault}{\updefault}{\color[rgb]{0,0,0}$\geod(x,y)$}%
}}}}
\put(2161,-3256){\makebox(0,0)[lb]{\smash{{\SetFigFont{12}{14.4}{\familydefault}{\mddefault}{\updefault}{\color[rgb]{0,0,0}$x_1$}%
}}}}
\put(5311,-3211){\makebox(0,0)[lb]{\smash{{\SetFigFont{12}{14.4}{\familydefault}{\mddefault}{\updefault}{\color[rgb]{0,0,0}$y_1$}%
}}}}
\put(2611,-4246){\makebox(0,0)[lb]{\smash{{\SetFigFont{12}{14.4}{\familydefault}{\mddefault}{\updefault}{\color[rgb]{0,0,0}$x_3$}%
}}}}
\put(1936,-4021){\makebox(0,0)[lb]{\smash{{\SetFigFont{12}{14.4}{\familydefault}{\mddefault}{\updefault}{\color[rgb]{0,0,0}$x_2$}%
}}}}
\put(4726,-4246){\makebox(0,0)[lb]{\smash{{\SetFigFont{12}{14.4}{\familydefault}{\mddefault}{\updefault}{\color[rgb]{0,0,0}$y_3$}%
}}}}
\put(5491,-4111){\makebox(0,0)[lb]{\smash{{\SetFigFont{12}{14.4}{\familydefault}{\mddefault}{\updefault}{\color[rgb]{0,0,0}$y_2$}%
}}}}
\end{picture}%
\caption{$\gamma(x,y)$ and $\bar\pp_{x,y}$ near $A_i\in\mc{H}_{x,y}$.}
\label{f:waitjump1}
\end{center}
\end{figure}
\begin{figure}[htbp]
\begin{center}
\begin{picture}(0,0)%
\includegraphics{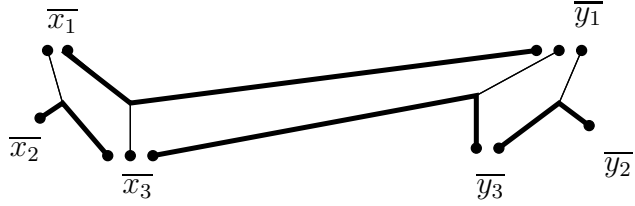}%
\end{picture}%
\setlength{\unitlength}{4144sp}%
\begingroup\makeatletter\ifx\SetFigFont\undefined%
\gdef\SetFigFont#1#2#3#4#5{%
  \reset@font\fontsize{#1}{#2pt}%
  \fontfamily{#3}\fontseries{#4}\fontshape{#5}%
  \selectfont}%
\fi\endgroup%
\begin{picture}(3585,1279)(1921,-4319)
\put(2161,-3256){\makebox(0,0)[lb]{\smash{{\SetFigFont{12}{14.4}{\familydefault}{\mddefault}{\updefault}{\color[rgb]{0,0,0}$\bar{x_1}$}%
}}}}
\put(5311,-3211){\makebox(0,0)[lb]{\smash{{\SetFigFont{12}{14.4}{\familydefault}{\mddefault}{\updefault}{\color[rgb]{0,0,0}$\bar{y_1}$}%
}}}}
\put(5491,-4111){\makebox(0,0)[lb]{\smash{{\SetFigFont{12}{14.4}{\familydefault}{\mddefault}{\updefault}{\color[rgb]{0,0,0}$\bar{y_2}$}%
}}}}
\put(4726,-4246){\makebox(0,0)[lb]{\smash{{\SetFigFont{12}{14.4}{\familydefault}{\mddefault}{\updefault}{\color[rgb]{0,0,0}$\bar{y_3}$}%
}}}}
\put(2611,-4246){\makebox(0,0)[lb]{\smash{{\SetFigFont{12}{14.4}{\familydefault}{\mddefault}{\updefault}{\color[rgb]{0,0,0}$\bar{x_3}$}%
}}}}
\put(1936,-4021){\makebox(0,0)[lb]{\smash{{\SetFigFont{12}{14.4}{\familydefault}{\mddefault}{\updefault}{\color[rgb]{0,0,0}$\bar{x_2}$}%
}}}}
\end{picture}%
\caption{Comparison tripods for the triangles in Figure \ref{f:waitjump1}.}
\label{f:waitjump2}
\end{center}
\end{figure}

It remains to define $f_2$ on the intervals $K_1 = \Big[ t_2, t_2 + 
(x_1,y_1)_{x_3} \Big]$ and $K_2 = \Big[ t_3 - (x_3,y_1)_{y_3}, t_3
\Big]$, and $g_2$ on the intervals $K_3 = \Big[s_1,s_1 + 
(x_3,y_1)_{x_1} \Big]$ and $K_4 = \Big[ s_2 - (x_3,y_3)_{y_1},
s_2 \Big]$.  Since $d(x_1,x_2)\le 2\delta$ and $d(x_2,x_3) \le K+2\delta$, we have $d(x_1,x_3)\le K+4\delta$.  This implies that the lengths
of $K_1$ and $K_3$ are at most $K+4\delta$.  The same is true
of the lengths of $K_2$ and $K_4$.

If $t \in K_1$ then define $f_2(t) = s_1 + (x_3,y_1)_{x_1}$.  If
$t \in K_2$ then define $f_2(t) = s_2 - (x_3,y_3)_{y_1}$.
If $s \in K_3$ then define $g_2(s) = t_2 + (x_1,y_1)_{x_3}$.
If $s \in K_4$ then define $g_2(s) = t_3 - (x_3,y_1)_{y_3}$.

If
$J_1$ and $J_2$ are empty, we must define $f_2$ and $g_2$ slightly
differently:
\begin{equation*}
\begin{split}
f_2(\Big[t_2,t_2+(y_1,y_3)_{x_3}\Big]) & =s_1\\
f_2(\Big[t_2+(y_1,y_3)_{x_3},t_4\Big]) & =s_2\\ 
g_2(\Big[s_1,s_1+(y_1,x_3)_{x_1}\Big]) & =t_2\\
g_2(\Big[s_1+(y_1,x_3)_{x_1},s_2\Big]) & =t_3
\end{split}
\end{equation*}
In this case, the intersections of $\gamma(x,y)$ and $\bar\pp_{x,y}$
with $A_i^{L_1-K}$ are quite small.

We have now defined $f_2$ and $g_2$ in the case that
$x,y \in D^{-1}[0,L_2]$.
It is not difficult to see that the functions $f  = f_2 \circ f_1$ and
$g = g_1 \circ g_2$ satisfy the requirements of the proposition.
In fact (and this is used below), we have shown that they
satisfy the requirements of the proposition with the constants
reduced by $6\delta$.

Finally, we now turn to the case that one or both of $x,y$ lie
in $\Hbound$.  

Suppose $x \in \Hbound$ and $y \in D^{-1}[0,L_1]$.  As in the 
discussion above the proposition, let $x'$ be the first point on
$\pp_{x,y}$ so that $x' \in D^{-1}(L_1)$.  Let $x''$ be the first
point on $\geod_{x,y}$ so that $x'' \in D^{-1}(L_1)$.  Then
$\mc{I}_{x,y} = (-\infty,0] \cup \mc{I}_{x',y}$ and
$\mc{J}_{x,y} = (-\infty,0] \cup \mc{J}_{x'',y}$.  Note also
that $\pp_{x',y} = \pp_{x,y}|_{\mc{I}_{x',y}}$, by the definition
of preferred paths.

We define $f_{xy}$ to be the identity on $(-\infty,0]$.
If we can prove that $d(x',x'') \le 6\delta$, then there are obvious
maps $\phi \co \mc{J}_{x',y} \to \mc{J}_{x'',y}$ and
$\psi \co \mc{J}_{x'',y} \to \mc{J}_{x',y}$ which identify the long 
side of the comparison tripod $Y_{x,x',y}$, 
and collapse the short sides to the center of the tripod.  If we have
such functions $\phi$ and $\psi$, then we define 
\[f_{x,y}|_{\mc{I}_{x',y}} =  \phi \circ f_{x',y} \co \mc{I}_{x',y} \to
\mc{J}_{x'',y},	\]
and
\[	g_{x,y}|_{\mc{J}_{x'',y}} = g_{x',y} \circ \psi \co \mc{J}_{x'',y}
\to \mc{I}_{x',y}	.	\]
These will satisfy the
requirements of the proposition so long as $d(x',x'') \le 6\delta$.

To see this, we argue as follows:  Let $\sigma$ be that part
of $\pp_{x,y}$ from $x$ until the second horoball in $\mc{H}_{x,y}$
(or to be all of $\pp_{x,y}$ if $\mc{H}_{x,y}$ is a singleton).  It is not
difficult to see that $\sigma$ must be a geodesic, because that
part after $x'$ travels away from the horoball as quickly as 
possible (by the definition of preferred paths).  Let $z$ be the endpoint
of $\sigma$ which is not $x'$.  By Axiom (A2) of Section \ref{s:horoballs},
and Lemma  \ref{l:tube}, there is a point
$z' \in \geod(x,y)$ which lies within $K$ of the horoball, $A$, which contain $z$.  The geodesic between $z$ and $z'$ lies entirely within the
convex set $A^{L_1-K}$.
The (partially ideal) triangle with vertices $x,z,z'$ is $3\delta$-slim,
which proves that $d(x',x'') \le 6\delta$ (because  both $\pp_{x,y}$
and $\geod(x,y)$ are vertical on the interval $(-\infty,3\delta)$.

We now suppose that $x,y \in \Hbound$.  Then there are points
$x',y' \in \pp_{x,y}$ and $x'',y'' \in \geod(x,y)$ so that
\[	\mc{I}_{x,y} = (-\infty,0] \cup \mc{I}_{x',y'} \cup [|\pp_{x',y'}|,\infty),	\]
and 
\[	\mc{J}_{x,y} = (-\infty,0] \cup \mc{J}_{x'',y''} \cup
[d(x'',y''),\infty)	.	\]
As before, $d(x',x''), d(y'y,'') \le 6\delta$.  Thus there are obvious functions
$\phi' \co \mc{J}_{x',y'} \to \mc{J}_{x'',y''}$ and
$\psi' \co \mc{J}_{x'',y''} \to \mc{J}_{x',y'}$ and we make the
following definitions:
\begin{eqnarray*}
f_{xy}(t) = \left\{ \begin{array}{cc}
t, & \mbox{if $t \in (-\infty,0)$} \\
\phi'(f_{x',y'}(t)), & \mbox{if $t \in \mc{I}_{x',y'}$, and} \\
t + (d(x'',y'') - |\pp_{x',y'}|) & \mbox{otherwise,}\\
\end{array} \right.
\end{eqnarray*}
and
\begin{eqnarray*}
g_{xy}(t) = \left\{ \begin{array}{cc}
t, & \mbox{if $t \in (-\infty,0)$} \\
\psi'(f_{x',y'}(t)), & \mbox{if $t \in \mc{I}_{x',y'}$, and} \\
t + (|\pp_{x',y'}| - d(x'',y'')) & \mbox{otherwise.}\\
\end{array} \right.
\end{eqnarray*}

\end{proof} 

We now state a couple of corollaries of Proposition \ref{p:async}
\begin{corollary}\label{c:ppneargeodesic}
For any $a$, $b\in X\cup\Hbound$, the Hausdorff distance between
$\pp_{a,b}$ and $\geod(a,b)$ is at most $\Displace$.  The Hausdorff
distance between $\pp_{a,b}[t_1,t_2]$ and $\geod(a,b)[f(t_1),f(t_2)]$
is at most $\Hdist$.
\end{corollary}
\begin{proof}
The first assertion is obvious from Proposition
\ref{p:async}.\eqref{displace}.  Suppose that $t\in [t_1,t_2]$.
Since the function $f$ from Proposition \ref{p:async} is monotone,
there is some $s\in [f(t_1),f(t_2)]$ so that $f(t)=s$.  By Proposition
\ref{p:async}.\eqref{displace},
$d(\geod(a,b)(s),\pp_{a,b}(t)\leq\Displace$.  Conversely, suppose that
$s\in [f(t_1),f(t_2)]$.  It follows from monotonicity and Proposition
\ref{p:async}.\eqref{waitjump} that there is a point 
$s'\in f([t_1,t_2]$ with $|s-s'|\leq \frac{\Jumpbound}{2}$.  It
follows that $\geod(a,b)(s)$ is no further than
$\Displace+\frac{\Jumpbound}{2}\leq\Hdist$ from $\pp_{a,b}[t_1,t_2]$.
\end{proof}

\begin{corollary}\label{c:quasigeodesic}
The path $\pp_{x,y}$ is a $(\lambda,\epsilon)$-quasi-geodesic, for 
\[\lambda=2 \mbox{, and}\] 
\[\epsilon=\Pepsilon.\]
\end{corollary}
\begin{proof}
Let $t_1$ and $t_2$ be points in $[0,\mathrm{length}(\pp_{x,y})]$, and
let $x_1=\pp_{x,y}(t_1)$ and $x_2=\pp_{x,y}(t_2)$.
Since $\pp_{x,y}$ is parametrized by arc length, 
$|t_2-t_1|\geq d(x_1,x_2)$ is automatic.

Let $x_i'=\geod(x,y)(f(t_i))$, for $i\in\{1,2\}$ and $f$ as described
in Proposition \ref{p:async}.  Then $d(x_1',x_2')\leq
d(x_1,x_2)+2(\Displace)$. 
The arc length $|t_2-t_1|$ is at most
$|f(t_2)-f(t_1)|=d(x_1',x_2')$ plus the amount of time in $[t_1,t_2]$
that $f'(t)=0$.  The number of horoballs in $\Hpairs_{x,y}$ that
$\pp_{x,y}$ passes within $\Hexradius$ of between time $t_1$ and
$t_2$ is at most $2+\frac{|t_2-t_1|}{L_1}$; the amount of time
$f'(t)=0$ when $p_{x,y}(t)$ is near a horoball is at most
$\Waitbound<\frac{1}{2}L_1$, by part (4) of Proposition \ref{p:async}.  Thus 
\[
|t_2-t_1| \leq d(x_1',x_2') +
 (2+\frac{|t_2-t_1|}{L_1})(\Waitbound),\]
which implies
\begin{equation*}
\begin{split}
  \frac{1}{2}|t_2-t_1|& \leq(1-\frac{\Waitbound}{L_1})|t_2-t_1|\\
                      &\leq d(x_1',x_2')+2(\Waitbound),
\end{split}
\end{equation*}
and so 
\begin{equation*}
\begin{split}
|t_2-t_1|  &\leq 2 d(x_1',x_2')+4(\Waitbound)\\
 &\leq 2 d(x_1,x_2)+4(\Displace)+4(\Waitbound)\\
& = 2 d(x_1,x_2)+ \Pepsilon
\end{split}
\end{equation*}
\end{proof}

We now make some some other useful observations about 
preferred paths.  

\begin{lemma}
Preferred paths are $G$-equivariant.
\end{lemma}

\begin{remark}
It follows from Corollary \ref{c:quasigeodesic} that sub-paths of
preferred paths are quasi-geodesic, and hence form 
$\delta''$-slim triangles for \emph{some} $\delta''$.  However,
the $\delta''$ coming from $(\lambda,\epsilon)$-quasi-geodesicity is
{\em much} worse than the $\delta'$ below, and our previous
choice of $L_1$ secretly takes into account the particular bound given
in the next proposition.
\end{remark}

\begin{proposition}\label{p:deltaprime}
There exists a constant 
\[\delta' = \Deltaprime \le \frac{1}{4}L_1\] 
so that any triangle or bigon whose sides are 
sub-paths of preferred paths is $\delta'$-slim.

This includes triangles any or all of whose vertices are
in $\Hbound$.
\end{proposition}

\begin{proof}
Let $P$, $Q$, and $R$ be three points in $X\cup\Hbound$, joined by
parts of preferred paths $\pp_{ab}$, $\pp_{cd}$ and $\pp_{ef}$.
Specifically, suppose that $\pp_{ab}[t_1,t_2]$ joins $P$ to
$Q$, $\pp_{cd}[t_3,t_4]$ joins $Q$ to $R$, and
$\pp_{ef}[t_5,t_6]$ passes from $R$ to $P$, where $t_i\in
[-\infty,\infty)$ for $i$ odd and $t_i\in (-\infty, \infty]$ for $i$
even.  

By Corollary \ref{c:ppneargeodesic},
the side $\pp_{ab}[t_1,t_2]$ is Hausdorff distance at most $\Hdist$
from the geodesic segment $\geod(a,b)[f(t_1),f(t_2)]$, where
$f=f_{ab}$ is the function from Proposition \ref{p:async}.
Furthermore, the
endpoints $\geod(a,b)(f(t_1))$ and $\geod(a,b)(f(t_2))$ of this
segment are (if not ideal) at most $\Displace$ from $P$ and $Q$
respectively, by Proposition \ref{p:async}.\eqref{displace}.  Similar
statements can be made about the other two sides.  If none of the
endpoints of the geodesic segments are ideal, 
they may be joined together by segments of length
at most $2(\Displace)$, to form a geodesic hexagon, which must be
$4\delta$-slim.  Otherwise, 
they may be partially joined up to form a pentagon with one ideal
vertex, a quadrilateral with two ideal vertices, or simply an ideal
triangle.  Call the (possibly partially ideal) polygon so obtained $H$.
Subdividing and using Lemma \ref{l:idealthin}, it is easy
to see that $H$ is $6\delta$-slim.

Suppose that $x\in \pp_{ab}[t_1,t_2]$.  There is some point $x'$ on
$\geod(a,b)[f(t_1),f(t_2)]$ at most $\Hdist$ away from $x$.
 (see Figure
\ref{f:deltaprime}.)
\begin{figure}[htbp]
\begin{center}
\begin{picture}(0,0)%
\includegraphics{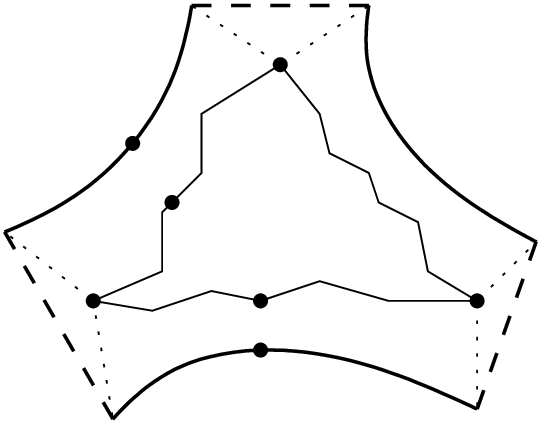}%
\end{picture}%
\setlength{\unitlength}{4144sp}%
\begingroup\makeatletter\ifx\SetFigFont\undefined%
\gdef\SetFigFont#1#2#3#4#5{%
  \reset@font\fontsize{#1}{#2pt}%
  \fontfamily{#3}\fontseries{#4}\fontshape{#5}%
  \selectfont}%
\fi\endgroup%
\begin{picture}(2474,1981)(564,-3235)
\put(1711,-2491){\makebox(0,0)[lb]{\smash{{\SetFigFont{12}{14.4}{\familydefault}{\mddefault}{\updefault}{\color[rgb]{0,0,0}$x$}%
}}}}
\put(1711,-3166){\makebox(0,0)[lb]{\smash{{\SetFigFont{12}{14.4}{\familydefault}{\mddefault}{\updefault}{\color[rgb]{0,0,0}$x'$}%
}}}}
\put(991,-1861){\makebox(0,0)[lb]{\smash{{\SetFigFont{12}{14.4}{\familydefault}{\mddefault}{\updefault}{\color[rgb]{0,0,0}$z'$}%
}}}}
\put(1441,-2266){\makebox(0,0)[lb]{\smash{{\SetFigFont{12}{14.4}{\familydefault}{\mddefault}{\updefault}{\color[rgb]{0,0,0}$z$}%
}}}}
\end{picture}%
\caption{The generic case is pictured.  It is also possible for $z'$
to lie on one of the dashed segments, in which case it is even closer
to the other two sides.}
\label{f:deltaprime}
\end{center}
\end{figure}
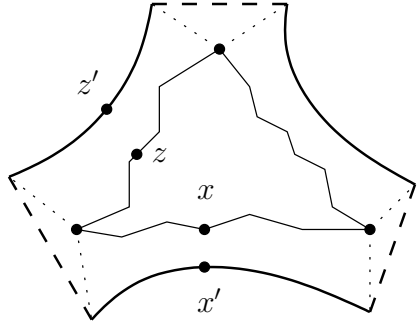
Since $H$ is $6\delta$-slim, there is some point $z'$ on another side
of $H$ so that $d(z',x')\leq 6\delta$.  If $z'$ is on
$\geod(c,d)[f(t_3),f(t_4)]$ or $\geod(e,f)[f(t_5),f(t_6)]$, then
there is a point $z$ on $\pp_{cd}[t_3,t_4]$ or $\pp_{ef}[t_5,t_6]$
with $d(z,z')\leq \Hdist$.  Otherwise, $z'$ is at most $\Displace + \delta$
from one of the vertices of our original triangle.  In either case $x$
is within $2(\Hdist)+6\delta=\Deltaprime$ 
of one of the other two sides of the triangle.
\end{proof}

\begin{torsionremark} \label{tr:bonus}
Any of the `possible' preferred paths referred to in Remark \ref{r:tf3}
will satisfy all of the above properties, and an identical proof
suffices to establish this.
\end{torsionremark}

\subsection{Preferred triangles} \label{s:preftri}

\begin{definition}\label{d:preftri}
Let $\Delta$ be a $2$-simplex.
A \emph{preferred triangle} is a continuous map
\[\psi\co \partial \Delta \to X \cup \Hbound\]
so that if $e$ is one of the three sides of $\Delta$,
then $\psi|_{e}$ is an embedding, whose image is a preferred path.
\end{definition}

\begin{remark}
If $\psi\co \partial \Delta \to X \cup \Hbound$ is a preferred
triangle, then $\psi|_{\Delta^{(0)}}$ is injective.
\end{remark}

\begin{assumption}
For ease of exposition, we will always assume that if $\psi\co
\partial \Delta \to X \cup \Hbound$ is a preferred triangle, then
$\psi$ does not send the vertices of $\Delta$ to $D^{-1}(L_2)$.  In
our applications, this will always be the case.  The following
definitions can be extended to the case that vertices are sent by
$\psi$ to $D^{-1}(L_2)$, but at the expense of making the statements a
bit more cumbersome.
\end{assumption}
For the rest of this section, we fix a preferred triangle
\[\psi\co \partial \Delta \to X \cup \Hbound\]
so that $\psi(v)\neq L_2$ for $v\in\Delta^{(0)}$.  We refer
to the elements of $\Delta^{(0)}$ as {\em corners}.

We next define a \emph{skeletal filling} of a preferred triangle.
This filling will take the form of a $1$-complex (with $3$ different
kinds of edges) inscribed on the
$2$-simplex $\Delta$.

\begin{definition}
The points in
$(D\circ\psi)^{-1}(L_2)$ will be referred to as {\em $L_2$-vertices}.
\end{definition}

We now describe the different types of preimages of
$L_2$-horoballs which can occur in a preferred triangle.

\begin{figure}[htbp]
\begin{center}
\begin{picture}(0,0)%
\includegraphics{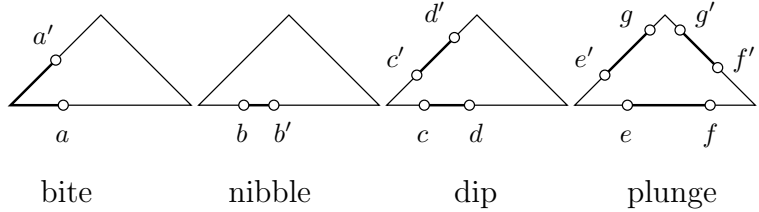}%
\end{picture}%
\setlength{\unitlength}{4144sp}%
\begingroup\makeatletter\ifx\SetFigFont\undefined%
\gdef\SetFigFont#1#2#3#4#5{%
  \reset@font\fontsize{#1}{#2pt}%
  \fontfamily{#3}\fontseries{#4}\fontshape{#5}%
  \selectfont}%
\fi\endgroup%
\begin{picture}(4489,1288)(384,-899)
\put(676,-466){\makebox(0,0)[lb]{\smash{{\SetFigFont{10}{12.0}{\familydefault}{\mddefault}{\updefault}{\color[rgb]{0,0,0}$a$}%
}}}}
\put(541,119){\makebox(0,0)[lb]{\smash{{\SetFigFont{10}{12.0}{\familydefault}{\mddefault}{\updefault}{\color[rgb]{0,0,0}$a'$}%
}}}}
\put(1981,-466){\makebox(0,0)[lb]{\smash{{\SetFigFont{10}{12.0}{\familydefault}{\mddefault}{\updefault}{\color[rgb]{0,0,0}$b'$}%
}}}}
\put(586,-826){\makebox(0,0)[lb]{\smash{{\SetFigFont{12}{14.4}{\familydefault}{\mddefault}{\updefault}{\color[rgb]{0,0,0}bite}%
}}}}
\put(1711,-826){\makebox(0,0)[lb]{\smash{{\SetFigFont{12}{14.4}{\familydefault}{\mddefault}{\updefault}{\color[rgb]{0,0,0}nibble}%
}}}}
\put(3061,-826){\makebox(0,0)[lb]{\smash{{\SetFigFont{12}{14.4}{\familydefault}{\mddefault}{\updefault}{\color[rgb]{0,0,0}dip}%
}}}}
\put(4096,-826){\makebox(0,0)[lb]{\smash{{\SetFigFont{12}{14.4}{\familydefault}{\mddefault}{\updefault}{\color[rgb]{0,0,0}plunge}%
}}}}
\put(2881,254){\makebox(0,0)[lb]{\smash{{\SetFigFont{10}{12.0}{\familydefault}{\mddefault}{\updefault}{\color[rgb]{0,0,0}$d'$}%
}}}}
\put(2656,-16){\makebox(0,0)[lb]{\smash{{\SetFigFont{10}{12.0}{\familydefault}{\mddefault}{\updefault}{\color[rgb]{0,0,0}$c'$}%
}}}}
\put(4051,254){\makebox(0,0)[lb]{\smash{{\SetFigFont{10}{12.0}{\familydefault}{\mddefault}{\updefault}{\color[rgb]{0,0,0}$g$}%
}}}}
\put(4726,-16){\makebox(0,0)[lb]{\smash{{\SetFigFont{10}{12.0}{\familydefault}{\mddefault}{\updefault}{\color[rgb]{0,0,0}$f'$}%
}}}}
\put(3151,-466){\makebox(0,0)[lb]{\smash{{\SetFigFont{10}{12.0}{\familydefault}{\mddefault}{\updefault}{\color[rgb]{0,0,0}$d$}%
}}}}
\put(4546,-466){\makebox(0,0)[lb]{\smash{{\SetFigFont{10}{12.0}{\familydefault}{\mddefault}{\updefault}{\color[rgb]{0,0,0}$f$}%
}}}}
\put(3781,-16){\makebox(0,0)[lb]{\smash{{\SetFigFont{10}{12.0}{\familydefault}{\mddefault}{\updefault}{\color[rgb]{0,0,0}$e'$}%
}}}}
\put(2836,-466){\makebox(0,0)[lb]{\smash{{\SetFigFont{10}{12.0}{\familydefault}{\mddefault}{\updefault}{\color[rgb]{0,0,0}$c$}%
}}}}
\put(1756,-466){\makebox(0,0)[lb]{\smash{{\SetFigFont{10}{12.0}{\familydefault}{\mddefault}{\updefault}{\color[rgb]{0,0,0}$b$}%
}}}}
\put(4501,254){\makebox(0,0)[lb]{\smash{{\SetFigFont{10}{12.0}{\familydefault}{\mddefault}{\updefault}{\color[rgb]{0,0,0}$g'$}%
}}}}
\put(4051,-466){\makebox(0,0)[lb]{\smash{{\SetFigFont{10}{12.0}{\familydefault}{\mddefault}{\updefault}{\color[rgb]{0,0,0}$e$}%
}}}}
\end{picture}%
\caption{Possibilities for the preimage of a single $L_2$-horoball are
shown in bold. For any letter $x$, the vertices $x$ and $x'$ form a
{\em pair} (Definition \ref{d:pairs}).}
\label{f:bitenibbledipplunge}
\end{center}
\end{figure}

\begin{definition}\label{d:nibblesetc}
Suppose that $P\in\mc{H}$ 
is such that 
$P^{L_2} \cap \psi(\partial\Delta) \neq \emptyset$.
Recall that the
unique point in $\partial X$ which is a limit of points in $P$ is
called $e_P$.
\begin{enumerate}
\item If $\psi^{-1}(P^{L_2}\cup e_P)$ contains a single corner of $\partial\Delta$,
then we say that $\psi^{-1}(P^{L_2}\cup e_P)$ is a {\em \bite}.
\item
If $\psi^{-1}(P^{L_2}\cup e_P)$ intersects exactly one side of the
triangle $\partial \Delta$, then $\psi^{-1}(P^{L_2}\cup e_P)$ is a
{\em nibble}.
\item
If $\psi^{-1}(P^{L_2}\cup e_P)$ intersects exactly two sides of the
triangle $\partial \Delta$ and is not a \bite, then
$\psi^{-1}(P^{L_2}\cup e_P)$ is a {\em \dip}.
\item
If $\psi^{-1}(P^{L_2}\cup e_P)$ intersects all three sides of the
triangle $\partial \Delta$, then it is a {\em \plunge}.
\end{enumerate}
See Figure \ref{f:bitenibbledipplunge}.
\end{definition}

\begin{lemma}
Any $L_2$-vertex is in the boundary of a unique \corner, \dip,
\nibble, or \plunge.
\end{lemma}

\begin{lemma}\label{l:pairs}
There is a partition of the $L_2$-vertices into sets, $\{S_i\}$, of cardinality
$2$ , so that if $S_i=\{x,y\}$, then
$d(\psi(x),\psi(y))\leq 1$.
\end{lemma}
\begin{proof}
Suppose that $x$ is an $L_2$-vertex.  Then $x$ is in the boundary
of a unique \corner, \dip, \nibble, or \plunge.  We consider each
of there in turn, in order to define the set $S_i$ of which $x$ is a 
member.  Let $P^{L_2}$ be the $L_2$-horoball containing
$\phi(x)$.  We delay the proof that if $S_i = \{ x,y \}$ then
$d(\psi(x),\psi(y)) \le 1$ until later in the proof.

Suppose that $x$ is contained in the boundary of a \bite.  Then
there is a unique $L_2$-vertex $y \neq x$ in $\partial\Delta$ so that
$\psi(y) \in P^{L_2}$.  Define $S_i = \{ x,y \}$ in this case.

Suppose that $x$ is contained in the boundary of a \dip.  Then
the set of $L_2$-vertices on the boundary of this \dip\ has 
cardinality $4$. We partition this into two sets by pairing the
$L_2$-vertices on different sides which are closest to the common
vertex of $\partial \Delta$.

Suppose that  $x$ is contained in a \nibble.  Then the boundary
of the \nibble\ has cardinality $2$, and we define the set containing
$x$ to be this boundary.

Finally, suppose that $x$ is contained in a \plunge.  The boundary of
the plunge has cardinality at most $6$; we partition the boundary 
into sets of
size $2$ by pairing $L_2$-vertices which are closest to a given vertex
of $\partial \Delta$.

The procedure described above (and also in Figure
\ref{f:bitenibbledipplunge}) defines a partition of the $L_2$-vertices
into sets $\{ S_i \}$ of cardinality $2$.  It remains to investigate
the distance between $\psi(x)$ and $\psi(y)$ where $\{ x,y \} = S_i$.

We consider first the case that $S_i = \{ x,y \}$, and that $x$ and $y$
lie in different sides of $\partial \Delta$.  This covers the cases
when $x$ and $y$ are contained in the boundary of a \corner, a \dip,
or a \plunge.

Now, $\psi(x) \in D^{-1}(L_2) \cap \pp_{a,b}$, say.  Suppose that $\psi(y)
\in D^{-1}(L_2) \cap \pp_{a,c}$ (we can arrange both of these by 
relabeling $a,b$ and $c$ if necessary).  Preferred paths are chosen
to be vertical from depth $L_1$ to at least $L_2+1$ (see Definition
\ref{d:sigma}).  In fact, because preferred paths travel as quickly
as possible towards $L_1$-horoballs, they must be vertical from
depth $2\delta$ to at least $L_2+1$.

Let $z_x$ be the point on $\pp_{a,b}$ which is at depth
$L_1$ and lies closest to $\psi(x)$, and let $z_y$ be defined similarly
in relation to $\psi(y)$.

Preferred paths are $\frac{1}{4}L_1$-slim by Proposition 
\ref{p:deltaprime}.  Thus, there is a point $w \in \pp_{a,c} \cup \pp_{b,c}$
which is within $\frac{1}{4}L_1$ of $z_x$.  

Suppose first that $w \in \pp_{a,c}$.  Then because preferred paths
are vertical between depths $2\delta$ and $L_2+1$, it is not difficult
to see that $d(z_x,z_y) \le \frac{1}{2}L_1$.  Because of the geometry
of combinatorial horoballs, distance between points on vertical paths
strictly decreases with increasing depth.  Thus, by depth 
$1\frac{1}{2}L_1$ in $P$, the paths $\pp_{a,b}$ and $\pp_{a,c}$ must be
within $1$ of each other.  It is now easy to see that  $d(\psi(x), \psi(y))
\le 1$.

Suppose then that $w \in \pp_{b,c}$.  Then, arguing as above,
there is a point $w_x \in D^{-1}(L_1) \cap P \cap \pp_{b,c}$ so that
$d(z_x,w_x) \le \frac{1}{4}L_1$.

Similarly, there is a point $w' \in \pp_{a,b} \cup \pp_{b,c}$ within
$\frac{1}{4}L_1$ of $z_y$.  If $w' \in \pp_{a,b}$ then we see that
$d(\psi(x),\psi(y)) \le 1$, just as above.  Therefore, suppose that
$w' \in \pp_{b,c}$, in which case we have some $w_y \in
D^{-1}(L_1) \cap P \cap \pp_{b,c}$ so that 
$d(z_y,w_y) \le \frac{1}{4}L_1$.

Note that the part of $\pp_{b,c}$ between $w_x$ and $w_y$ is
$\sigma(w_x,w_y)$.
We consider the maximum depth $\sigma(w_x,w_y)$.  First note 
that because of the
choice of $x$ and $y$, it is must be that $\pp_{b,c} \cap P^{L_2}
= \emptyset$.   Therefore, the maximum depth of $\sigma(w_x,w_y)$
is less than $L_2$.
 If this depth is at most $L_2 - 2$, then
at depth $L_2 - 2$ the preferred paths $\pp_{a,b}$ and $\pp_{a,c}$ are
at most $3$ apart, which implies that at depth $L_2$ they are at most
$1$ apart, as required.  Therefore, we are left with the case that the
maximum depth of $\sigma(w_x,w_y)$ is exactly $L_2 - 1$.  In this
case, at depth $L_2-2$, the preferred paths $\pp_{a,b}$ and $\pp_{b,c}$
are distance at most $1$ apart, as are the preferred paths 
$\pp_{a,c}$ and $\pp_{b,c}$.  Because $\sigma(w_x,w_y)$ has
maximum depth $L_2 -1$, at depth $L_2-2$ the two vertices in
$\pp_{b,c}$ are distance at most $2$ apart.  This implies that the 
distance (at depth $L_2-2$) between the vertices in 
$\pp_{a,b}$ and $\pp_{a,c}$ at depth $L_2-2$ is at most $4$.  This
implies that $d(\psi(x),\psi(y)) \le 1$.

We are now left with the case that $x$ and $y$ are contained in
the same side of $\partial\Delta$, so in the boundary of a \nibble.
Suppose that $\psi(x), \psi(y) \in \pp_{a,b}$.  Define the points
$z_x,z_y \in D^{-1}(L_1) \cap P \cap \pp_{a,b}$.  There are points
$w_x,w_y \in D^{-1}(L_1) \cap P \cap (\pp_{b,c} \cup \pp_{a,c})$ so
that $d(z_x,w_x), d(z_y,w_y) \le \frac{1}{2}L_1$.  Here there are
two cases to consider, depending on whether $w_x$ and $w_y$ are
contained in the same preferred path, or not.  In any case,
The preferred paths $\pp_{a,c}$ and $\pp_{b,c}$ do not intersect
$P^{L_2}$ (because we are considering a \nibble).  However, by
Proposition \ref{p:deltaprime}, one of $\pp_{a,c}$ or $\pp_{b,c}$ must
intersect $P^{L_2- \frac{1}{4}L_1}$ nontrivially.  

There are a number of cases to consider, depending on the maximum
depth in $P$ of $\pp_{a,c}$ and $\pp_{b,c}$.  We consider the
most complicated, and leave the remainder as exercises for the 
reader.  

Suppose that $w_x \in \pp_{a,c}$, $w_y \in \pp_{b,c}$ and that
both $\pp_{a,c}$ and $\pp_{b,c}$ intersect $P^{L_2-1}$ nontrivially.
Consider the points in 
$\pp_{a,c} \cap P \cap D^{-1}(L_1)$.  There are two of these points,
$w_x$, and $u_x$, say.  Similarly, let 
$\pp_{b,c} \cap P \cap D^{-1}(L_1) = \{ w_y, u_y \}$.

It is not difficult to see (using $\frac{1}{4}L_1$-slim triangles and the
properties of the paths $\sigma(r,s)$) that 
$d(u_x,u_y) \le \frac{1}{2}L_1$.

Therefore, at depth $L_2-3$, we have points $z_x', z_y' \in \pp_{a,b}$,
$w_x', u_x' \in \pp_{a,c}$ and $w_y', u_y' \in \pp_{b,c}$ so that
a primed point is directly beneath its unprimed counterpart.
Now, $d(z_x',w_x'), d(u_x',u_y'), d(w_y',z_y') \le 1$. Also,
the $(L_2-3)$-distance between $w_x'$ and $u_x'$ is at most
$4$ (because the maximum depth of $\pp_{a,c}$ in $P$ is $L_2-1$).
Similarly, the $L_2-3$ distance between $w_y'$ and $u_y'$ is at most
$4$.  This shows that the $L_2-3$ distance between $z_x'$ and $z_y'$
is at most $7$.  This implies (because $7 < 2^3$) that
$d(\psi(x),\psi(y)) \le 1$, as required.

This completes the proof of Lemma \ref{l:pairs}.
\end{proof}

\begin{definition}\label{d:pairs}
An element of the partition in Lemma \ref{l:pairs} will be called a
{\em pair}.
\end{definition}

\begin{definition} \label{d:skel} 
We define a $1$-complex $\skel(\psi)$ which we call the {\em skeletal
filling} of $\psi$.
The vertex set of $\skel(\psi)$ is equal to 
$\Delta^0\cup (D\circ\psi)^{-1}(L_2)$.  There are three kinds of
edges:
\begin{enumerate}
\item  The first kind come from the subdivision of $\partial \Delta$
by the vertex set.  
\item  If $\{x,y\}$ are a pair of $L_2$-vertices, in the sense of
Definition \ref{d:pairs}, and $\psi(x)=\psi(y)$, then we connect $x$
and $y$ by an edge which we call a {\em ligament}.  
\item  If $\{x,y\}$ are a pair of $L_2$-vertices so that
$\psi(x)\neq\psi(y)$, then we connect $x$ and $y$ by an edge which we
call a {\em rib}.
\end{enumerate}
\end{definition}  

\begin{lemma}\label{l:inscribe}
The identity map on $\partial \Delta$ extends to an embedding of
$\skel(\psi)$ into $\Delta$.
\end{lemma}

\begin{remark}\label{r:extend}
By Lemma \ref{l:pairs}, the map $\psi$ also extends to a map $\ddot{\psi} \co \skel(\psi) \to
X \cup \Hbound$ which sends each
ligament to a point and each rib to a single horizontal edge.
\end{remark}

\begin{definition}
Suppose there are (possibly degenerate)
subintervals $\sigma_1$ and $\sigma_2$ of sides of
$\partial \Delta$ so that $\psi(\sigma_1)=\psi(\sigma_2)$, and suppose
that these are chosen to be maximal such intervals with endpoints in
the vertex set of $\skel(\psi)$.  Then the minimal subcomplex of
$\skel(\psi)$ containing all edges with both endpoints in
$\sigma_1\cup\sigma_2$ is called a {\em leg}.  A leg contains one or
two {\em distinguished} ligaments, joining the endpoints of $\sigma_1$ to
the endpoints of $\sigma_2$.
\end{definition}

We note that $\psi$ collapses each leg to a subsegment of a preferred
path.

\begin{lemma}\label{l:ThreeLegs}
Each leg of $\skel(\psi)$ contains exactly on corner of $\Delta$.
In particular, $\skel(\psi)$ has at most $3$ legs.
\end{lemma}
\begin{proof}
Consider a corner $v$ of $\partial\Delta$, and the two sides $e_1$
and $e_2$ containing $v$.  Let $x_1 \in e_1$ and $x_2 \in e_2$ be
the endpoints of the ligament joining $e_1$ and $e_2$ 
 which is furthest
from $v$.

Let $\psi(e_1) = \pp_{a,b}$ and $\psi(e_2) = \pp_{a,c}$, so
$\psi(v) = a$.  Suppose that $\psi(x_1) = \psi(x_2) \in D^{-1}(L_2)
\cap P$.  Then $P \in \mc{H}_{a,b} \cap \mc{H}_{a,c}$. 
By Axiom (A6) and the definition of preferred paths
those parts of $\pp_{a,b}$ and $\pp_{a,c}$
between $a$ and $\psi(x_1)=\psi(x_2)$ are identical.
This implies that the segments between $v$ and $x_1$ and
$v$ and $x_2$, and the ligaments joining them, form a leg.

If there is no ligament joining the edges $e_1$ and $e_2$ then there
is no leg intersecting both of these sides.
\end{proof}

\begin{definition}
Let $\Gamma$ be a graph and let $e$ be an edge of $\Gamma$ which
separates $\Gamma$ into two components.  We say $\Gamma_1$ and
$\Gamma_2$ are obtained from $\Gamma$ by  {\em surgery along $e$} if
there are edges $e_1$ in $\Gamma_1$ and $e_2$ in $\Gamma_2$ which can
be identified to give a graph isomorphic to $\Gamma$.
\end{definition}

\begin{definition}
Successive surgery along the distinguished ligaments of $\skel(\psi)$
yields a collection of at most $4$ graphs, at most three of which come
from legs, and exactly one of which contains pieces of all three edges
of $\partial \Delta$.  We call the graph which contains pieces of all three
edges of $\partial\Delta$ the {\em middle} of $\skel(\psi)$. \footnote{
In an earlier version of this paper, preferred triangles could have {\em
feet}, as well as legs and a middle.  The improved construction
in Section \ref{s:horoballs} does away with the need for feet.}
 (See Figure \ref{f:triskelion}.)
\begin{figure}[htbp]
\begin{center}
\begin{picture}(0,0)%
\includegraphics{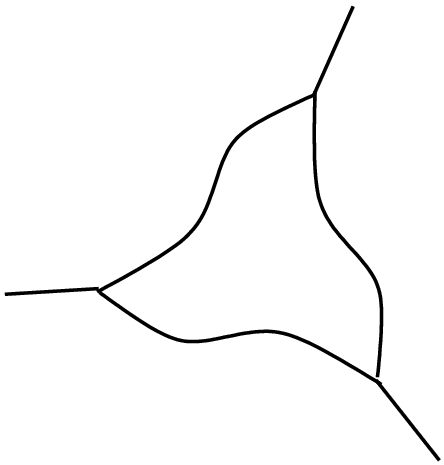}%
\end{picture}%
\setlength{\unitlength}{4144sp}%
\begingroup\makeatletter\ifx\SetFigFont\undefined%
\gdef\SetFigFont#1#2#3#4#5{%
  \reset@font\fontsize{#1}{#2pt}%
  \fontfamily{#3}\fontseries{#4}\fontshape{#5}%
  \selectfont}%
\fi\endgroup%
\begin{picture}(2030,2119)(1267,-3913)
\end{picture}%
\caption{A somewhat stylized picture of the image of a preferred triangle.}
\label{f:triskelion}
\end{center}
\end{figure}
\end{definition}

\begin{lemma} \label{l:InTwoHSys}
Let $\psi \co \partial \Delta \to X \cup \Hbound$ be a preferred triangle,
with edges $e_1,e_2,e_3$, and suppose that $\psi(e_1) = \pp_{a,b}$,
$\psi(e_1) = \pp_{a,c}$ and $\psi(e_3) = \pp_{b,c}$.  Suppose further
that  $P$ is a horoball which is such that $\psi(e_1) \cap 
P^{L_2} \ne \emptyset$.  Then $P \in \mc{H}_{a,b}$ and either
$P \in \mc{H}_{a,c}$ or $P \in \mc{H}_{b,c}$.
\end{lemma}
\begin{proof}
Follows from the fact that the image $\psi(\partial\Delta)$ is a 
quasi-geodesic triangle which is $\frac{1}{4}{L_1}$-slim and that
$L_2 - \frac{1}{4}L_1 = 3L_1 - \frac{1}{4}L_1 > L_1$.
\end{proof}

\begin{definition} \label{d:assoc}
Suppose that $\psi \co \partial\Delta \to X \cup \Hbound$ is 
a preferred triangle, with corners $v_1, v_2, v_3$ so that
$\psi(v_1) = a, \psi(v_2) = b$ and $\psi(v_3) =
c$.  Let $e_1$ be the edge of $\partial\Delta$ joining
$v_1$ to $v_2$.   If $\psi^{-1}(P^{L_2}) \neq\emptyset$ and
$P \in \mc{H}_{a,b} \cap \mc{H}_{a,c}$ then
we say that $P$ is {\em associated to $v_1$}.  We define
horoballs associated to $v_2$ and $v_3$ analogously.

If $P$ is a horoball so that $\psi^{-1}(P^{L_2})$ is nonempty and
contained in a leg, $l$, then we say that $P$ is {\em associated
to $l$}.
\end{definition}

\begin{remark}
 By Lemma \ref{l:InTwoHSys}, if $P^{L_2} \cap
\psi(e_1) \neq 0$, then $P$ is contained in two of
$\mc{H}_{a,b}, \mc{H}_{a,c}$ and $\mc{H}_{b,c}$.

Therefore, for any horoball $P$, if $\psi^{-1}(P^{L_2})$ is nonempty
then $P$ is associated to some corner of $\Delta$.
\end{remark}

It could be that a horoball $P$ as in Definition \ref{d:assoc}
is contained in all three of $\mc{H}_{a,b}, \mc{H}_{a,c}$ and 
$\mc{H}_{b,c}$, in which case $P$ is associated to all three corners
of $\partial \Delta$.  However

\begin{lemma} \label{l:oneplunge}
There is at most one horoball associated to all three corners of 
$\partial\Delta$.
\end{lemma}

\begin{proposition} \label{p:sixhoroballs}
Let $\psi \co \partial\Delta \to X \cup \Hbound$ be a preferred
triangle and let $v$ be a vertex of $\partial\Delta$.

There is at most one horoball which is associated to $v$ and
is not also associated to a leg in $\skel(\psi)$.
\end{proposition}
\begin{proof}
Suppose that the horoballs $A$ and $B$ are associated to
$v$.  

Let $a,b,c$ be the images of the vertices of $\partial\Delta$, with
$a = \psi(v)$.  Then $A,B \in \mc{H}_{a,b} \cap \mc{H}_{a,c}$ and suppose that $A < B$ in the order coming from
$\mc{H}_{a,b}$.We will prove that
$A$ is associated to a leg in $\partial\Delta$.

The preferred paths $\pp_{a,b}$ and $\pp_{a,c}$ coincide at
least until they reach the horoball $B$
and in particular, the intersections of these paths with $A$ is
identical.  Since one of $A^{L_2} \cap \pp_{a,b}$ and $A^{L_2}
\cap \pp_{a,c}$ is nonempty (since $A$ is associated to $v$),
both of these sets are non-empty.

It is not difficult to see that $\psi^{-1}(A^{L_2})$ is contained
in a leg.
\end{proof}

\begin{corollary} \label{c:ribsbound}
Let $\psi \co \partial\Delta \to X \cup \Hbound$ be a preferred
triangle.  The total number of ribs and ligaments in the feet and
middle of $\skel(\psi)$ is at most $6$.  In particular, the number of ribs in $\skel(\psi)$ is at most $6$.
\end{corollary}
\begin{proof}
Each $L_2$-vertex in $\skel(\psi)$ maps to a horoball 
associated to some corner of $\partial\Delta$.  If $u$ is an
$L_2$-vertex in a foot or middle of $\skel(\psi)$ and $P$
is the horoball containing $\psi(u)$ then $P$ is not associated
to a leg.

For each horoball, $A$, the set of $L_2$-vertices in $\skel(\psi)$
which map to $A^{L_2}$ has cardinality at most $6$, and at most
$4$ unless $\psi^{-1}(A^{L_2} \cup e_A)$ is a plunge.  Thus
there are $3$ ribs or ligaments for a plunge, and at most $2$
for other horoballs.

Lemma \ref{l:oneplunge} implies that there is at most one plunge
and Proposition \ref{p:sixhoroballs} implies that there are at most
$3$ horoballs not associated to legs.  The result now
follows easily.
\end{proof}

\begin{corollary} \label{c:SkelBound}
Let $\psi \co \partial\Delta \to X \cup \Hbound$ be a preferred
triangle and let $M$
be the union of the feet and middles of $\skel(\psi)$. Then
\[ | M^{(0)} | \le 15. \]
\end{corollary}
\begin{proof}
There are $3$ vertices for the corners of $\partial \Delta$.

Also, there are at most $6$ ribs or ligaments in the feet and
middle of $\skel(\psi)$, and each of these contributes two vertices.
\end{proof}

We now proceed to decompose preferred triangles.  A key
feature of preferred triangles is Proposition \ref{p:geodsinthick} below.

\begin{definition} [Sub-pictures]
Suppose that $\psi \co \partial \Delta \to X \cup \Hbound$ is a preferred 
triangle.
Successive surgery along all of the ribs and ligaments of $\skel(\psi)$
yields a collection of graphs.  We call these graphs {\em sub-pictures}.
\end{definition}

The following result follows immediately from the fact that the depth
function $D$ is continuous.

\begin{lemma}
Suppose that $\bf Pic$ is a sub-picture of $\skel(\psi)$.  Then either
$\psi({\bf Pic}) \subset D^{-1}([0,L_2])$ or $\psi({\bf Pic}) \subset
D^{-1}([L_2,\infty))$.
\end{lemma}

\begin{definition}\label{d:thick}
[Thick-thin decomposition of preferred triangles]
Suppose that $\psi \co \partial\Delta\ to X \cup \Hbound$ is a preferred
triangle, and that $\bf Pic$ is a sub-picture of $\skel(\psi)$.

If $\psi({\bf Pic}) \subset D^{-1}([0,L_2])$ then we call $\bf Pic$ a
{\em thick sub-picture}. If $\psi({\bf Pic}) \subset D^{-1}([L_2,\infty))$
then we call $\bf Pic$ a {\em thin sub-picture}.
\end{definition}

\begin{proposition} \label{p:geodsinthick}
If $\psi \co \partial\Delta \to X \cup \Hbound$ is a 
preferred triangle, and $u$ and $v$ are vertices in the
same thick sub-picture of $\skel(\psi)$ then any geodesic
$\geod(\psi(u),\psi(v))$ joining $\psi(u)$ to $\psi(v)$ does not
penetrate any $(L_2 + L_1)$-horoball.
\end{proposition}
\begin{proof}
Let $\bf Pic$ be a thick sub-picture of $\skel(\psi)$.
Then $\psi({\bf Pic}) \subset D^{-1}([0,L_2])$.
Furthermore, the image of $\bf Pic$ under $\psi$
consists of subsegments of preferred paths $\pp_{a,b}, \pp_{a,c}$
and $\pp_{b,c}$, say, together with single edges from $D^{-1}(L_2)$,
corresponding to the images of ribs in $\bf Pic$.

If $\pp_{a,b}$ and $\pp_{a,b}$ both intersect $\psi({\bf Pic})$
nontrivially, then they either intersect in $\psi({\bf Pic})$ (in case
their preimages are joined by a ligament) or else are
joined by an edge of length $1$ in $D^{-1}(L_2)$ (in case their
preimages are joined by a rib).

Thus, any point in $\psi({\bf Pic})$ lies within distance at most
$\frac{1}{2}$ from $\pp_{a,b} \cap \pp_{a,c} \cap \pp_{b,c}$.

Suppose $\pp_{a,b} \cap \psi({\bf Pic}) \neq \emptyset$. Let
$m_{a,b}$ be the point in $\pp_{a,b} \cap \psi(\bf Pic)$ closest
to $a$, and let $n_{a,b}$ be the point closest to $b$. Define
points $m_{a,c}, n_{a,c}, m_{b,c}$ and $n_{b,c}$ analogously.
Not all of these points need be defined, if the intersection of some 
preferred path with $\psi({\bf Pic})$ is empty.  Also, it is certainly
possible for some of these points to be the same, if they lie
in the image of a ligament, or if they are the image of a corner
of $\partial\Delta$.

Suppose that $m_{a,b} = \pp_{a,b}(t)$.  Let $\hat{m}_{a,b}
= \geod(a,b)(f(t))$, where $f$ is the function from 
Proposition \ref{p:async}.

Use the function $f$ to define points $\hat{n}_{a,b}$,
$\hat{m}_{a,c}$, etc.

Note that $d(m_{a,b},\hat{m}_{a,b}) \le \Displace$, and 
similarly for the other points and their hatted counterparts.

The points $m_{a,b}, n_{a,b}, \ldots$ are
partitioned into some of the pairs defined earlier in this 
section, and in particular each pair consists of one or two
points which are distance at most $1$ apart.

Thus, the geodesics between the points $\hat{m}_{a,b},
\hat{n}_{a,b}, \ldots$, consisting of the subsegments of
$\geod(a,b), \geod(a,c), \geod(b,c)$, and paths joining the
points whose un-hatted counterparts correspond to a `pair'.

Therefore, depending on whether there are two or three preferred
paths which intersect $\bf Pic$ nontrivially, the points $\hat{m}_{a,b}$,
$\hat{n}_{a,b}$, etc. form either a hexagon or a quadrilateral.  If
it is a hexagon, then three of the sides have length at most 
$2(\Hdist) + 1$, and if it is a quadrilateral two of the sides of length at
most $2(\Hdist)+1$.  Call the sides $[\hat{m}_{a,b},\hat{n}_{a,b}]$ etc.,
the `long' sides of this polygon.
Any point on a long side of this polygon
is distance at most $2(\Hdist)+1 + 3\delta$ from one of the other
long sides.

Now consider any vertices $u, v \in \bf Pic$.  As already noted,
there are points $u', v'$ so that $u'$ and $v'$ are contained in
preferred paths and $d(u,u'), d(v,v') \le \frac{1}{2}$.  

By Proposition \ref{p:deltaprime} there are points $u_1, v_1$
on different preferred paths to $u'$ and $v'$, respectively
so that $d(u',u_1), d(v',v_1) \le \delta' + 1$.  Since there
are only three preferred paths, at least two of $u',u_1, v',
v_1$ must lie on the same preferred path.  Since $u'$ and 
$u_1$ lie on different preferred paths, as do $v'$ and $v_1$,
one of the two points lies within $\delta'$ of $u'$ and one within
$\delta'$ of $v'$.  Let this pair of points on the same
preferred path be $\overline{u}$ and $\overline{v}$, so that
$d(\overline{u},u'), d(\overline{v},v') \le \delta'$.  That
part of the preferred path which lies between $\overline{u}$
and $\overline{v}$ lies in a thick sub-picture, and so does not
penetrate any $(L_2+2)$-horoball.  Let $\mu$ be this part
of the preferred path.

By Proposition \ref{p:async}, there are points $\hat{u}, \hat{v}$
on a long side of the polygon described above, so that
$d(\overline{u},\hat{u}), d(\overline{u},\hat{u}) \le \Hdist$.
The geodesic between $\hat{u}$ and $\hat{v}$ (which lies
on a long side of the polygon) lies at Hausdorff distance at most
$\Hdist$ from $\mu$, and therefore does not penetrate
any $(L_2 + (\Hdist) + 2)$-horoball.

Now, consider the geodesic quadrilateral formed by
$u,v, \hat{u}, \hat{v}$.  We have 
\begin{eqnarray*}
d(u,\hat{u}) & \le & d(u,u') + d(u',\overline{u})
+ d(\overline{u},\hat{u})\\
& \le & \delta' + (\Hdist) + 1\frac{1}{2},
\end{eqnarray*}
and, similarly, $d(v,\hat{v}) \le \delta' + (\Hdist) + 1\frac{1}{2}$.

Therefore, any point on the geodesic between $u$ and $v$
lies at most $\delta' + (\Hdist) + 1\frac{1}{2} + 2\delta$ of some
point on the geodesic between $\hat{u}$ and $\hat{v}$.  

This implies that the geodesic between $u$ and $v$ does not 
penetrate any $(L_2 + 2(\Hdist) + 2\delta + 3\frac{1}{2})$-horoball.
Since $L_2 + 2(\Hdist) + 2\delta + 3\frac{1}{2} < L_2 + L_1$, the proof
is finished.
\end{proof}

\begin{torsionremark} \label{r:tf4}
In the context of the ``quantum" preferred paths, when $G$ is not
torsion-free, not all choices of triples of preferred paths will
have the controlled properties of preferred triangles in this section.
In applications, it is possible to make consistent choices.  We will
say more about this in Part 2.
\end{torsionremark}

\def\qtogeod {K + 25\delta + 9} 

\def\qConst {6000\delta^2}

\section{A homological bicombing}\label{section:bicombing}

In this section we describe how the results from previous
sections in this paper may be combined with results of
Mineyev from \cite{min:str} to construct a quasi-geodesic
homological bicombing of $X \cup \Hbound$ (in the sense of
Definitions \ref{d:homologicalextension} and \ref{d:quasigeodesic}).

\subsection{Mineyev's bicombing}

In this subsection, we briefly recall a construction of Mineyev
from \cite{min:str}.
We need a slightly more general statement than appears
in \cite{min:str}, and we explain how Mineyev's proof
implies Theorem \ref{t:Mineyev} below.

Suppose that $\Gamma$ is a locally finite graph
which is $\delta$-hyperbolic for some integer $\delta \ge 1$.

\begin{remark}
In his construction, Mineyev further assumes that $\Gamma$ 
has bounded valence.  This is important for the area bounds, but 
not for making the definitions.
\end{remark}

Suppose that the group $G$ acts freely on $\Gamma$.
Let $\geod$ be a $G$-equivariant geodesic bicombing
on $\Gamma$ (see Definition \ref{d:geodesicbicombing}).  Let
$P$ be the homological bicombing induced by $\geod$ (see Remark
\ref{r:geodtohom}).

For each vertex $a$ in $\Gamma$, define
\[ pr_a : \Gamma^{(0)} \to \Gamma^{(0)},	\]
as follows:
\begin{itemize}
\item $pr_a(a) = a$; and
\item if $b \neq a$ then $pr_a(b) = \geod(a,b)(r)$, where $r$
is the largest (integral) multiple of $10\delta$ which is strictly
less than $d(a,b)$.
\end{itemize}

For vertices $a,b$ in $\Gamma$, the {\em flower at $b$ with
respect to $a$} is the set
\[	Fl(a,b) = S(a,d(a,b)) \cap B(b,\delta) \subset \Gamma^{(0)}.	\]

Now, for each pair of vertices $a,b \in \Gamma$ define a 
$0$-chain $f(a,b)$ in $\Gamma$ inductively on
the distance $d(a,b)$ as follows:
\begin{itemize}
\item if $d(a,b) \le 10\delta$ then $f(a,b) = b$;
\item if $d(a,b) > 10\delta$ and $d(a,b)$ is not an integral
multiple of $10\delta$ then $f(a,b) = f(a,pr_a(b))$; and
\item if $d(a,b) > 10\delta$ and $d(a,b)$ is an integral multiple
of $10\delta$ then
\[	f(a,b) = \frac{1}{\# Fl(a,b)} \sum_{x \in Fl(a,b)} f(a,pr_a(x)).	\]
\end{itemize}

For each vertex $a \in \Gamma^{(0)}$, define a $0$-chain, 
$star(a)$ by
\[ star(a) = \frac{1}{\# B(a,7\delta)} \sum_{x \in B(a,7\delta)} x.	\]
By linearity, $star(a)$ is defined for any $0$-chain $a$.

Now define, for $a,b \in \Gamma^{(0)}$,
\[	\bar{f}(a,b) = star(f(a,b)).	\]

We now define a homological bicombing $Q'$ on $\Gamma$.
First note that, by linearity, $P_{a,b}$ makes sense when $a$ is any $0$-chain.
The $1$-chain $Q'_{a,b}$ is defined inductively on $d(a,b)$, as follows:
if $d(a,b) \le 10\delta$ then $Q'_{a,b} = P_{a,b}$.  Suppose now that
$d(a,b) > 10\delta$.  By \cite[Proposition 7(2)]{min:str}
\[	\text{supp}(\bar{f}(a,b)) \subseteq B(P_{a,b}(10\delta),8\delta).	\]
Note that Mineyev's proof of this does not use the bounded valence
assumption.
Therefore, for each $x \in \text{supp}(\bar{f}(b,a))$ we have
$d(a,x) < d(a,b)$, so $Q'_{a,x}$ is defined by induction.  Define
$Q'_{a,\bar{f}(b,a)}$ by linearity over the second variable and define
\[	Q'_{a,b} = Q'_{a,\bar{f}(b,a)} + P_{\bar{f}(b,a),b}.	\]
Finally, we define 
\[	Q_{a,b} = \frac{1}{2} (Q'_{a,b} - Q'_{b,a}),\]
so that $Q$ is anti-symmetric.

Mineyev proves that when $\Gamma$ has bounded valence,
the bicombing $Q$ has {\em bounded area}, in the sense
of Theorem \ref{t:Mineyev} below.  When $\Gamma$ does
not have bounded valence, these area bounds break down 
completely.  However, we are only going to use the 
bicombing from \cite{min:str} on a subset of $X$ of uniformly 
bounded depth, and on such a subspace the valence {\em is}
uniformly bounded, and Mineyev's techniques apply.

The proof of the area bound for a given triangle occurs
entirely in the $60\delta$-neighborhood of the union
of the three geodesic sides.  Moreover, if we have a 
bound on the valence in some part of a graph, then we can 
calculate a bound on the number of vertices in any given
ball which lies entirely within the chosen part.
The proof from \cite{min:str} now applies directly to prove the
following theorem.

\begin{theorem}(Mineyev) \label{t:Mineyev}
There is a function $T=T(\delta,v)$ so that:  For any finite
valence
$\delta$-hyperbolic graph $\Gamma$, and group $G$ acting
freely on $\Gamma$,
there is a ($18\delta$)-quasi-geodesic 
$G$-equivariant antisymmetric homological
bicombing $Q$ so that 
\[|Q_{ab}+Q_{bc}+Q_{ca}|_1\leq T(\delta,v)\]
whenever $a$, $b$, $c$ are vertices of $\Gamma$ spanning a geodesic
triangle 
$\mathcal{T}_{abc}$ so that every vertex in the
$60\delta$-neighborhood of $\mathcal{T}_{abc}$ has valence less than
$v$. 
\end{theorem}
Theorem \ref{t:Mineyev} is a key ingredient in the construction 
in this section.  The other key ingredients are the construction of 
preferred paths and the analysis of preferred triangles from
Section \ref{section:pp}

\subsection{The bicombing $q$}

In the remainder of  this section we define our homological bicombing 
$q$ of 
$(X \cup \Hbound) \times (X \cup \Hbound)$, which
uses preferred paths and Mineyev's bicombing $Q$
from Theorem \ref{t:Mineyev}.
(See Definition \ref{d:homologicalextension} for the definition
of homological bicombing which allows some points to be ideal.)

Suppose that $a,b \in X \cup \Hbound$.  Let
$\pp_{a,b}$ be the preferred path between $a$ 
and $b$.

Decompose $\pp_{a,b}$ into subintervals, oriented
consistently with $\pp_{a,b}$, which
lie either entirely within $D^{-1}([0,L_2])$ or entirely
within $D^{-1}([L_2,\infty))$, and so that the endpoints
of these subintervals lie in $\{ a,b \} \cup D^{-1}(L_2)$.

By the way that preferred paths were defined, there is
a unique way of performing this decomposition.

Suppose that $\mu$ is a subinterval in the decomposition
of $\pp_{a,b}$ so that $\mu$ lies in $D^{-1}([0,L_2])$, and
let $x$ and $y$ be the endpoints of $\mu$.

\begin{lemma} \label{l:NotTooDeep1}
The geodesic between $x$ and $y$ does not intersect
any $(L_2 + 2(\Hdist) + 2\delta)$-horoball.
\end{lemma}
\begin{proof}
By Proposition \ref{p:async}, there exist $x', y' \in \geod(a,b)$
so that $d(x,x'), d(y,y') \le \Hdist$, and so that the Hausdorff
distance between $\pp_{a,b}[x,y]$ and $\geod(a,b)[x',y']$ is at
most $\Hdist$. Thus, $\geod(a,b)[x',y']$ does not penetrate
any $(L_2 + (\Hdist))$-horoball.

The geodesic quadrilateral with vertices $x,x',y,y'$ has two sides
of length at most $\Hdist$.  Therefore, the geodesic between
$x$ and $y$ lies in a $(\Hdist) + 2\delta$ neighborhood of 
$\geod(a,b)[x',y']$, and hence does not penetrate any
$(L_2 + 2(\Hdist) + 2\delta)$-horoball, as required.
\end{proof}

\begin{definition}\label{d:q1}{\em (Definition of $q$)}
Suppose that $a,b \in X \cup \Hbound$ are distinct.
Let $\pp_{a,b}$ be the preferred path between $a$ and $b$,
and let $\underline{\pp_{a,b}}$ be the induced $1$-chain.
The decomposition of $\pp_{a,b}$ described above induces
a decomposition of $\underline{\pp_{a,b}}$.
Let $\mu$ be an element of the decomposition of
$\underline{\pp_{a,b}}$ for
which $\text{supp}(\mu) \subset D^{-1}[0,L_2]$.
Let $\partial\mu = \mu_+ - \mu_-$.

Taking the sum over all such $\mu$ we define
\[	q_{a,b} = \underline{\pp_{a,b}} - \sum_\mu \big(
Q(\mu_-, \mu_+) - \mu \big)	.	\]
\end{definition}

Because $\partial\mu = \partial Q(\mu_-,\mu_+) =
\mu_+ - \mu_-$, we have
$\partial\underline{\pp_{a,b}} = \partial q_{a,b}$.  We claim that
$q$ is a homological bicombing on $X \cup \Hbound$
in the sense of Definition \ref{d:homologicalextension}
and that furthermore it is $\qConst$-quasi-geodesic 
in the sense
of Definition \ref{d:quasigeodesic}.  It also has nice properties
analogous to those in Theorem \ref{t:Mineyev} above.
See Theorem \ref{t:Xfillings} below for the precise
statements about the bicombing $q$.

\begin{proposition} \label{p:qneargeod}
For any $a, b \in X \cup \Hbound$,
$\text{supp}(q_{a,b})$ lies in a $\big( \qtogeod \big)$-neighborhood of any
geodesic between $a$ and $b$.
\end{proposition}
\begin{proof}
By Corollary \ref{c:ppneargeodesic}, $\pp_{a,b}$ lies in 
a $(K+7\delta + 9)$-neighborhood of $\geod(a,b)$.

By construction and Theorem \ref{t:Mineyev},
$\text{supp}(q_{a,b})$ lies in an $18\delta$-neighborhood
of $\pp_{a,b}$.
\end{proof}

Since $K + 25\delta + 9 < \qConst$, Proposition 
\ref{p:qneargeod} proves the first of the
two required statements for $q$ to be $\qConst$-quasi-geodesic.
We now prove the remaining requirement.

\begin{proposition} \label{p:qisqg}
$q$ is a $\qConst$-quasi-geodesic homological
bicombing on $X \cup \Hbound$.
\end{proposition}
\begin{proof}
By Proposition \ref{p:qneargeod}, it remains to prove statement (2) of Definition \ref{d:quasigeodesic}.
Let $a,b \in X^{(0)}$ be distinct.  By Theorem \ref{t:Mineyev}
and the definition of $q$, we have 
\[	|q(a,b)|_1 \le 18\delta |\pp_{a,b}|_1.	\]
By Corollary \ref{c:quasigeodesic}, the length of $\pp_{a,b}$
is at most $2d(a,b) + \Pepsilon < 325\delta d(a,b)$.
\end{proof}
\begin{remark} 
It is also possible to prove a suitable refinement of 
statement (2) as alluded to in Remark \ref{r:wrongdef}.
\end{remark}

\subsection{Bounded thick area}

The main result of Section \ref{section:bicombing}
is Theorem \ref{t:Xfillings} below.

\begin{definition} \label{d:c_abc}
Suppose $a,b,c \in X^{(0)} \cup \Hbound$.  We define
a $1$-cycle $c_{abc}$ as follows:

Let $\phi \co \Delta \to X \cup \Hbound$ be the preferred
triangle associated to the triple $(a,b,c)$.
Associated to $\phi$ is the graph $\skel(\phi)$
(see Definition \ref{d:skel}), which has associated
thick sub-pictures (see Definition \ref{d:thick}).
Let ${\bf Pic}$ be a thick sub-picture of 
$\skel(\phi)$.  The vertices of ${\bf Pic}$
inherit a circular order $(v_1, \ldots , v_n)$ from the order $(a,b,c)$.

Define\[	c_{\bf Pic} = \sum_{i \in \Z/m} Q(\ddot{\phi}(v_i),\ddot{\phi}(v_{i+1})).	\]
Finally, define
\[	c_{abc} = \sum_{\bf Pic} c_{\bf Pic},	\]
where the sum is over all thick sub-pictures of $\skel(\phi)$.
\end{definition}

\begin{observation}
If ${\bf Pic}$ is a thick sub-picture lying in a leg of $\skel(\phi)$
then $c_{\bf Pic} = 0$.
\end{observation}

The following is a key theorem for our proof of Theorem \ref{t:Zhlii},
one of the major steps in proving Theorem \ref{t:rhds}:
\begin{theorem}\label{t:Xfillings}
There exists a constant $T_1$, depending only on $X$, so that 
for all $a,b,c \in X \cup \Hbound$ there is a $1$-cycle $c_{abc}$,
as described in Definition \ref{d:c_abc} above then
\[	|c_{abc}|_1 \le T_1.	\]
Also,
\[	\text{supp} \big( q(a,b) + q(b,c) + q(c,a) - c_{abc} \big)
\subset D^{-1}[L_2,\infty).	\]
\end{theorem}
\begin{proof}
We have already noted that the thick sub-pictures which lie inside
a leg do not contribute anything to $c_{abc}$.
Thus we may concentrate on thick sub-pictures lying in the feet
or middle of $\skel(\phi)$.  By Corollary \ref{c:SkelBound}, there are
at most $15$ vertices in total in all such sub-pictures.  Therefore,
we can triangulate all of these sub-pictures with at most $13$ triangles,
whose vertices all appear as a vertex in one of the thick sub-pictures.

Let $u$ and $v$ be such vertices.  By Proposition \ref{p:geodsinthick},
the geodesic between $u$ and $v$ does not intersect any 
$(L_2 + L_1)$-horoball.

Let $v$ be the maximum valence of any vertex in 
$D^{-1}[0,L_2 + L_1 + 18\delta]$, and let $T(\delta,v)$
be as in Theorem \ref{t:Mineyev}.  Define $T_1 = 13T(\delta,v)$.

We can express $c_{abc}$ as the sum of at most $13$ $1$-chains
of the form $Q_{u_1,u_2} +Q_{u_2,u_3} + Q_{u_3,u_1}$, where
$u_1, u_2$ and $u_3$ are the images under $\ddot{\phi}$ of 
vertices in thick sub-pictures of $\skel(\phi)$ which lie in the feet
or middle.

The result now follows from Theorem \ref{t:Mineyev}.
\end{proof}

The following is immediate from Theorem \ref{t:Xfillings} and
Theorem \ref{t:tfae}.
\begin{corollary} \label{c:Xfillings}
For all $a,b,c \in X \cup \Hbound$ there exists a $2$-chain $\omega_{abc}$ so that
\begin{enumerate}
\item $\partial \omega_{abc} = c_{abc}$; and
\item $|\omega_{abc}|_1 \le M_XT_1$,
\end{enumerate}
where $M_X$ is the constant for the linear homological isoperimetric
function for $X$.
\end{corollary}

\begin{torsionremark} \label{r:tf5}
In the presence of torsion, the bicombing $q$, and the
$2$-chains $\omega_{abc}$, can be defined
without much difficulty using ``averaged" preferred paths, and the
ideas already contained in this section.
\end{torsionremark}

\newpage
\part{Dehn filling in relatively hyperbolic groups}

\section{Dehn filling in groups}\label{s:dehnintro}
In Part 2 of this paper, we provide an application of
the constructions from Part 1.  

\begin{definition}\label{d:slopelength}
Let $G$ be a group, and $P$ a subgroup.  Suppose that
$G$ is generated by $S$ and $P$ is generated by $P \cap S$.
If $K \unlhd P$ is a normal subgroup of $P$ then we define
\[ |K|_P = \inf_{k \in K \smallsetminus \{ 1\} } |k|_{P\cap S} ,\]
where $|k|_{P \cap S}$ is the distance from $k$ to the identity
in the Cayley graph $\Gamma (P , P \cap S)$.
By convention $| \{ 1 \} |_P = \infty$.
\end{definition}

 In the special case that $P$ is
free abelian and $K$ the cyclic group generated by $\kappa\in K$, 
$|K|_P$ is just the length of $\kappa$ in $P$.

The main result of Part 2 is the following theorem.

\begin{theorem} \label{t:rhds}
Let $G$ be a torsion-free group, which is hyperbolic
relative to a collection $\mc{P} = \{ P_1, \ldots , P_n \}$
of finitely generated subgroups.  Suppose that $S$
is a generating set for $G$ so that for each $1 \le i \le n$ 
we have $P_i = \langle P_i \cap S \rangle$.

There exists a constant $B$ depending only on 
$(G, \mc{P})$ so that for any collection $\{ K_i \}_{i=1}^n$ 
of subgroups satisfying
\begin{itemize}
\item $K_i \unlhd P_i$; and
\item $| K_i |_{P_i} \ge B$,
\end{itemize}
then the following hold, where $K$ is the normal closure in $G$ of
$K_1\cup \cdots \cup K_n$. 
\begin{enumerate}
\item The map $P_i/K_i\xrightarrow{\iota_i}G/K$ given by $pK_i\mapsto pK$
  is injective for each $i$.
\item $G/K$ is hyperbolic relative to the collection
  $\mc{Q}=\{\iota_i(P_i/K_i)\mid 1\leq i\leq n\}$.
\end{enumerate}
\end{theorem}

In fact, much more than this is true.  For example, for $i \neq j$ we have
\[ \iota_i(P_i/K_i) \cap \iota_j(P_j/K_j) = \{ 1 \},	\]
(see Corollary \ref{c:ijdisjoint}) and each $\iota_i(P_i/K_i)$ is malnormal in $G/K$ (see Corollary
\ref{c:conjugatesdisjoint}).  Also, if $G$ is non-elem\-en\-ta\-ri\-ly
hyperbolic relative to $\{ P_i \}$, then $G/K$ is non-elem\-en\-ta\-ri\-ly
hyperbolic relative to $\{ \iota_i(P_i/K_i) \}$ (see Theorem 
\ref{t:nonelementary}).

The remainder of the paper is devoted to the proof
of Theorem \ref{t:rhds}, and the subsidiary assertions mentioned
above.  

Theorem \ref{t:rhds}
clearly holds if $n=1$ and $G = P_1$, so we henceforth assume
(without mention) that this is not the case.

\section{Equations involving parabolics and skeletal
fillings of surfaces.}\label{s:dyingwords} 
In this section we suppose that $G$ is hyperbolic relative to
$\mc{P}=\{P_1,\ldots,P_n\}$, and that $X=X(G,\mc{P},S,\mc{R})$ is the
cusped space associated to some compatible set of generators $S$, and
some collection of relators $\mc{R}$, as described in Section
\ref{section:cusped}.  Finally $H\unlhd G$ is an arbitrary normal
subgroup of $G$.

In order to use the geometry of $X$ to study the quotients of $G$, we
will
need to turn equations in $G/H$ or (represented as maps of compact
planar surfaces into $X/H$) into ``pleated surfaces'' in $X/H$.
Exactly what this means depends on the context, and will become clear
as we proceed.

\subsection{Lifting and straightening}\label{s:extensions} 
Let $\Gamma_H = \Gamma(G,S)/H\subset X/H$.
Notice that $\Gamma_H$ is a Cayley graph for $G/H$.
We first show how to extend maps of compact surfaces
with boundary components in $\Gamma_H$ to proper maps of non-compact
surfaces.  We then say what we mean by ``lifting and straightening''
such proper maps.

Let $\Sigma$ be a compact planar surface, and let 
\[\phi\co \Sigma\to X/H\]
be a cellular map so that $\phi|_{\partial\Sigma}\subset \Gamma_H$.
\begin{definition}\label{d:parabolicboundary}
If $P\in \mc{P}$, let $\partial_P\Sigma$ be the union of those
boundary components $c$ so that $\phi|_c$ lifts to an arc in $X$
which lies in a single left coset of $P$.
Let 
\[\partial_{\mc{P}}\Sigma=\cup_{P\in\mc{P}}\partial_P\Sigma.\]
be the union of those boundary
components of $\Sigma$, each of which is sent into the image of a single
$0$-horoball of $X$ (see Definition \ref{d:horoball}).  We refer to
$\partial_{\mc{P}}$ as the {\em parabolic boundary} of $\Sigma$.
\end{definition}
\begin{definition}\label{d:extend}
Let 
\[\extend{\Sigma} = \Sigma
\cup_{\partial_{\mc{P}}\Sigma} \partial_{\mc{P}}\Sigma\times[0,\infty)\] 
be the surface obtained from $\Sigma$ by attaching a half-open annulus
to each component of $\partial_{\mc{P}}\Sigma$. We extend $\phi$ to
a proper map 
\[\extend{\phi}\co \extend{\Sigma}\to X/H\]
as follows:  Let $c\cong S^1$ be a component of
$\partial_{\mc{P}}\Sigma$.  
The map $\phi$ sends $c\times\{0\}$ to a loop
$\gamma_0$ which is contained in $\Gamma(P,P\cup S)/P$ for some
$P\in \mc{P}$.  This loop lifts to some path $\twid{\gamma_0}$ in $X$,
each edge of which is the top edge of some vertical square.  Let
$\gamma_1$ be the loop which is the projection of the path obtained by
traversing the bottom edges of those squares, and define
$\extend{\phi}|_{c\times[0,1]}$ to be a homotopy across the images
of those squares in $X/G$ so that
$\extend{\phi}|_{c\times\{1\}}=\gamma_1$.  Similarly define
$\extend{\phi}|_{c\times[k,k+1]}$ for each $k\geq 1$ so that
$D(\extend{\phi}(c\times \{t\}))=t$ for each $t\geq 0$.
\end{definition}

A {\em peripheral path} in $X/H$ is a path in the $1$-skeleton which
lies in the image of the Cayley graph of some $P_i \in \mc{P}$. 

\begin{definition}\label{d:reduce}
Let $\Sigma$ be a compact planar surface, and
let $\phi\co\Sigma\to X/H$ be as above.  A \emph{reducing arc} 
for $\phi$ is an essential,
properly embedded interval $\sigma\co I\to\Sigma$ so that
$\sigma(\partial I)\subset \partial_P\Sigma$ for some $P\in \mc{P}$
and $\phi\circ\sigma$ is homotopic rel endpoints to a peripheral
path.
\end{definition}

\begin{definition}
Analogously we may define reducing arcs for proper maps into $X/G$:
Suppose that $\Xi$ is obtained from a compact planar surface by
removing finitely many points.
Let $f: \Xi \to X/H$ be a proper map, and suppose $\sigma: \R\to
\Xi$ is an essential proper arc.  We say that $\sigma$ is a
reducing arc if $f\circ\sigma$ is properly homotopic into
$D^{-1}[L,\infty)$ for some (and hence for any) $L>0$.   
\end{definition}

\begin{lemma}
The map $\phi \co \Sigma \to X/H$ has a reducing
arc if and only if $\extend{\phi} \co \extend{\Sigma} \to X/H$
has a reducing arc.
\end{lemma}

\begin{lemma}\label{l:straighten}
Let $\Xi$ be a surface of finite type (possibly with boundary),
and let $\mc{T}$ be a triangulation of $\Xi$ (partially ideal
if appropriate).  Let $\psi \co \Xi \to X/H$ be a proper map
satisfying:
(i) $\psi$ has no reducing arcs; and (ii) if $e$ is an edge of $\mc{T}$
which limits on a puncture $p$ of $\Xi$ then there is a 
neighborhood $U$ of $p$ so that $\psi |_{e \cap U}$ is vertical.

Then there is a proper homotopy from $\psi$ to a map
\[\straightedge{\psi}\co\Xi\to X/H \]
so that if $e$ is any edge of $\mc{T}$, then
$\straightedge{\psi}|_{e}$ lifts to a preferred path in $X$.
\end{lemma}
\begin{proof}
Let $e$ be an edge of $\mc{T}$.
It suffices to show that $\psi|_{e}$ is homotopic to a preferred path,
and that this homotopy is level-preserving near the ends of $e$.

Choose a lift $\twid{\psi|_e}$ of $\psi|_e$ to $X$. The map
$\twid{\psi|_e}$ extends to a map from $I$ to $X \cup \Hbound$.
Let $a$ be the image of $0$ and $b$ the image of $1$.  Since
$e$ is not a reducing arc, $a \neq b$ and there is a preferred path
$\pp_{a,b}$ between these points.  Consider $\pp_{a,b}$ as a map
from $e$ to $X$.

Suppose that $a$ is ideal.  Then both $\twid{\psi|_e}$ and $\pp_{a,b}$
are vertical on some initial segment, so we may reparametrize so 
that $D \circ \twid{\psi|_e}$ and $D \circ \pp_{a,b}$ agree on this
initial segment.  Therefore, there is an obvious horizontal homotopy
from $\twid{\psi|_e}$ to $\pp_{a,b}$ on this initial segment.  We consider
the projection of this homotopy to $X/H$.

In case $b$ is ideal, we may similarly perform a homotopy on a terminal
subsegment of $e$.

We are now left to deal with a compact loop, formed by the paths
$\twid{\psi|_e}$ and $\pp_{a,b}$ (or the sub-paths with which we
have not yet dealt).  The space $X$ is simply-connected.
\end{proof}

\begin{remark}
Suppose that $\phi \co \Sigma \to X/H$ is as at the beginning of
this subsection.  The surface $\extend{\Sigma}$ and the map 
$\extend{\phi} \co \extend{\Sigma} \to X/H$ satisfy condition
(ii) of the hypothesis of Lemma \ref{l:straighten}.
\end{remark}

\begin{torsionremark} \label{r:tf6}
In the presence of torsion it is not, in general,
possible to find a map
$\straightedge{\psi}$ as in Lemma \ref{l:straighten} so that each
triangle is mapped to something which is combinatorially controlled
the way that preferred triangles are.  Having no reducing arcs is
not sufficient: An additional hypothesis is required.
\end{torsionremark}

\subsection{The skeleton of a map}
In this subsection, we define the {\em skeleton} of the map of a
surface into $X/H$, assuming the surface has been triangulated by
edges which are sent to preferred paths.

\begin{definition}\label{d:skeleton}
Suppose that $\Xi$ is a surface with a (possibly partially
ideal) triangulation $\mc{T}$, and that $\theta\co\Xi\to X/H$ sends
each edge of $\mc{T}$ to a non-degenerate 
path which lifts to a preferred path 
between points in $X\cup\partial_{\mc{H}} X$.  
Further assume that $\theta$ sends no vertex of $\mc{T}$ to a point in
an $L_2$-horoball.
Let $\bar{\Xi}$ be
the compact surface obtained from $\Xi$ by filling in the
punctures; $\mc{T}$ induces a triangulation $\bar{\mc{T}}$ of
$\bar\Xi$.
Then $\theta$ extends to a map 
\[\bar\theta\co\bar\Xi\to X/H \cup (\Hbound)/H.\]
If $\Delta$ is a triangle of $\bar{\mc{T}}$, then 
$\bar\theta|_{\partial\Delta}$ lifts to a preferred
triangle in $X$ (Definition \ref{d:preftri}).  The skeletal filling 
$\skel{\theta|_{\partial\Delta}}$ can then 
be inscribed on $\bar\Xi$ (Lemma \ref{l:inscribe}). 
The {\em skeleton of $\theta$}, $\skel{\theta}$, 
is the $1$-complex in $\bar\Xi$ which is
the union of the $\skel{\theta|_{\partial\Delta}}$ for 
$\Delta \in \bar{\mc{T}}$.
\end{definition}

\begin{remark}
A typical application of the above definition and the lemmas below
is the situation where $\phi \co \Sigma \to X/H$ is as at the beginning
of the section (with no reducing arcs).  In this situation, we will
take $\Xi$ to be $\extend{\Sigma}$ and $\theta$ to be
$\straightedge{\extend{\phi}}$.
\end{remark}

\begin{remark}\label{r:maponskel}
Let $\theta\co\Xi\to X/H$ and $\mc{T}$ be as described in Definition
\ref{d:skeleton}, and write $\Xi^{(1)}$ for the union of the edges
of $\mc{T}$.
By Remark \ref{r:extend}, we can extend $\theta|_{\Xi^{(1)}}$ to a
map 
\[\ddot{\theta}\co\skel{\theta}\to X/H\cup (\Hbound)/H\]
which collapses each ligament to a point, and sends each rib to a
horizontal edge.  Observe that $D(\ddot{\theta}(x))=L_2$ for any vertex
$v$ of $\skel{\theta}$ not coming from a vertex of $\bar{\mc{T}}$.
\end{remark}

\begin{definition}\label{d:link}
Let $\theta\co\Xi\to X/H$ and $\mc{T}$ be as described in Definition
\ref{d:skeleton}, and let $\ddot{\theta}\co\skel{\theta}\to
X/H\cup(\Hbound)/H$ be as in Remark \ref{r:maponskel}.
Let $\bar{D}\co\skel{\theta}\to[0,\infty]$ be
as follows:
\[\bar D(x) = \left\{\begin{array}{ll}  
\infty & \mbox{ if $x$ is a vertex coming from a puncture.}\\
D(\ddot\theta (x)) & \mbox{ otherwise.}
\end{array}\right.\]
Let $v$ be a vertex of $\skel{\theta}$ coming from
a puncture of $\Xi$.  Let $E(v)\subset\bar{\Xi}$ be the smallest closed
disk containing the component of $\bar D^{-1}[L_2,\infty]$ which also
contains $v$.  The {\em link of $v$}, $\Lk(v)$, is the boundary of
$E(v)$.  
\end{definition}

\begin{remark}
Note that $\Lk(v)$ is contained in the skeleton of $\theta$.  If $\Xi$
is a punctured sphere, then $\Lk(v)$ is a circle made up entirely of
ribs and ligaments; otherwise it may contain parts of
edges of $\mc{T}$.  Also notice that $E(v)$ contains no vertex of
$\bar{\mc{T}}$ other than $v$.  In Figure \ref{fig:lollipop}, one can
see an example showing both kinds of links.
\end{remark}

The next two lemmas follow easily from the fact that $\theta$ restricted
to a triangle of $\bar{\mc{T}}$ lifts to a preferred triangle in $X$,
and that both $X$ and $L_2$-horoballs in $X$ are simply connected.
\begin{lemma}\label{l:skelstraighten}
Let $\theta\co\Xi\to X/H$ and $\mc{T}$ be as described in Definition
\ref{d:skeleton}, and let $\ddot{\theta}\co\skel{\theta}\to
X/H\cup(\Hbound)/H$ be as in Remark \ref{r:maponskel}.  Let 
$\iota\co \skel{\theta}\to \bar\Xi$ be an inclusion which is the
identity on edges in $\bar{\mc{T}}$.  The map $\theta$ is properly
homotopic to a map $\theta'$ so that $\theta'\circ\iota =\ddot \theta$.
\end{lemma}

\begin{lemma} \label{l:link}
Let $\theta\co\Xi\to X/H$ and $\mc{T}$ be as described in Definition
\ref{d:skeleton}, and let $\ddot{\theta}\co\skel{\theta}\to
X/H\cup(\Hbound)/H$ be as in Remark \ref{r:maponskel}.  Let $\sigma\co
S^1\to\Xi$ be a loop surrounding a puncture $x$, so that 
$\theta\circ\sigma(S^1)$ lies entirely in the component of 
$(\theta\circ D)^{-1}[L_2,\infty)$ surrounding $x$.  Let $\gamma\co
  S^1\to \Lk(x)$ be a homeomorphism.
Then
  $\theta\circ\sigma$ is homotopic to $\ddot\theta\circ\gamma$ inside
  $D^{-1}[L_2,\infty)$. 
\end{lemma}

\def\surjects{\twoheadrightarrow}

\section{Punctured spheres and disks} \label{s:spheres}

In this section we investigate relations in $G/K$ amongst
the images of the parabolic elements of $G$.  In particular,
we prove assertion (1) of Theorem \ref{t:rhds}.  We also prove
Theorem \ref{t:dontintersect}.

For this section, we make the following standing assumptions:
\begin{enumerate}
\item $G$ is a finitely generated, torsion-free group, which is hyperbolic
relative to a collection $\mc{P} = \{ P_1, \ldots , P_n \}$ of finitely generated subgroups;
\item $S$ is a finite compatible generating set for $G$ with respect
to $\mc{P}$ (in the sense of Definition \ref{d:compatible});
\item $X(G,\mc{P},S)$ is $\delta$-hyperbolic.
\end{enumerate}

\subsection{From parabolic equations to surfaces and punctured spheres}\label{s:equations}
In this paragraph, we explain how any equation amongst
parabolic words (and their conjugates) may be turned into
a map of a compact planar surface with boundary into $X/G$; 
we can then extend this to
a proper map of a punctured sphere into $X/G$ using Definition 
\ref{d:extend}.

Suppose that in $G$ there is an equality of the form:
\begin{equation} \label{paraproduct}
1 = \prod_{i=1}^m g_i p_i g_i^{-1},	
\end{equation}
where, for each $i$, $p_i \in P_{j_i}$ and $g_i \in G$.
Choosing words for each $g_i$ and words in $S \cap P_{j_i}$
for each $p_i$, we find a map $\tilde{\phi} \co \tilde{\Sigma} \to X$
of a disk, sending the boundary to a loop representing the equation.

Projection gives a map $\pi \circ \tilde{\phi} \co \tilde{\Sigma} \to X/G$.  

\begin{figure}[htbp]
\begin{center}
\begin{picture}(0,0)%
\includegraphics{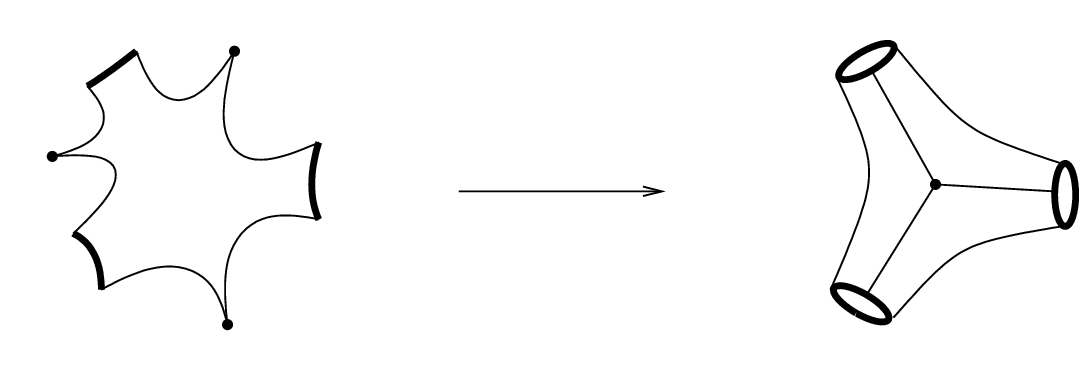}%
\end{picture}%
\setlength{\unitlength}{4144sp}%
\begingroup\makeatletter\ifx\SetFigFont\undefined%
\gdef\SetFigFont#1#2#3#4#5{%
  \reset@font\fontsize{#1}{#2pt}%
  \fontfamily{#3}\fontseries{#4}\fontshape{#5}%
  \selectfont}%
\fi\endgroup%
\begin{picture}(4942,1689)(2461,-1548)
\put(2476,-627){\makebox(0,0)[lb]{\smash{{\SetFigFont{9}{10.8}{\familydefault}{\mddefault}{\updefault}{\color[rgb]{0,0,0}$1$}%
}}}}
\put(3566, 14){\makebox(0,0)[lb]{\smash{{\SetFigFont{9}{10.8}{\familydefault}{\mddefault}{\updefault}{\color[rgb]{0,0,0}$g_1k_1g_1^{-1}$}%
}}}}
\put(3534,-1492){\makebox(0,0)[lb]{\smash{{\SetFigFont{9}{10.8}{\familydefault}{\mddefault}{\updefault}{\color[rgb]{0,0,0}$g_1k_1g_1^{-1}g_2k_2g_2^{-1}$}%
}}}}
\end{picture}%
\caption{Identifying those parts of the disk $\tilde{\Sigma}$ labelled
by the conjugating elements yields a compact planar surface which maps
into $X/G$.}
\label{f:realize}
\end{center}
\end{figure}
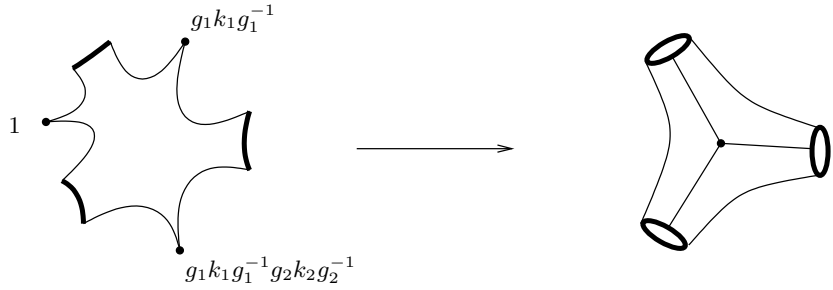

Let $\Sigma$ be the surface obtained from the disk $\tilde{\Sigma}$
by identifying those parts of the boundary corresponding to the $g_j$
in pairs (Figure \ref{f:realize}).

The map $\pi \circ \tilde{\phi}$ factors through the quotient map from
$\tilde{\Sigma}$ to $\Sigma$.  Let $\phi \co \Sigma \to X/G$ be the 
resulting map.

The $m$ boundary components of $\Sigma$ are sent by $\phi$ to 
peripheral loops
in $X/G$ whose labels are the words chosen for the $p_i$ above.

Conversely, given a compact planar surface with all but (possibly) one boundary component labelled by elements of the $K_i$, we can reconstruct an equality like \eqref{paraproduct} for the word representing the other boundary component, by choosing the $g_j$
via some paths through the surface.

\subsection{The groups $P_i/K_i$ inject into $G/K$}

\begin{theorem} \label{t:injects}
There is a constant $R=R(\delta)\leq 12\cdot 2^{3000\delta}$ so that
the following holds:
If $\{ K_1, \ldots , K_n \}$ are so that
$K_i\unlhd P_i$ and $|K_i|_{P_i} > R$
then the natural map $\iota_i \co P_i/K_i\to G/K$ is an injection,
where $K$ is the normal closure in $G$ of 
$K_1 \cup \cdots \cup K_n$.
\end{theorem}
\begin{proof}(cf. \cite[Proof of Theorem 3.1]{lackenby:whds})

We will show that if, for some $i$, the map $\iota_i$ is
not injective, then $|K_l|_{P_l}$ must be small for some $l$.

Let $\alpha\in P_i\setminus
K_i$ be an element of $K$. Thus there is some equation
\begin{equation}\label{alphaproduct}
\alpha = \prod_{j=1}^m g_j k_j g_j^{-1}
\end{equation}
for some finite sequence of $g_j\in G$ and $k_j\in K_{i_j}$.
We say that such an equation {\em represents the death of $\alpha$ in
$G/K$}.
We may suppose that the product in \eqref{alphaproduct} is minimal in
the following sense:  If $\alpha'$ is \emph{any} element of
$\cup_i P_i\setminus \cup_i K_i$, and 
\[\alpha' = \prod_{r=1}^{m'} g_r k_r g_r^{-1}\] for some collection of
$g_r\in G$ and $k_r\in K_{i_r}$, then $m'\geq m$.  (In other words, the
expression in \eqref{alphaproduct} is minimal in length not only for
$\alpha$, but over all such equations with the left hand side an
element of $\cup_i P_i\setminus \cup_i K_i$.)

The equality \eqref{alphaproduct} is realized by a map $\phi\co 
\Sigma\to X/G$, as in Subsection \ref{s:equations}.

\begin{claim} \label{Claim1}
The map $\phi$ has no reducing arcs.
\end{claim}
\begin{proof}
Let $\sigma$ be a reducing arc.  Each of the boundary components
of $\Sigma$ has a corresponding word, and thus a corresponding base point (the starting point of this word).  We may assume that $\sigma$ starts and finishes at one of these base points.  The path $\sigma$
determines some element  $p$ of $G$.  Because $\sigma$ is homotopic
to a peripheral path, $p$ is contained in some $P_k$.  We suppose 
$\phi$ has been homotoped so that $\phi \circ \sigma$ is a peripheral 
path.

There are five cases to consider, depending on the endpoints of $\sigma$.

\begin{case} \label{ktok}
Suppose that the endpoints of $\sigma$ are on different boundary components, neither of which corresponds to $\alpha$.  Cutting along
$\sigma$ yields a new surface $\Sigma'$ with fewer boundary components
than $\Sigma$ (see Figure \ref{f:ktok});  
\begin{figure}[htbp]
\begin{center}
\begin{picture}(0,0)%
\includegraphics{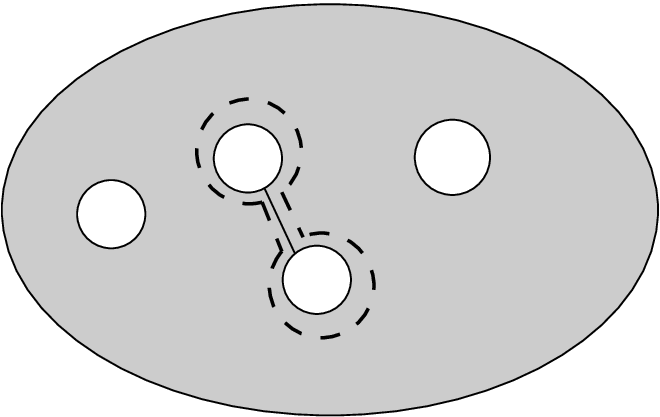}%
\end{picture}%
\setlength{\unitlength}{4144sp}%
\begingroup\makeatletter\ifx\SetFigFont\undefined%
\gdef\SetFigFont#1#2#3#4#5{%
  \reset@font\fontsize{#1}{#2pt}%
  \fontfamily{#3}\fontseries{#4}\fontshape{#5}%
  \selectfont}%
\fi\endgroup%
\begin{picture}(3016,1894)(1143,-3698)
\end{picture}%
\caption{Reducing the number of boundary components labelled with
elements of $K_i$.}
\label{f:ktok}
\end{center}
\end{figure}
two have been removed, and the new one is sent to a
word representing $k_ipk_jp^{-1}$ for some $i,j$. This contradicts the
minimality assumption. 
\end{case}
\begin{case}\label{alphatok}
Suppose next that the initial point of $\sigma$ lies on the boundary component corresponding to $\alpha$, and that the terminal point
lies on the
boundary component corresponding to $k_i$ for some $i$.   Cutting along $\sigma$ again yields a surface $\Sigma'$ with one fewer
boundary component than $\Sigma$.  The boundary components
corresponding to $\alpha$ and $k_i$ have been removed.  The new boundary component has label $\alpha' = \alpha p k_i p^{-1}$.  Since 
$K_i \unlhd P_i$, and $\alpha \not\in K_i$, we also have $\alpha'
\not\in K_i$.  This again contradicts the minimality of equation 
\eqref{alphaproduct}.
\end{case}
\begin{case} \label{alphatoalpha}
Suppose that both endpoints of $\sigma$ lie on the boundary
component corresponding to $\alpha$.  Cutting $\Sigma$ along 
$\sigma$ yields two compact planar surfaces, each with fewer
boundary components than $\Sigma$.  One of them represents
the death of $p$ in $G/K$, whilst the other represents the death
of $p\alpha$.

In case $p \not\in K_i$, the equation for $p$ contradicts minimality.  
Thus we may suppose that $p \in K_i$, in which case $p\alpha \not\in
K_i$, and again we get a contradiction to minimality.
\end{case}
\begin{case} \label{ktosamek}
Finally, suppose that $\sigma$ is a loop whose endpoints are
on the boundary component corresponding to $k_i$ (as in Figure
\ref{f:ktosamek}).  
\begin{figure}[htbp]
\begin{center}
\begin{picture}(0,0)%
\includegraphics{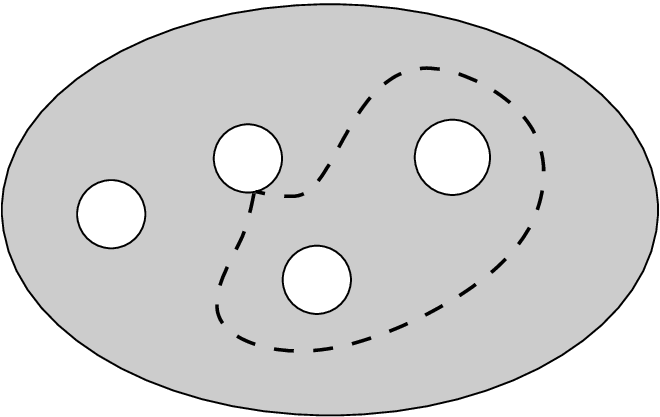}%
\end{picture}%
\setlength{\unitlength}{4144sp}%
\begingroup\makeatletter\ifx\SetFigFont\undefined%
\gdef\SetFigFont#1#2#3#4#5{%
  \reset@font\fontsize{#1}{#2pt}%
  \fontfamily{#3}\fontseries{#4}\fontshape{#5}%
  \selectfont}%
\fi\endgroup%
\begin{picture}(3016,1894)(1143,-3698)
\put(2648,-2164){\makebox(0,0)[lb]{\smash{{\SetFigFont{12}{14.4}{\familydefault}{\mddefault}{\updefault}{\color[rgb]{0,0,0}$\sigma$}%
}}}}
\end{picture}%
\caption{Case \ref{ktosamek} of Claim \ref{Claim1}.}
\label{f:ktosamek}
\end{center}
\end{figure}
Cutting along
$\sigma$ yields a pair of surfaces, both of which have fewer boundary
components than $\Sigma$.  If $p \not\in K_{j_i}$ then the surface $\Sigma'$
not containing $\alpha$ corresponds to an equation representing the
death of $p$ or $pk_i$  (depending on whether $k_i$ is contained in $\Sigma'$).  On the other hand, if 
$p \in K_{j_i}$ then we can replace the surface not containing $\alpha$
by a single puncture, to yield a surface $\Sigma''$ with fewer
punctures than $\Sigma$, representing the death of $\alpha$.

In either case, this contradicts the minimality of $m$.
\end{case}
This proves Claim \ref{Claim1}
\end{proof}

As in Section \ref{s:extensions}, let 
$\extend{\Sigma} = \Sigma\cup \partial\Sigma\times[0,\infty)$, and let
\[\extend{\phi}\co \extend{\Sigma}\to X/K\]
be the extension from Definition \ref{d:extend}.

We now choose an ideal triangulation $\mc{T}$ of the punctured sphere
$\extend{\Sigma}$.  As $\extend{\Sigma}$ is an $m+1$-times punctured
sphere, $\mc{T}$ contains $2m-2$ triangles.
Let
\[\straightedge{\extend{\phi}}\co\extend{\Sigma}\to X/K\]
be the map from Lemma \ref{l:straighten}, which sends each edge of
$\mc{T}$ to a preferred path, and let 
\[ \straightedge{\ddot{\extend{\phi}}} \co \skel(\straightedge{\extend{\phi}}) \to X/G \cup (\Hbound)/G, \] 
be as in Remark \ref{r:maponskel}.

If $T\in\mathcal{T}$, then $\phi|_{\partial T}$ lifts to a preferred
triangle $\twid{\phi_T}\co\partial T\to X$.  Let $R(T)$ be the number
of ribs in $\skel(\twid{\phi_T})$, and note that this number does not
depend on the lift chosen.
Corollary \ref{c:ribsbound} implies that $R(T)\leq 6$.

Let 
\[A(\phi) = \sum_{T\in \mathcal{T}} R(T).\]
Corollary \ref{c:ribsbound} immediately implies
\begin{equation}\label{upperbound}
A(\phi)\leq 6 (2m-2)< 12 m.
\end{equation}

Let $x$ be one of the punctures of $\extend{\Sigma}$ not corresponding
to $\alpha$.  By Lemma \ref{l:link}, $\ddot{\extend{\phi}}|_{Lk(x)}$ is
a loop at the $L_2$-level which, considering the $L_2$-level to be the image of a Cayley
graph for $P_i$, represents an element $k$ of $K_i$.

It follows that $Lk(x)$ must contain at least $2^{-L_2}|k|_{P_i}$ ribs.

Thus 
\[	A(\phi) \ge 2^{-L_2} \Big( \min_{s} |K_{i_s}|_{P_{i_s}} \Big) m .
\]

Therefore, by \eqref{upperbound},
\[	\min_s |K_{i_s}|_{P_{i_s}} \le 2^{L_2} \frac{A(\phi)}{m} 
< 12 \cdot 2^{L_2}.
\]
\end{proof}

\subsection{On the structure of the quotients $G/K$}

The theme in this subsection is that, by choosing large enough
algebraic slope lengths, we can preserve much of the structure
of $G$ in its quotients.  See \cite{GMO} for more results along
these lines.

\begin{theorem} \label{t:dontintersect}
Suppose $G$ is torsion-free and that $|K_i|_{P_i} > 12\cdot 2^{L_2}$ for each $i$.  Let $P$ and $P'$ be conjugate into $\mc{P}$
and suppose that whenever $P^g = P'$ then $g \not\in K$.  Then
the images of $P$ and $P'$ in $G/K$ intersect trivially.
\end{theorem}
\begin{proof}
Suppose that $P$ is conjugate to $P_j \in \mc{P}$, and let $K_P$ be the conjugate of $K_j$
in $P$.  Define $K_{P'}$ similarly.
If the theorem is false, then there is an equation in $G$ of the form
\begin{equation}\label{annuluseq}
	q = q' \prod_{i=1}^m g_i k_i g_i^{-1}	,	
\end{equation}
where $q \in P \smallsetminus K_P$ and $q' \in P' \smallsetminus
K_{P'}$.  Suppose we have chosen such an equation with $m$ minimal
over all such equations (over all choices of $P$ and $P'$).

This gives rise to a map $\phi$ of a compact planar surface with $m+2$ boundary
components into $X/G$. Once again, we claim that this map has no
reducing arcs.

\begin{figure}[htbp]
\begin{center}
\begin{picture}(0,0)%
\includegraphics{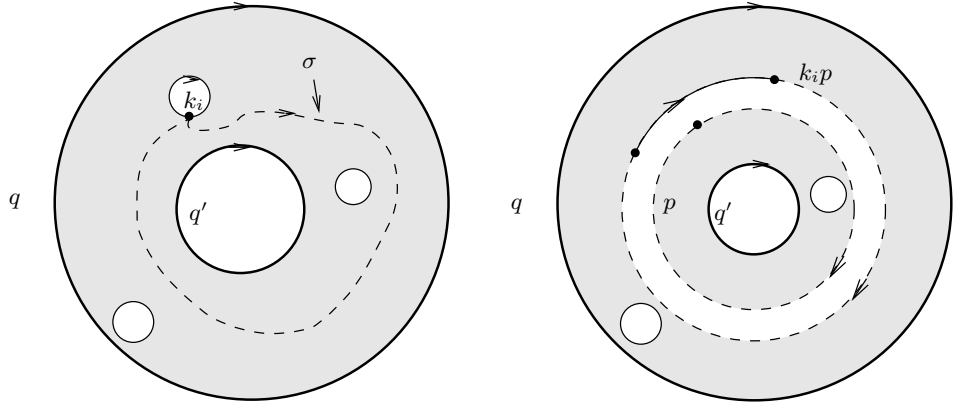}%
\end{picture}%
\setlength{\unitlength}{4144sp}%
\begingroup\makeatletter\ifx\SetFigFont\undefined%
\gdef\SetFigFont#1#2#3#4#5{%
  \reset@font\fontsize{#1}{#2pt}%
  \fontfamily{#3}\fontseries{#4}\fontshape{#5}%
  \selectfont}%
\fi\endgroup%
\begin{picture}(5667,2408)(1651,-6698)
\put(5885,-5607){\makebox(0,0)[lb]{\smash{{\SetFigFont{9}{10.8}{\familydefault}{\mddefault}{\updefault}{\color[rgb]{0,0,0}$q'$}%
}}}}
\put(2746,-5607){\makebox(0,0)[lb]{\smash{{\SetFigFont{9}{10.8}{\familydefault}{\mddefault}{\updefault}{\color[rgb]{0,0,0}$q'$}%
}}}}
\put(2712,-4932){\makebox(0,0)[lb]{\smash{{\SetFigFont{9}{10.8}{\familydefault}{\mddefault}{\updefault}{\color[rgb]{0,0,0}$k_i$}%
}}}}
\put(3421,-4696){\makebox(0,0)[lb]{\smash{{\SetFigFont{9}{10.8}{\familydefault}{\mddefault}{\updefault}{\color[rgb]{0,0,0}$\sigma$}%
}}}}
\put(4670,-5540){\makebox(0,0)[lb]{\smash{{\SetFigFont{9}{10.8}{\familydefault}{\mddefault}{\updefault}{\color[rgb]{0,0,0}$q$}%
}}}}
\put(1666,-5506){\makebox(0,0)[lb]{\smash{{\SetFigFont{9}{10.8}{\familydefault}{\mddefault}{\updefault}{\color[rgb]{0,0,0}$q$}%
}}}}
\put(5581,-5540){\makebox(0,0)[lb]{\smash{{\SetFigFont{9}{10.8}{\familydefault}{\mddefault}{\updefault}{\color[rgb]{0,0,0}$p$}%
}}}}
\put(6391,-4764){\makebox(0,0)[lb]{\smash{{\SetFigFont{9}{10.8}{\familydefault}{\mddefault}{\updefault}{\color[rgb]{0,0,0}$k_ip$}%
}}}}
\end{picture}%
\caption{This kind of reducing arc gives rise to at least one new
equation of type \eqref{annuluseq}, if $p\not\in K_{j_i}$.}
\label{f:annulusreduce}
\end{center}
\end{figure}

Let $\sigma$ be a reducing arc which is a peripheral arc, and let $p$
be the corresponding element of $G$, contained in $K_i$, say.  Most of 
the cases are entirely analogous to those in the proof
of Theorem \ref{t:injects}; they all lead to a contradiction to minimality, or 
to injectivity.  We deal with the most interesting case, when $\sigma$ is
a loop whose endpoints lie on a puncture corresponding to some
$k_i$, and so that cutting along $\sigma$ separates the puncture
corresponding to $q$ from the puncture corresponding to $q'$.  See
Figure \ref{f:annulusreduce}.

By Theorem \ref{t:injects}, $q$ and $q'$ are nontrivial in $G/K$.  
Therefore, $p \not\in K_i$.  Now, the hypothesis of the theorem
implies that it cannot be that there are $k, k' \in K$ so that $P = P_i^k$ and $P' = P_i^k$.  Therefore, one of the two diagrams
obtained by cutting along $\sigma$ yields a contradiction to
the minimality of $m$.  This, and the omitted cases, shows that there
are no reducing arcs.

We now proceed as in the proof of Theorem \ref{t:injects}.  The only
difference is that we now have $2m$ triangles, rather than $2m-2$.
However, it is still the case that $\skel(\straightedge{\extend{\phi}})$
has at most $12m$ ribs, and therefore
\[\min_s |K_{i_s}|_{P_{i_s}} \le 12 \cdot 2^{L_2},\]
as required.
\end{proof}

\begin{corollary} \label{c:ijdisjoint}
Under the hypotheses of Theorem \ref{t:dontintersect}, if $i \neq j$ then
\[	\iota_i(P_i/K_i) \cap \iota_j(P_j/K_j) = \{ 1 \}	.	\]
\end{corollary}

\begin{corollary} \label{c:conjugatesdisjoint}
Under the hypotheses of Theorem \ref{t:dontintersect}, $\iota_i(P_i/K_i)$
is malnormal in $G/K$, for each $i$.
\end{corollary}

\begin{proposition} \label{p:thickdontdie}
Suppose that $G$ is torsion-free and that $|K_i|_{P_i} >
12\cdot 2^{L_2}$ for each $i$.
Suppose that $g \in G \smallsetminus \{ 1 \}$ is such that there exists $x \in X$ so that
the preferred path $\pp_{x,g.x}$ lies entirely within $D^{-1}([0,L_2-1])$.
Then $g \not\in K$.
\end{proposition}
\begin{proof}
As usual, we suppose that the theorem is false, build a surface in
$X/K$, and use its geometry to derive a contradiction.

If $g \in K$ then there is an equality in $G$ of the form
\begin{equation} \label{e:gexpression}
1 = g^{-1} \prod_{i=1}^m g_i k_i g_i^{-1}	,
\end{equation}
where $k_i \in K_{j_i}$ and $g_i \in G$.  We choose such an 
expression for $g$ which minimizes $m$. Since $g \neq 1$ in $G$,
we have $m \ge 1$.

Choosing words for each $g_i$, each $k_i$ and $g$,
the expression for $g$ in \eqref{e:gexpression} gives a loop in $X$,
beginning at $1$.
In turn, this induces a map $\tilde\phi \co \tilde \Sigma \to X$ of a disk, as
described in Subsection \ref{s:equations}.  Also as in Subsection
\ref{s:equations}, we may glue $\tilde\Sigma$ along parts of the boundary corresponding to the $g_i$, to obtain a compact
planar surface, $\Sigma$, together with a map $\phi' \co \Sigma \to
X/K$.  The surface $\Sigma$ has one distinguished boundary component which is labelled by a word representing $g$, and we
call this component the {\em $g$-boundary} of $\Sigma$.

We claim that there are no reducing arcs for $\phi$ whose endpoints
do not lie on the $g$-boundary.  This follows as before: Any
such reducing arc either contradicts Theorem \ref{t:injects} or
else the minimality of $m$.

Let $\xi$ be a simple path in $X$ from $1$ to $x$, and consider the
loop in $X$ which is the concatenation $\xi \cdot \pp_{x,g x} \cdot
\mu \cdot \eta$, where $\mu$ is the path $g \xi$ traversed backwards,
and $\eta$ is a lift to $X$ of the image under $\phi'$ of the $g$-boundary
of $\Sigma$.

This loop may be filled with a disk in $X$, which projects to an
annulus in $X/K$, which has one boundary component the 
image of the $g$-boundary of $\Sigma$.  

Homotoping $\phi'$ over the annulus gives a new map 
$\phi \co \Sigma \to X$ which maps the $g$-boundary to the
image of the preferred path $\pp_{x,g.x}$ in $X/K$. 
See Figure \ref{f:glove}.
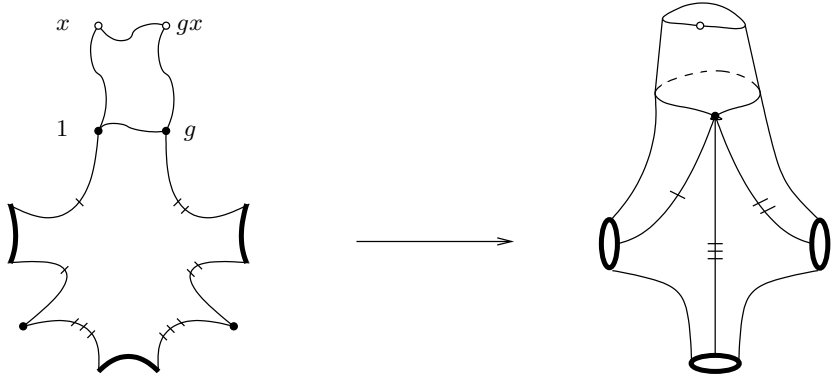
\begin{figure}[htbp]
\begin{center}
\begin{picture}(0,0)%
\includegraphics{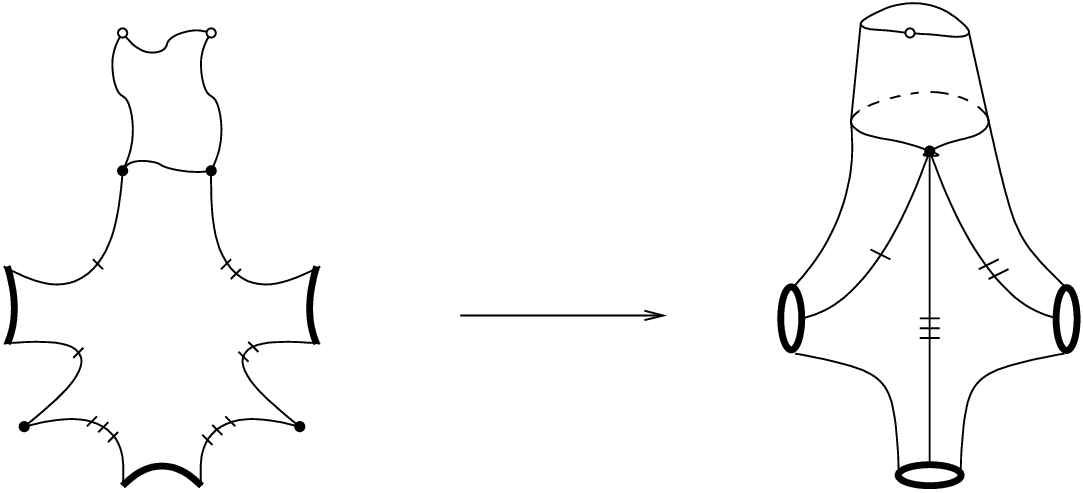}%
\end{picture}%
\setlength{\unitlength}{4144sp}%
\begingroup\makeatletter\ifx\SetFigFont\undefined%
\gdef\SetFigFont#1#2#3#4#5{%
  \reset@font\fontsize{#1}{#2pt}%
  \fontfamily{#3}\fontseries{#4}\fontshape{#5}%
  \selectfont}%
\fi\endgroup%
\begin{picture}(4947,2251)(2456,-1534)
\put(2746,-106){\makebox(0,0)[lb]{\smash{{\SetFigFont{9}{10.8}{\familydefault}{\mddefault}{\updefault}{\color[rgb]{0,0,0}$1$}%
}}}}
\put(3511,-106){\makebox(0,0)[lb]{\smash{{\SetFigFont{9}{10.8}{\familydefault}{\mddefault}{\updefault}{\color[rgb]{0,0,0}$g$}%
}}}}
\put(2746,524){\makebox(0,0)[lb]{\smash{{\SetFigFont{9}{10.8}{\familydefault}{\mddefault}{\updefault}{\color[rgb]{0,0,0}$x$}%
}}}}
\put(3466,524){\makebox(0,0)[lb]{\smash{{\SetFigFont{9}{10.8}{\familydefault}{\mddefault}{\updefault}{\color[rgb]{0,0,0}$gx$}%
}}}}
\end{picture}%
\caption{The surface $\Sigma$.}
\label{f:glove}
\end{center}
\end{figure}

We now define a surface $\extend\Sigma$ almost as in Definition
\ref{d:extend} by attaching a half-open annulus to each
component of $\partial\Sigma$ {\em except} for the $g$-boundary.
This gives a map
\[	\extend\phi \co \extend\Sigma \to X/K	.	\]
The $g$-boundary of $\extend\Sigma$ is defined in the obvious way.

Choose a (partially ideal) triangulation $\mc{T}$ of $\extend\Sigma$ which has one vertex which is the preimage of $x$, and all other 
vertices ideal, and one edge the $g$-boundary of $\extend\Sigma$.

Homotope $\extend\phi$ to 
\[	\straightedge{\extend\phi} \co \extend\Sigma \to X/K	.\]
as in Lemma \ref{l:straighten}.

Because no edge in $\mc{T}$ is a reducing arc for 
$\straightedge{\extend{\phi}}$, the image of each triangle $T \in \mc{T}$
lifts to a preferred triangle $T_{abc}$ in $X$.
Therefore, we can define the map
\[ \straightedge{\ddot{\extend{\phi}}} \co \skel(\straightedge{\extend{\phi}}) \to X/G \cup (\Hbound)/G, \] 
as in Remark \ref{r:maponskel}.

We now proceed as in the proof of Theorem \ref{t:injects}.  By Euler
characteristic, there are $2m-1$ triangles.  Therefore,
$\skel(\straightedge{\extend{\phi}})$ has less than $12m$ ribs.

Because $\pp_{x,g.x}$ lies entirely within $D^{-1}([0,L_2-1])$, none
of the ribs intersect the $g$-boundary of 
$\skel(\straightedge{\extend{\phi}})$ (where the $g$-boundary of
$\skel(\straightedge{\extend{\phi}})$ is defined in the obvious way).
By Lemma \ref{l:link} there must be at least 
\[	\min_s |K_{i_s}|_{P_{i_s}} m,	\]
ribs in $\skel(\straightedge{\extend{\phi}})$.  Thus we have
\[	\min_s |K_{i_s}|_{P_{i_s}} \le 12 \cdot 2^{L_2} , \]
which is the required contradiction.
\end{proof}

This proposition has the following corollary:

\begin{corollary} \label{c:injectonfinite}
Let $F\subset G$ be a finite set.  Then there is a constant $C=C(F)$
so that if each $|K_i|_{P_i}>C$, then the quotient map $G\to G/K$ is
injective on $F$.
\end{corollary}
\begin{proof}
We prove the equivalent statement that there exists some
$C'=C'(F)$ so that $G\to G/K$ sends no element of
$F\smallsetminus\{1\}$ to $1$.  (To see the equivalence notice that
$C(F)\leq C'((F\cup F^{-1})^2)$.)

Without loss of generality, we may assume that any element of $F$
which is conjugate into some $P_i$ actually lies in $P_i$.  Write
$F=F_1\cup F_2$, where $F_1=F\cap(\cup\mc{P})$ is the set of parabolic
elements of $F$, and $F_2 = F\smallsetminus F_1$ is the set of
non-parabolic elements of $F$.  Let $S' = S\cup F_2$.  Note that since
$S$ is a compatible generating set, so is $S'$.  We may thus form the
space $X' = X(G,\mc{P},S')$.  By Theorem \ref{t:tfae}, $X'$ is
$\delta'$-hyperbolic for some $\delta'>0$.  

Since $S'\cap P_i=S\cap
P_i$, the meaning of $|K_i|_{P_i}$ does not change in moving from $X$
to $X'$.  We may thus apply Proposition \ref{p:thickdontdie} in the
context of $X'$ rather than $X$.  The preferred path joining $1$ to $f$
in $X'$ is a single edge for any $f\in F_2$.  Let 
$C'' = 24\cdot 2^{3000\delta'}$. 
Proposition \ref{p:thickdontdie} implies that if $|K_i|_{P_i}>C''$ for
each $i$, then no element of $F_2$ is sent to $1\in G/K$.  Each $f\in
F_1$ is contained in some $P_{i_f}$.
Let 
\[C'''=\max\{|f|_{P_{i_f}\cap S}\mid f\in F_1\}\]
and let $C' =\max\{ C'',C'''\}$.
\end{proof}

Using exactly analogous arguments to those in Proposition
\ref{p:thickdontdie} we can obtain the following result.\begin{proposition}\label{p:thicknotparabolic}
Suppose that $G$ is torsion-free, and that $|K_i|_{P_i}> 18\cdot
2^{L_2}$ for each $i$.  Suppose that $g\in G$ is not conjugate into
any $P_i\in \mc{P}$, and that there is some $x\in X$ with
$\pp_{x,gx}\subset D^{-1}[0,L_2-1]$.  Then there is no $h\in G$, $p\in
P\in \mc{P}$ so that $hgh^{-1}p^{-1}\in K$.
\end{proposition}
\begin{proof}
The constant
for $|K_i|_{P_i}$ changes because we consider a disk with several
punctures.  However, we have no control of the length of one of the
punctures (that corresponding to $p$).

Therefore, if
there are $m$ punctures corresponding to elements of the $K_i$,
then there are $2(m+1)-1 = 2m+1 \le 3m$ triangles.

Otherwise, the proof is just as before.
\end{proof}

\begin{remark}
Once we have proved that $G/K$ is hyperbolic relative to the images of
the $P_i$ in $\mc{P}$, this will imply that $g$ does not project to a
parabolic element of $G/K$.
\end{remark}

\begin{torsionremark} \label{r:tf7}
In the presence of torsion, many of the arguments in this section
(and later sections) become more difficult to implement.  In
particular, the notion of minimality for surfaces needs to be refined.
Also, some of the results in this section need to be reformulated
in the presence of torsion.  The main issue here is that parabolics
in $G$ may already intersect nontrivially, and this causes a number
of problems.
\end{torsionremark} 

\section{The surgered space}\label{surgered}
In this section we make the following assumptions:
$G$ is a torsion-free group which is 
hyperbolic relative to a collection $\mc{P} =
\{ P_1, \ldots , P_n \}$ of finitely generated
subgroups.  The finite set $S$ is a compatible generating set for
$G$ with respect to $\mc{P}$.
Finally, $\langle S,\mc{P}\mid\mc{R}\rangle$ is a finite relative
presentation for $G$.

Recall from Theorem \ref{t:tfae}, that under these
assumptions, the space $X(G,\mc{P},S)$  is
$\delta$-hyperbolic for some $\delta$.  Moreover,
given the finite relative presentation 
$\langle S,\mc{P}\mid\mc{R}\rangle$, we may adjoin $2$-cells to
$X(G,\mc{P},S)$ 
to form a two-complex
$X=X(G,\mc{P},S,\mc{R})$ which satisfies a linear combinatorial
isoperimetric inequality.

We further suppose that, for $1 \le i \le n$, we have 
$K_i \unlhd P_i$.  Let $K\unlhd G$ be the normal closure of the union
of the $K_i$.  In order to
construct a nice model space for $\Gmk$ it is essential that the
$P_i/K_i$ inject into $\Gmk$.  We thus make the standing assumption
that:
\begin{assumption}\label{K}
For each $1\le i\le n$, $|K_i|_{P_i}> 12\cdot 2^{L_2}$.
\end{assumption}

The space $X$ contains
a copy of the Cayley graph $\Gamma (G,S)$ of $G$, and
has an associated depth function $D\co X \to \R_+$ so that
$\Gamma (G,S) \subset D^{-1}(0)$.

We now describe a ``neutered'' version of $X$, and how to modify it to
get a model for $\Gmk$.
\begin{definition}\label{d:Y}
Let $Y = D^{-1}([0,L_2])$.  
The {\em boundary of $Y$} is $\partial Y = 
D^{-1}(L_2)$.   If $H<G$ then the boundary of $Y/H$ is the image of
$\partial Y$ in $Y/H$.
\end{definition}

\begin{remark}
If the parabolics are not finitely presented, then $Y$ will not be
simply connected.  Its fundamental group is generated by those loops in
$D^{-1}(L_2)$ which cannot be filled in $D^{-1}(L_2)$.  
\end{remark}

The boundary of $Y/G$ has $n$ connected components
which correspond to the subgroups $P_1, \ldots ,
P_n$, as described in the next few paragraphs.

Let $\bar{1} \in X/G$ be the image of the vertices of 
$\Gamma (G,S)$ in $X/G$.  For each $1 \le i \le n$ 
there is a unique $L_2$-horoball in $X$ stabilized by $P_i$,
which we denote $H_i$.
There is a unique vertical path $\tilde{\gamma}_i$ 
joining $1$ to $(i,1,1,L_2) \in H_i$.

Let $\gamma_i$ be the image of $\tilde{\gamma}_i$
in $Y/G \subset X/G$, and let $\gamma^{-1}_i$ be
$\gamma_i$ in the opposite direction.  
Let $T_i \subset Y/G$ be the 
component of $\partial (Y/G)$ containing the image $y_i$
of $(i,1,1,L_2)$.  The vertical path $\gamma_i$ induces an inclusion
of $\pi_1(T_i,y_i)$ into $\pi_1(Y/G,\bar{1})$.

Any loop in $Y/G \subset X/G$ based at $\bar{1}$ lifts to a unique path
in $X$ starting at $1$ and ending at some group element.  This gives a 
well-defined homomorphism from $\pi_1(Y/G,\bar{1})$ onto $G$, which
maps $\pi_1(T_i)$ onto $P_i$.

The next two lemmas follow from Assumption \ref{K}, together with
Theorem \ref{t:injects}.

\begin{lemma} \label{l:noshortloops}
If $c$ is any loop in a boundary component of $Y/K$ of length less than 
$12 \le |K_i|_{P_i}/2^{L_2}$, then $c$ lifts to a loop in $\partial
Y$.  

In particular, there is no loop in $\partial (Y/K)$ consisting of a
single edge.
\end{lemma}

\begin{definition} \label{d:surgered}
Let $Z=Z(K)$ be the $2$-complex obtained from $Y/K$ by gluing a
combinatorial horoball
onto each component of the boundary of $Y/K$.
\end{definition}

\begin{remark}
The depth function on $Y/K$ naturally extends to a depth function $D$
on $Z$.
\end{remark}

\begin{lemma}
If each $|K_i|_{P_i}> 2^{L_2}$,
the complex $Z$ is simply connected.
\end{lemma}
\begin{proof}
The fundamental group of $Y$ is generated by peripheral loops.
Note that we have an exact sequence:
\[ 1\longrightarrow \pi_1(Y)\longrightarrow \pi_1(Y/K) \longrightarrow K \longrightarrow 1\]
Thus, the fundamental group of $Y/K$ is generated by the fundamental
group of $Y$ together with peripheral loops in $Y/K$ which represent
elements of $\pi_1(Y/K)$ which map to normal generators of $K$.

Any peripheral loop eventually dies in a horoball, by Lemma
\ref{l:noshortloops} and Proposition \ref{p:combhorolii}.  
\end{proof}

\begin{lemma}
$G/K$ acts freely and properly on $Z$.
\end{lemma}

We thus have the following diagram of spaces:
\cdlabel{spaces}{
X\ar[d] & Y\ar[l]\ar[d]\\
X/K\ar[d] & Y/K\ar[r]\ar[l]\ar[d] & Z\ar[d] \\                       
X/G  & Y/G \ar[r]\ar[l] & Z/(G/K)}
Each horizontal arrow in \eqref{spaces} is an inclusion, and each
vertical arrow is a covering map.

We will show in Section \ref{section:lollipop} that $Z$
satisfies a linear isoperimetric inequality if the $|K_i|_{P_i}$
are sufficiently large.  Together with
Proposition \ref{ZisanX} below, this will imply that $\Gmk$ is
hyperbolic relative to $\mc{P}'=\{P_1/K_1,\ldots,P_n/K_n\}$, and
hence complete the proof of Theorem \ref{t:rhds}.

\begin{lemma}\label{ZintoXmodK}
There is a $\Gmk$-equivariant embedding $\rho\co Z^{(1)}\to (X/K)^{(1)}$.
\end{lemma}
\begin{proof}
The spaces $Z$ and $X/K$ are identical at depths less than or equal
to $L_2$.  Thus there is an obvious map at these depths.  It is
also obvious that vertical edges in $Z$ can be uniquely associated
to vertical edges in $X/K$.

A horizontal edge $e$ in $Z$ at depth $L_2 + L$ corresponds to a path $p$ at depth $L_2$ in $Z$ of length at most $2^L$.  This path has already been mapped to a path $p'$ in $X/K$ at depth $L_2$ (still of length at most $2^L$).
The path $p'$ lifts to a path $\tilde{p'}$ in $X$ (which is still of length
at most $2^L$).  The path $\tilde{p'}$ lies above an edge $\tilde{e'}$ in 
$X$ at depth $L_2+L$, which projects to an edge $e'$ in $X/K$.  This is
$\rho(e)$.
\end{proof}

Let $\mc{P}' = \{ \iota_i(P_i/K_i) \mid 1\leq i\leq n\}$.

\begin{observation} \label{o:ccomp}
$D^{-1}(0) \subset Y/K \subset X/K$ is a relative Cayley complex for
$G/K$ with respect to the finite relative presentation 
$\langle S, \mc{P}' \mid \mc{R} \rangle$.
\end{observation}

Recall the notation $\mc{H}_N$ from Definition \ref{d:toppart}, where
$\mc{H}$ is a combinatorial horoball.
\begin{lemma} \label{l:shell}
Let $\tilde A$ be a $0$-horoball in $X$, and $A$ its projection in $X/K$.
Then the intersection $A \cap Y/K$ is isomorphic to $\mc{H}_{L_2}$, where $\mc{H} = \mc{H}(A \cap D^{-1}(0))$.
\end{lemma}
\begin{proof}
We define maps $\eta_1 \co A \cap Y/K \to \mc{H}_{L_2}$ and 
$\eta_2 \co \mc{H}_{L_2} \to A \cap Y/K$.

There are obvious bijections on the $0$-skeleta, which extend to
isomorphisms at the $0$level, and the vertical edges.

The first thing to note is that horizontal edges in $A \cap Y/K$ are not loops, by Assumption \ref{K} and Theorem \ref{t:injects}.  It is also true
that horizontal edges in $\mc{H}_{L_2}$ are not loops.

Consider a horizontal edge $e \in A \cap Y/K$, at depth $L$.  This lifts to a horizontal
edge $\tilde e$ in $X$ at depth $L$.  This can be pushed up to a path $\tilde p$
of length at most $2^L$ in $D^{-1}(0) \cap \tilde A$, which project to a path $p$ in $D^{-1}(0) \cap A$ above $e$.  We have already defined
$\eta_1(p)$, and this path lies above an edge in $\mc{H}_{L_2}$ at depth
$L$.  This edge is $\eta_1(e)$.  Thus we have defined $\eta_1$ on the
$1$-skeleton of $A \cap Y/K$.

We now define $\eta_2$ on the $1$-skeleton of $\mc{H}_{L_2}$.  Let
$e'$ be an edge at depth $L$ in $\mc{H}_{L_2}$.   There is a 
path $p'$ above $e'$ at the $0$-level of $\mc{H}_{L_2}$.  The
path $\eta_2(p')$ is already defined.  The path $\eta_2(p')$ lifts
to a path in $D^{-1}(0) \cap \tilde A$, and lies above an edge $\tilde e$
in $D^{-1}(L) \cap \tilde{A}$, which in turn projects to an edge $e$ in 
$A \cap Y/K$.  Set $\eta_2(e') = e$.

We have now defined $\eta_1$ and $\eta_2$ on the $1$-skeleta, and we leave it as an exercise to prove that they are mutually inverse.

The map $\eta_1$ obviously extends over the $2$-skeleta.
Consider then a $2$-cell $c$ in $\mc{H}_{L_2}$.  Then
$\eta_2(\partial c)$ is a loop in $A \cap Y/K$ of length at most $5$.
Suppose $\eta_2(\partial c)$ lifts to a path $\sigma$ which is not a loop.  
Let $k \in K_i$ be the element which sends one endpoint of $\sigma$ to
the other.  Then $|k|_{P_i \cap S} \le 5. 2^{L_2}$, a contradiction.
Therefore, $\eta_2(\partial c)$ {\em does} lift to a loop in $Y$, so there
is a $2$-cell filling $\eta_2(\partial c)$.  This $2$-cell is $\eta_2(c)$.
\end{proof}

The following proposition follows easily from Observations \ref{o:2step}
and \ref{o:ccomp} and Lemma \ref{l:shell}.

\begin{proposition}\label{ZisanX}
If $|K_i|_{P_i}\geq 12\cdot 2^{L_2}$ for all $i\in \{1,\ldots,n\}$,
then $Z$ is equivariantly isomorphic to $X(G/K,\mc{P}',{S},{\mc{R}})$.
\end{proposition}

\section{A linear isoperimetric inequality}  \label{section:lollipop}

We have reduced the proof of Theorem \ref{t:rhds} to proving
that the space $Z$ satisfies a homological isoperimetric inequality.
This is proved in Theorem \ref{t:Zhlii} below.  The proof of this
result breaks neatly into two pieces: a combinatorial piece 
(Proposition \ref{prop:keypoint}), and a
homological piece (which becomes Theorem \ref{t:Zhlii} below).  

We continue to assume that $G$ is a group which is hyperbolic
relative to $\mc{P}=\{P_1,\ldots,P_n\}$, and that $K_i\lhd P_i$
for each $i$. 
We let $X=X(G,\mc{P},S,\mc{R})$, $Y$, and $Z$ be as described in
section \ref{surgered}.  In order to slightly simplify the proof of 
Proposition \ref{prop:keypoint}, we replace Assumption \ref{K} with
the slightly stricter: 
\begin{assumption} \label{K'}
\[|K_i|_{P_i}\ge 24 \cdot 2^{L_2}=24\cdot 2^{3000\delta}\]
for each $i$.
\end{assumption}

We show (Theorem \ref{t:Zhlii})
that under these conditions, the surgered space $Z$
described in Section \ref{surgered} satisfies a linear homological
isoperimetric inequality.
It follows via Theorem \ref{t:tfae} and Proposition \ref{ZisanX}
that $G/K$ is hyperbolic
relative to the images of the subgroups in $\mc{P}$.

Let
\[\rho\co Z^{(1)}\to X/K^{(1)}\]
be the map from Lemma \ref{ZintoXmodK}.

\begin{lemma}\label{missing}
Let $w$ be a $2$-cell in $Z$ so that $\rho(\partial w)$ does not bound a
$2$-cell in $X/K$.  Then $w\subset D^{-1}[L_2+1,\infty)$.
\end{lemma}
\begin{proof}
Certainly $\rho(\partial w)$ must lie in $D^{-1}([L_2,\infty))$ in $X/K$, since $X/K$ and $Z$ are identical between depth $0$ and $L_2$. 

If $\rho(\partial w)$ does not surround a $2$-cell then it does not lift
to a loop in $X$.  Thus, since the length of $\rho(\partial w)$ is at
most $5$, it lifts to a path in $X$ of length $5$, whose endpoints are
in the same orbit under the action of $K$.  Theorem \ref{t:injects}
and Assumption \ref{K'} now imply that these points must lie in the same
orbit under the action of the stabilizer of the horoball in which they lie.
Thus there is some $k \in K_i$ so that $gkg^{-1}$ sends one endpoint
to the other for some $g$.

Let $L$ be the minimal depth of $\rho(\partial w)$.  There is a lift
of $\rho(\partial w)$ starting and ending at depth $L$ in $X$, and 
the endpoints of this lift are joined by a horizontal path of length at most
$4$ (as the diligent reader may readily verify).  Therefore, the points
at depth $0$ above these endpoints may be joined by a horizontal 
path of length at most $4 \cdot 2^L$.  The length of such a path is
an upper bound for $|k|_{P_i \cap S}$, so by  Assumption \ref{K'} we
have $4\cdot 2^L \geq 24 \cdot 2^{L_2}$.  In particular,
$L > L_2 +1$ as required.
\end{proof}
\begin{definition}
We refer to $2$-cells in $Z$ as described in the above lemma as
\emph{missing} $2$-cells.  Note that the map $\rho$ extends to those
$2$-cells of $Z$ which are not missing.
\end{definition}
\begin{definition}\label{d:pmp}
Suppose that $E$ is a 
cellulated disk, and that $\phi\co E\to Z$ is a
combinatorial map. A \emph{partly missing piece} of $\phi$ is a component
$C$ of $(D\circ\phi)^{-1}(L_2,\infty)$ so that $\phi(C)$ contains some
missing $2$-cell.  Let $E_*\subseteq E$ be the complement of the
partly missing pieces of $E$.

Let $\bar{E}$ be a closed regular neighborhood of $E_*\cup \partial E$
whose image under $\phi$ contains no missing $2$-cells (we may need
to adjust $\phi$ by a small homotopy rel $\partial E$ (and re-cellulate)
to ensure that $\bar{E}$ exists).
A \emph{reducing arc} for $\phi$ is an map $\sigma\co[0,1]\to E$
satisfying the following:
\begin{enumerate}
\item $\sigma$ is an essential arc in $\bar{E}$.
\item The endpoints of $\sigma$ lie in partly missing pieces of $E$.
\item $\sigma$ is homotopic rel endpoints into $D^{-1}(L_2,\infty)$,
  and this homotopy does not pass over any missing $2$-cell.
\end{enumerate}
\end{definition}
\begin{remark}
This is slightly different from the way reducing arcs were defined in
Section \ref{s:dyingwords} for a couple of reasons.
First, we need to deal with the
possibility that, for instance, two components of $E\smallsetminus
E_*$ have intersecting closures.  Second, we do not want to reduce
along arcs with an endpoint in $\partial E$, because this would change
the loop being filled.
\end{remark}

\begin{lemma}\label{l:simpconn}
Let $E$ be a disk.
Suppose that  $\gamma$ is a combinatorial loop in $Z$, and $\phi \co
E\to Z$ is a filling of $\gamma$.  Then there is some $\phi'\co E\to
Z$ so that each partly missing piece of $\phi'$ is simply connected,
and the number of partly missing pieces of $\phi'$ is no more than the
number of partly missing pieces of $\phi$.
\end{lemma}
\begin{proof}
Fix a partly missing piece $C$, and let $\xi$ be the component of the
boundary of $\bar{E}$ which lies in $C$ and contains all other
components of $\partial\bar{E}\cap C$.  The loop $\phi(\xi)$ lies
entirely in $Z\smallsetminus Y$, and can thus be filled there.  We
set $\phi'=\phi$ outside of $\xi$, and set $\phi'$ equal to this new
filling inside $\xi$.
\end{proof}

\begin{lemma}\label{l:pmpmin}
Let $E$ be a disk.
Suppose that  $\gamma$ is a combinatorial loop in $Z$, and $\phi \co
E\to Z$ is a
filling of $\gamma$ chosen so that the number of partly missing pieces
of $E$ is minimized.
Then there are no reducing arcs for $\phi$.
\end{lemma}
\begin{proof}
If there is a reducing arc $\sigma$, we may reduce the number of
partly missing pieces as follows:  

Suppose first that the arc joins distinct
partly missing pieces.  Then we may cut open $E$ along the arc and add
in two copies of the homotopy into $D^{-1}[L_2,\infty)$, thus
combining the two partly missing pieces into one.  Since the homotopy
passes over no missing $2$-cells, we have decreased the number of
partly missing pieces.

Next suppose that the arc $\sigma$ joins some partly missing piece $C$
to itself.  Again we can cut open along $\sigma$ and insert the
homotopy, thus creating a non-simply connected partly missing
pieces.  By Lemma \ref{l:simpconn}, this non-simply connected
component can be replaced with a simply connected one.  Since $\sigma$
was essential, it enclosed at least one partly missing piece other
than $C$, and so the modified map has fewer partly missing pieces than
it did before.
\end{proof}

\begin{definition} \label{d:diskextend}
If $\phi$ is a map of a disk $E$ into $Z$, then
$\rho\circ\phi|_{E_*\cup\partial E}$
can be extended to a proper map of a punctured disk
\[\extend{\phi}\co\extend{E_*}\to X/K.\]
The surface $\extend{E_*}$ can be obtained from $E$ by removing a
point from the interior of each partly missing piece.  The map
$\extend{\phi}$ is then defined to be equal to $\rho\circ\phi$ on
$E_*\cup\partial E$.  The complement of $E_*\cup\partial E$ in
$\extend{E_*}$ is a union of annuli;  the map $\extend{\phi}$ is
defined so that the image of these annuli consists entirely of
vertical squares, just as in Definition \ref{d:extend}.
\end{definition}

The proof of Lemma \ref{l:straighten} easily adapts to a proof
of the following:

\begin{lemma} \label{l:diskstraighten}
Let $E$ be a disk, and let $\phi\co E\to Z$ be a map with no reducing
arcs, so that $\phi$ has at least one partly missing piece and
$\phi(\partial E)\cap \Gamma_K$ is nonempty.  
Let $\mc{T}$ be a partially ideal triangulation of $\extend{E_*}$ with
a single vertex $v_0\in\partial E=\partial (\extend{E_*})$ satisfying
$\phi(v_0)\in \Gamma_K$.
Then $\extend{\phi}$ is properly homotopic to a map 
\[\straightedge{\extend{\phi}}\co \extend{E_*}\to X/K\]
which sends each edge of $\mc{T}$ to the image in $X/K$ 
of a preferred path in $X$.
\end{lemma}

The key combinatorial step in proving Theorem \ref{t:Zhlii}
is the following:
\begin{proposition}\label{prop:keypoint}
There is a constant $C = C(\delta) > 0$ so that the following holds:
Let $w\co S^1\to Z$ be a combinatorial loop, and suppose that 
\[\phi\co E \to Z\] is a filling of this loop by a disk $E$, and that
this filling has no reducing arcs.  If $\phi$ has $m\geq 1$ partly missing
pieces, then $|w|_1\geq C m$.
\end{proposition}
\begin{proof}
We will show that $C$ can be chosen equal to
$\min\{1,\frac{L_2}{2(\lambda+\epsilon)}\}$, where $\lambda$ and
$\epsilon$ are the constants of quasi-geodesicity from Corollary
\ref{c:quasigeodesic}.
If $w$ does not intersect $\Gamma_K=\Gamma(G,S)/K$, then $w$ can be
filled with a disk containing at most one partly missing piece.
We may therefore assume that $w(1)\in \Gamma_K \subset Z$.

Let $\gamma$ be a loop homotopic (rel $1$)
to the loop $\rho\circ w$ in $\Xgraph/K$, so that $\gamma$
lifts to a preferred path $\twid{\gamma}$
in $X$.  Let $g\in G$ be the unique group element which
sends the initial point of $\twid{\gamma}$ to the terminal point.

By Corollary \ref{c:quasigeodesic}, 
$|\gamma|_1=|\twid{\gamma}|_1\leq \lambda d(1,g)+\epsilon$, where
$d(1,g)$ is the
distance from $1$ to $g$ in $X$.  As $|w|_1$ is
bounded below by $d(1,g)$, we have
\[|w|_1\geq \frac{1}{\lambda+\epsilon} |\gamma|_1.\]
It thus suffices to bound $|\gamma|_1$ below linearly in terms of $m$.
We remark that $\gamma(1) = (\rho \circ w)(1) \in \Gamma_K \subset 
X/K$.

By assumption, we have a map $\phi \co E \to Z$ of a disk into $Z$
with no reducing arcs.  
Let $\extend{\phi} \co \extend{E_*} \to X/K$ be the proper map from
Definition \ref{d:diskextend} of an $m$-times
punctured disk into $X/K$. Note that $\extend{\phi}|_{\partial(\extend{E_*})} = \rho \circ w$.

Let $\mc{T}$ be a (partially ideal) triangulation
of $\extend{E_*}$ with one vertex the preimage of $(\rho \circ w)(1)$ and
all other vertices ideal, and one edge on the boundary.  This induces
an obvious triangulation $\mc{T}'$ of the disk $E$.

By Lemma \ref{l:diskstraighten}, $\extend{\phi}$ is homotopic to
$\straightedge{\extend{\phi}} \co \extend{E_*} \to X/K$ so that
$\straightedge{\extend{\phi}}|_{\partial \extend{E_*}} = \gamma$,
and all edges of $\mc{T}$ map to paths in $X/K$ which lift to
preferred paths in $X$.

We now consider the map
\[	\skelmap \co \skel(\straightedge{\extend{\phi}}) \to X/G \cup (\Hbound)/G	,	\]
as in Remark \ref{r:maponskel}.

Let $p$ be a puncture on $\extend{E_*}$.  The link $Lk(p)$ in 
$\skel(\straightedge{\extend{\phi}})$ is a loop and we have
\[	\skelmap|_{Lk(p)} \co S_1 \to D^{-1}([L_2,\infty)) \subset X/G.	\]
This loop represents a conjugacy class in $P_i$ for some $i$, and 
this class is contained in $K_i$, by Lemma \ref{l:link}.

A puncture $p$ is called {\em interior} if $Lk(p)$ is composed entirely
of ribs and ligaments, and {\em exterior} otherwise.  The puncture $p$
corresponds to a vertex of the triangulation $\mc{T}'$ of $E$, which
we also describe as {\em interior} or {\em exterior}.

See Figure \ref{fig:lollipop}.
\begin{figure}[htb!]
\begin{center}
\begin{picture}(0,0)%
\includegraphics{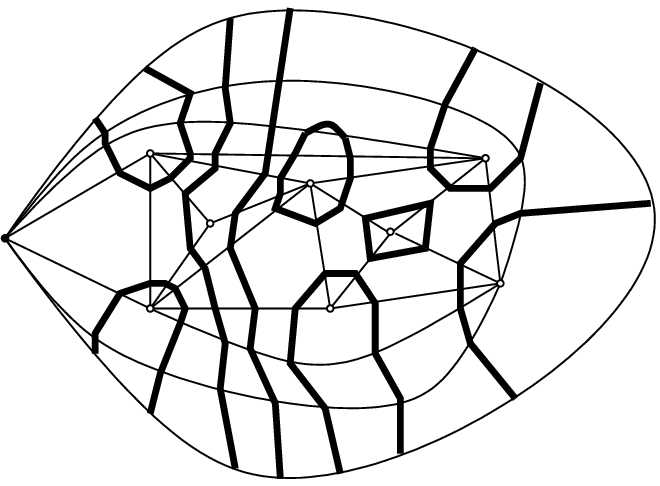}%
\end{picture}%
\setlength{\unitlength}{4144sp}%
\begingroup\makeatletter\ifx\SetFigFont\undefined%
\gdef\SetFigFont#1#2#3#4#5{%
  \reset@font\fontsize{#1}{#2pt}%
  \fontfamily{#3}\fontseries{#4}\fontshape{#5}%
  \selectfont}%
\fi\endgroup%
\begin{picture}(3008,2217)(3241,-2361)
\end{picture}%
\caption{This picture shows a possible (actually, impossible) picture of
  $\skel(\straightedge{\extend{\phi}})$, with two interior punctures,
  and the rest exterior.  Ribs and ligaments are shown in bold. 
  Each of the exterior punctures contributes
  at least $2$ points to $(D\circ\twid{\gamma})^{-1}(L_2)$.}
\label{fig:lollipop}
\end{center}
\end{figure}

Let $V_I$ be the number of interior vertices of $\mc{T}'$, and let
$V_\partial = m + 1 - V_I$ be the number of exterior vertices.  It is
clear from the definitions that the cardinality of 
$(D\circ\twid{\gamma})^{-1}(L_2)$ is at least $2 (V_\partial -1)$.
The set $D^{-1}(L_2)$ partitions $\twid{\gamma}$ into subsegments,
which alternate between lying in $D^{-1}[0,L_2]$ and lying in
$D^{-1}[L_2,\infty)$.  
  As
the initial and terminal points of $\twid{\gamma}$ lie in $D^{-1}(0)$,
those subsegments of $\twid{\gamma}$ with image in $D^{-1}[0,L_2]$ must
have length at least $L_2$  (In fact all but the initial and terminal
subsegments must have length at least $2 L_2$, as they must pass
between distinct $L_2$-horoballs.).  There are at least
$V_\partial$ such subsegments, and so 
\begin{equation} \label{e:sizeofgamma}
|\gamma|_1=|\twid{\gamma}|_1\geq V_\partial L_2 = (m+1-V_I)L_2 > (m-V_I)L_2 .
\end{equation}
In order to bound $|\gamma|_1$ below by a linear function of $m$, it
therefore suffices to show that $V_I$ is at most some definite
proportion of $m$ (bounded away from $1$).  

Let $R$ be the number of ribs in $\skel(\straightedge{\extend{\phi}})$.
Then $R \le 6(2m-1)$, by Corollary \ref{c:ribsbound}.  However,
$R \ge 24 V_I$, by Assumption \ref{K'} and Lemma
\ref{l:link}.

This implies that $V_I < \frac{1}{2} m$.  By \eqref{e:sizeofgamma}, we deduce that
\[	|\gamma|_1 \ge \frac{L_2}{2}m,	\]
and hence
\[	|w|_1\ge \frac{L_2}{2(\lambda + \epsilon)} m.	\]
\end{proof}

\subsection{Proof of Theorem \ref{t:rhds}}

We make the same assumptions about $G$, $\mc{P}$, $S$, $\mc{R}$ and
$X$ stated at the beginning of the last section.

\begin{theorem} \label{t:Zhlii}
If $\min_j \{ |K_j|_{P_j} \} > 24 \cdot 2^{L_2}= 24 \cdot
2^{3000\delta}$, then $Z$ satisfies a 
homological linear isoperimetric inequality in the sense of Definition
\ref{d:hlii}. 
\end{theorem}
\begin{proof}
The idea behind the proof is as follows:  Take a loop in $Z$, which
is filled by {\em some} disk.  Attempt to move this disk to $X/K$, fail,
and find a punctured disk mapping into $X/K$.  Triangulate, lift and straighten the triangles, project back to $X/K$, attempt to transfer
back to $Z$, fail and fill the failures.

By Theorem \ref{t:coherent}, it suffices to show that there is a
constant $M>0$ so that any simple loop $c$ bounds a rational 
$2$-chain $w$ with
\[ |w|_1\leq M |c|_1. \]

Let $c$ be a simple (combinatorial) loop in $Z$.  If $c$ lies entirely
inside a single horoball, we may fill with a disk of area at most $3
|c|_1$ by Proposition \ref{p:combhorolii}.  We may thus suppose that
$c(1)\in D^{-1}(0)\subset Z$.  Thus $\rho\circ c$ lifts to a path
$\twid{c}$ in $X$ between two (not necessarily distinct) 
vertices $g$ and $h$ of the Cayley graph $\Gamma(G,S)\subset X$.  
Consider the $1$-chain $\underline{\twid{c}}$ corresponding to
the path $\twid{c}$.  Then $c - q_{g,h}$ is a $1$-cycle.  By
Proposition \ref{p:qisqg}, $|q_{g,h}|_1 \le \qConst d(g,h) \le \qConst |c|_1$.
Hence, if $M_X$ is the homological isoperimetric constant for $X$,
there is a $2$-chain $\twid{\omega}$ so that
$\partial\twid{\omega} = \underline{\twid{c}} - q_{g,h}$, and
$|\twid{\omega}|_1 \le M_X \big( \qConst + 1 \big) |c|_1$.

As $Z$ is simply
connected, there is a combinatorial map of a cellulated disk $E$
\[\phi\co E\to Z\]
so that $\phi|_{\partial E}$ is $c$.  We may suppose that the map
$\phi$ has the minimal possible number of partly missing pieces in the
sense of Definition \ref{d:pmp}.  By Lemma \ref{l:pmpmin}, $\phi$ has
no reducing arcs.  By Proposition \ref{prop:keypoint}, the number of
partly missing pieces is at most $C^{-1} |c|$, for a $C>0$ which
depends only on $\delta$.  Thus we can triangulate the punctured
surface $\extend{E_*}$ defined in Definition \ref{d:diskextend} with
a triangulation $\mc{T}$ consisting of
fewer than $2C^{-1} |c|$ triangles. Let $\straightedge{\extend{\phi}}|_e$ 
be the map from Lemma \ref{l:diskstraighten}. 
This induces, for each triangle $T \in \mc{T}$, a preferred
triangle $\psi_T\co \partial\Delta\to X\cup\Hbound$.
  
Suppose that $T \in \mathcal{T}$, and that the image of the vertices of $T$ are $a$, $b$, and $c$ (with order coming from the orientation of $T$).  Let $c_T = c_{abc}$,
the $1$-chain defined in Definition \ref{d:c_abc}, and let
$\omega_T = \omega_{abc}$ be the $2$-chain as in Corollary 
\ref{c:Xfillings}.
Let $\twid{\xi} = \sum_{T\in \mc{T}}\omega_T$.
Each $\omega_T$ satisfies $|\omega_T|_1\leq M_X T_1$ by Corollary
\ref{c:Xfillings}.

Now let $\mu = \pi_\sharp(\twid{\omega}+\twid{\xi})$, and
let $\mu^{\mathrm{thick}}$ be the $2$-chain which comes from including
only those $2$-cells which lie entirely in $Y=D^{-1}[0,L_2]$.  Note that
\begin{eqnarray*}
|\mu^{\mathrm{thick}}|_1 & \leq & |\mu|_1\\
&  \leq & |\twid{\omega}|_1+|\twid{\xi}|_1\\
& \leq & \big[ (M_X(\qConst +1))+2C^{-1}T_1M_X \big] |c|_1
\end{eqnarray*}
Let $M_{\mu}=(M_X(\qConst +1)+2C^{-1}T_1M_X$.

Since $\mu^{\mathrm{thick}}$ is supported entirely in $Y$, it
determines a $2$-chain $\mu_Z$ in $Z$.  Furthermore,
$c^{\mathrm{thin}}:=c-\partial \mu_Z$ satisfies the following:
\begin{enumerate}
\item The support of $c^{\mathrm{thin}}$ lies entirely in
  $D^{-1}[L_2,\infty)\subset Z$, and
\item $|c^{\mathrm{thin}}|_1\leq |c|_1+|\partial\mu_Z|\leq
  (1+MM_\mu)|c|_1$, where $M$ is the maximum length of the boundary of
  a $2$-cell in $Z$ (which is the same as that maximum length in $X$).
\end{enumerate}
Thus by Proposition \ref{p:combhorolii} and Theorem \ref{t:coherent}
there is a $2$-chain $\zeta$ satisfying 
\begin{enumerate}
\item $\partial \zeta = c^{\mathrm{thin}}$, and
\item $|\zeta|_1\leq 3 |c^{\mathrm{thin}}|_1\leq 3(1+MM_\mu)|c|_1$.
\end{enumerate}
Finally, we note that $\partial(\mu_Z+\zeta)=c$ and
$|\mu_Z+\zeta|_1\leq |\mu_Z|_1+|\zeta|_1\leq (M_\mu+3+3MM_\mu)|c|_1$.
This completes the proof of Theorem \ref{t:Zhlii}, and hence also of
Theorem \ref{t:rhds}.
\end{proof}

We close by proving that $G/K$ is nonelementary.

\begin{theorem} \label{t:nonelementary}
Suppose that $G$ is torsion-free and that $|K_i|_{P_i} > 24\cdot 2^{L_2}$ for each $i$.
Then $G/K$ is non-elementary relatively hyperbolic (relative to
$\{ \iota_i(P_i/K_i) \}_i$. 
\end{theorem}
\begin{proof}
If all parabolics are finite, then $G/K$ is hyperbolic.  

By Theorem \ref{t:onethickZ}, there is a hyperbolic element
$g \in G$ with an axis which is entirely contained in 
$D^{-1}[0,19\delta]$.  If $x \in X$ is contained in this axis, then for all $j$, the preferred path $\pp_{x,g^j.x}$ is the geodesic between these points (since $19\delta + \delta \le L_1$).  Therefore, by Proposition
\ref{p:thickdontdie}, none of these elements die in $G/K$.  Therefore,
$G/K$ is infinite.  Furthermore, by Proposition \ref{p:thicknotparabolic},
$g$ does not project to a parabolic element of $G/K$.  Thus $G/K$
is not equal to any of the $\iota_i(P_i/K_i)$.  By Theorem \ref{t:dontintersect}, the intersection of two distinct parabolic subgroups 
of $G/K$ is trivial. 

Suppose $G/K$ is virtually cyclic.  Then $G/K$ has a finite normal subgroup $N$ with quotient either infinite cyclic or infinite dihedral.  In fact, $N$ must be contained in every parabolic; and so it is trivial.  However, the parabolic subgroups of $G/K$ have size (much)
greater than $2$.

We may now suppose that some $\iota_i(P_i/K_i)$ is infinite.  We have
already seen that $G/K$ does not equal $\iota_i(P_i/K_i)$ for any
$i$, so Theorem \ref{t:independentZ} implies that $G/K$ is
non-elementary relatively hyperbolic.
\end{proof}

\end{document}